\documentclass[mathpazo]{cicp}

\usepackage{subfig}
\usepackage{multirow}
\allowdisplaybreaks[4]

\begin{document}
\title{Finding similarity of orbits between two discrete dynamical systems via optimal principle}


 \author[Y. Chen and Y. Li]{Yuting Chen\affil{1} and
       Yong Li\affil{1,2}\comma\corrauth}
 \address{\affilnum{1}\ College of Mathematics,
          Jilin University,
          Changchun 130012, P.R. China. \\
           \affilnum{2}\ School of Mathematics and Statistics,
           and Center for Mathematics and Interdisciplinary Sciences,
           Northeast Normal University, Changchun 130024, P.R. China.}
 \emails{{\tt cyt1012@jlu.edu.cn} (Y.~Chen), {\tt liyong@jlu.edu.cn} (Y.~Li)}

\begin{abstract}
Whether there is similarity between two physical processes in the movement of objects and the complexity of behavior is an essential problem in science. How to seek similarity through the adoption of quantitative and qualitative research techniques still remains an urgent challenge we face. To this end, the concepts of similarity transformation matrix and similarity degree are innovatively introduced to describe similarity of orbits between two complicated discrete dynamical systems that seem to be irrelevant. Furthermore, we present a general optimal principle, giving a strict characterization from the perspective of dynamical systems combined with optimization theory. For well-known examples of chaotic dynamical systems, such as Lorenz attractor, Chua's circuit, R$\rm\ddot{o}$ssler attractor, Chen attractor,  L$\rm\ddot{u}$ attractor and hybrid system, with using of the homotopy idea, some numerical simulation results demonstrate that similarity can be found in rich characteristics and complex behaviors of chaotic dynamics via the optimal principle we presented.
\end{abstract}


\maketitle

\section{Introduction}
\label{sec1}
Discrete dynamical systems described by iteration of mappings appear everywhere, showing directive laws from physical science or result from simulations to better understand differential equations numerically. Generally, it is much more difficult but interesting to investigate how complex behavior happens to discrete dynamical systems than continuous dynamical systems after some iterations, since there are probably greater covered ranges and more ghost phenomena. Along with the development of computer technology, modeling problems by means of discrete dynamical systems mathematically has already been gained in different fields such as biology, economics, demography, engineering, and so on. It is universally acknowledge that no matter how different the various technologies develop as well as the objects appear in the research process, there are certain underlying similarities.

Similarity, in addition to being frequently encountered, is viewed as a fundamental concept in scientific research. The idea of similarity has gained widespread popularity in the era of big data and machine learning by various means. For instance, scale similarity is found in many natural phenomena in the universe \cite{Wang2017(1)}. An embedding-based vehicle method with deep representation learning drastically accelerates trajectory similarity computation \cite{Chen2022}. A novel brain electroencephalography (EEG) clustering algorithm not only handles the problem of unlabeled EEG, but also avoids the time-consuming task of manually marking the EEG \cite{Dai2022}. Based on the cosine similarity, a transductive long short-term memory model is developed for temperature forecasting \cite{Karevan2020}. Self-similar coordinates are investigated in Lattice Boltzmann equation, showing that the time averaged statistics for velocity and vorticity express self-similarity at low Reynolds \cite{Zarghami2012}. Many other applications include gene expression \cite{Arbela2018}, image registration \cite{Mang2017}, web pages and scientific literature \cite{Wang2017}, fuzzy linguistic term sets \cite{Liao2014}, collaborative filtering \cite{Liu2014}, pattern analysis \cite{Zhao2013} and  preferential attachment \cite{Papadopoulos2012}. Indeed, the ubiquitous similarity is attributed to facilitate prediction of indeterminate events by analyzing known data, being an essential task in many natural systems and phenomena of real life.

A core part of similarity search is the so-called similarity measure whose famous characteristic is able to assess how similar two sequences are, in other words, the degree to which a given sequence resembles another. Many researchers have paid great attention to devise a proper similarity measure and have achieved several valuable results, which can be roughly categorized into two sorts. One sort is based on the traditional measures, such as Euclidean distance, dynamic time warping, cosine and cotangent similarity measures and Pearson correlation coefficient \cite{Liao2005}. The other sort is some transform-based methods, such as singular value decomposition, principal component analysis, Fourier coefficients, auto-regression and moving average model \cite{Bartolini2005,Fu2011}. The cautious selection of similarity measure scheme has long been a research hotspot, affecting the accuracy of further data mining tasks directly, such as classification, clustering and indexing \cite{Shen2018,Torrente2021}.

Up to now, whether there is similarity between two physical processes and how to seek similarity through a mathematical principle are still remain a significant challenge. Several theoretical approaches are available to deal with this problem by taking into account asymptotic equivalence, synchronization and stability just as some kind of similarity. Furthermore, almost all similarity measure criteria are exploited according to application background and actual data, which can only be regarded as quantitative representations to estimate pairwise similarity of a given series resembles another under certain conditions.

Determining similarity between orbits derived from chaotic systems generally characterized by highly complex behavior is particularly difficult when the general similarity measure is employed. Having in view that, for any two chaotic dynamical  systems, what are their similarities and how do we find them are both fundamental but challenging subjects in science and engineering. For this purpose, we will try to touch these problems.

To the authors best knowledge, this is the first work to develop a mathematical framework, studying the connection of orbits between two discrete dynamical systems from the perspective of dynamical systems combined with optimization theory. Novel contributions and results of this paper include 1) proposing the concepts of similarity transformation matrix and similarity degree to describe what extent the orbits derived from two discrete dynamical systems are similar; 2) presenting a general optimality principle by employing variational method when orbits of two discrete dynamical systems are similar at some step; 3) constructing hybrid dynamical system with richer and more complex dynamical behavior via the idea of homotopy, applying to numerical simulation.

The remainder of this paper begins with review of several typical chaotic attractors related to this study. We formalize the similarity via optimization techniques and give the definitions of similarity transformation matrix and similarity degree mathematically to assess what extent the orbits of two dynamical systems are similar, followed by establishing the main results of this paper in Section \ref{sec3}. Section \ref{sec4} reports numerical simulation results of chaotic systems to support the theoretical findings. Conclusions are drawn in Section \ref{sec5}.

\section{Chaotic systems}
\label{sec2}
Among a broad variety of dynamical systems in the universe, we consider some typical chaotic systems such as Lorenz attractor, Chua' circuit, R$\rm\ddot{o}$ssler attractor, Chen attractor, L$\rm\ddot{u}$ attractor and their hybrid systems.

\begin{figure}[htbp]
\centering
\begin{minipage}{0.45\linewidth}
\includegraphics[scale=.42]{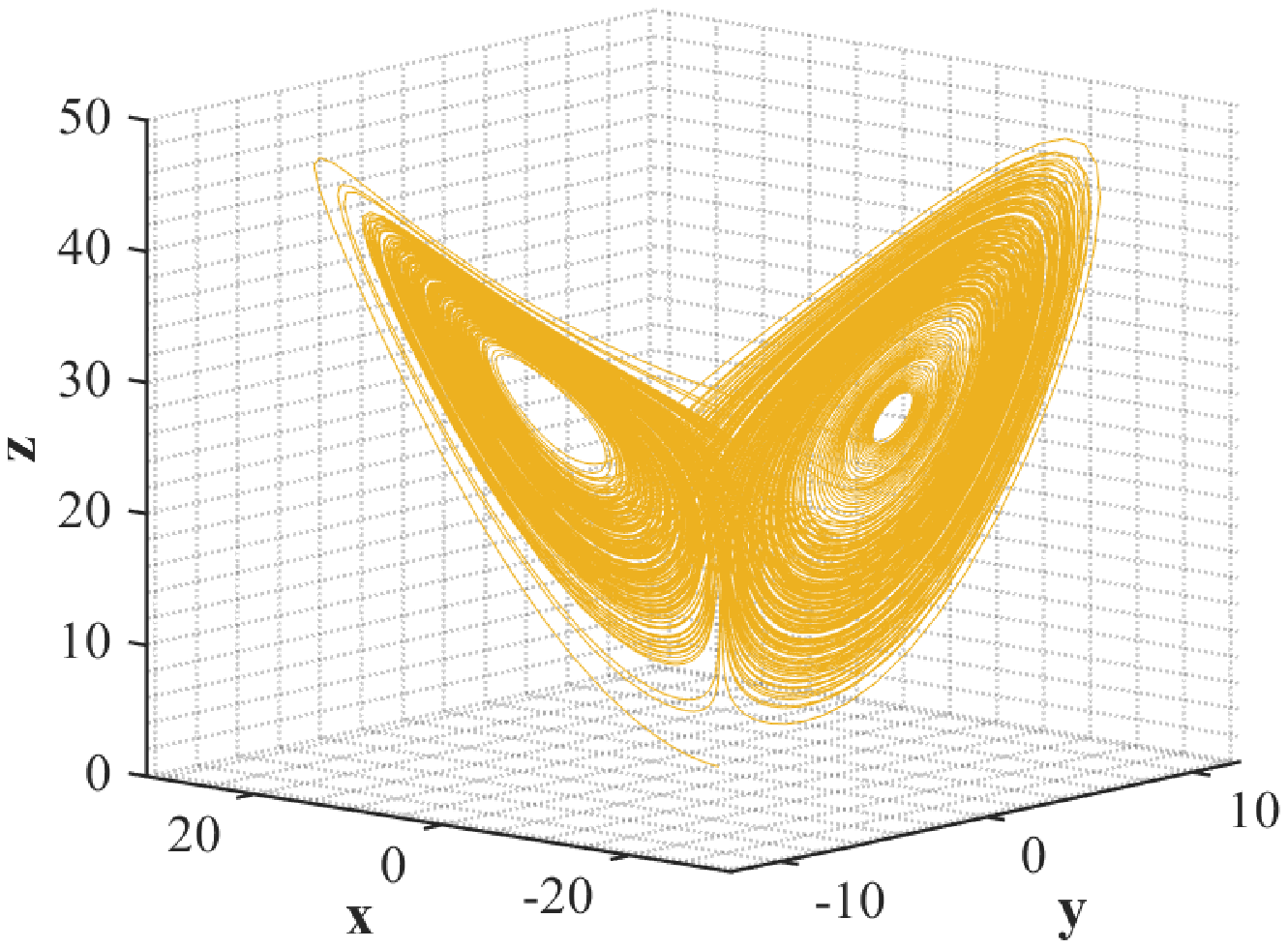}
\caption{Lorenz attractor.}
\label{Fig.1}
\end{minipage}
\begin{minipage}{0.45\linewidth}
\includegraphics[scale=.42]{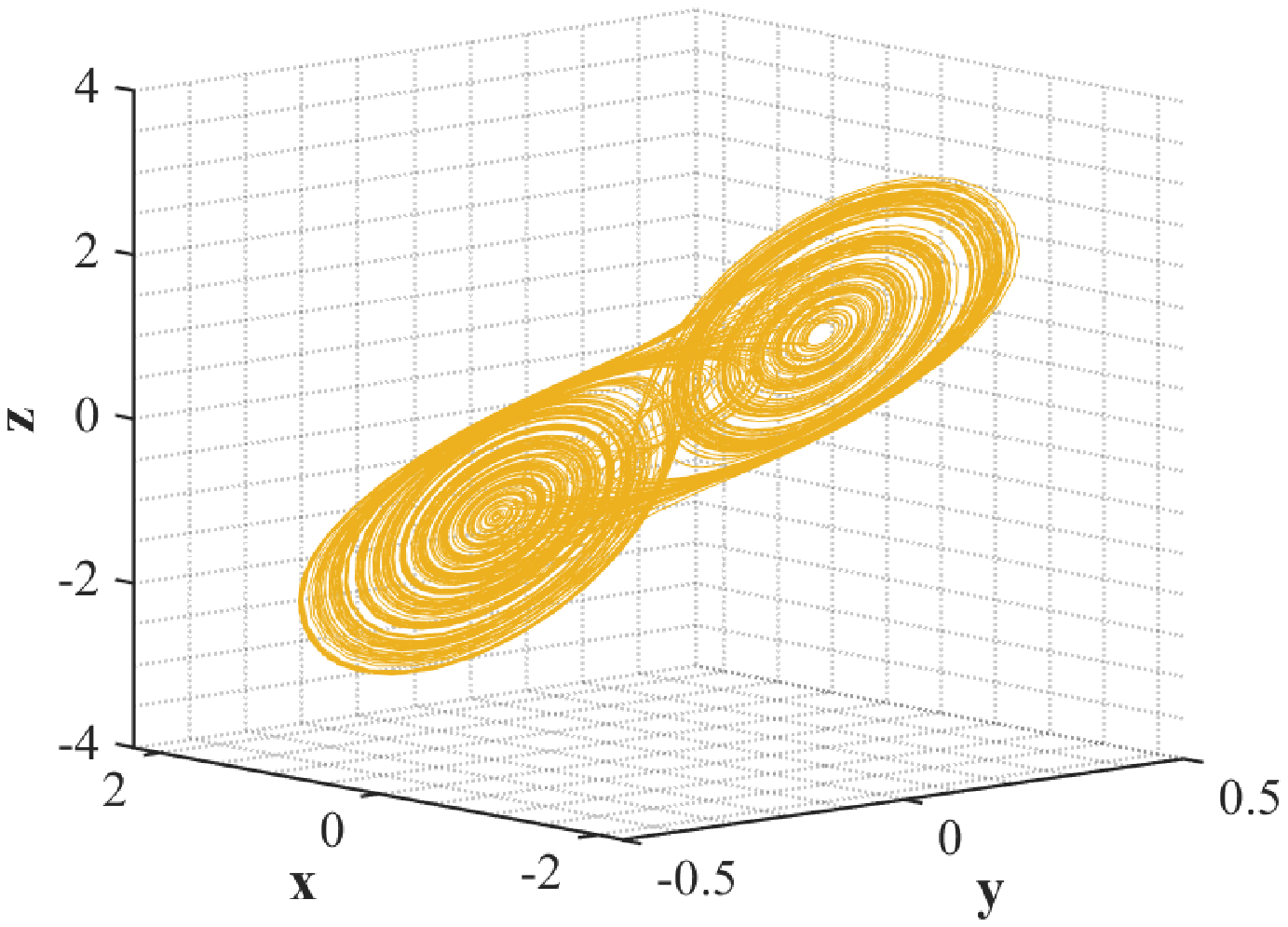}
\caption{Chua's circuit.}
\label{Fig.2}
\end{minipage}
\end{figure}

Lorenz attractor, the first chaotic dynamical system, was obtained by Lorenz in 1963 from simplified mathematical model developed for atmospheric convection while modelling meteorological phenomena \cite{Lorenz1963}. The chaotic system is a typical nonlinear system with three differential equations known as the Lorenz equations
\begin{equation}\label{1}
\dot{x}=-\sigma x+\sigma y,~~
\dot{y}=-xz+rx-y,~~
\dot{z}=xy-bz,
\end{equation}
where $x$, $y$, $z$ represent the system states and the system parameters are selected as $\sigma=10$, $b=8/3$, $r=28$. The initial conditions are $x(0)=0.1$, $y(0)=0.1$, $z(0)=0.1$, then the behavior of Lorenz attractor resembling a butterfly or figure eight is illustrated in Fig. \ref{Fig.1}.

Chua's circuit is the simplest electronic circuit known as nonperiodic oscillator \cite{Chua1992}. It has been confirmed by numerous experimental simulations and rigorous mathematical analysis that this circuit is able to produce an oscillating waveform exhibiting classic chaos behavior and many well-known bifurcation phenomena. Three ordinary differential equations are found as below in the analysis of Chua's circuit
\begin{equation}\label{2}
\dot{x}=\alpha[y-x-f(x)],~~
\dot{y}=x-y+z,~~
\dot{z}=-\beta y,
\end{equation}
where $x$, $y$ denote the voltage of capacities, $z$ represents inductance current and the parameters $\alpha$, $\beta$ are determined by the particular values of the circuit components. The function $f(x)$ is defined as a piece-linear function, describing the electrical response of nonlinear resistor
\begin{equation}\label{3}
f(x)=m_1x+\dfrac{1}{2}(m_0-m_1)(|x+1|-|x-1|).
\end{equation}
Fig. \ref{Fig.2} shows the double scroll attractor from Chua's circuit, in which the initial states are $x(0)=0.1$, $y(0)=0.1$, $z(0)=0.1$, and the parameters are determined as $\alpha=10$, $\beta=15$, $m_0=-1.2$, $m_1=-0.6$.

\begin{figure}[htbp]
\centering
\begin{minipage}{0.45\linewidth}
\includegraphics[scale=.42]{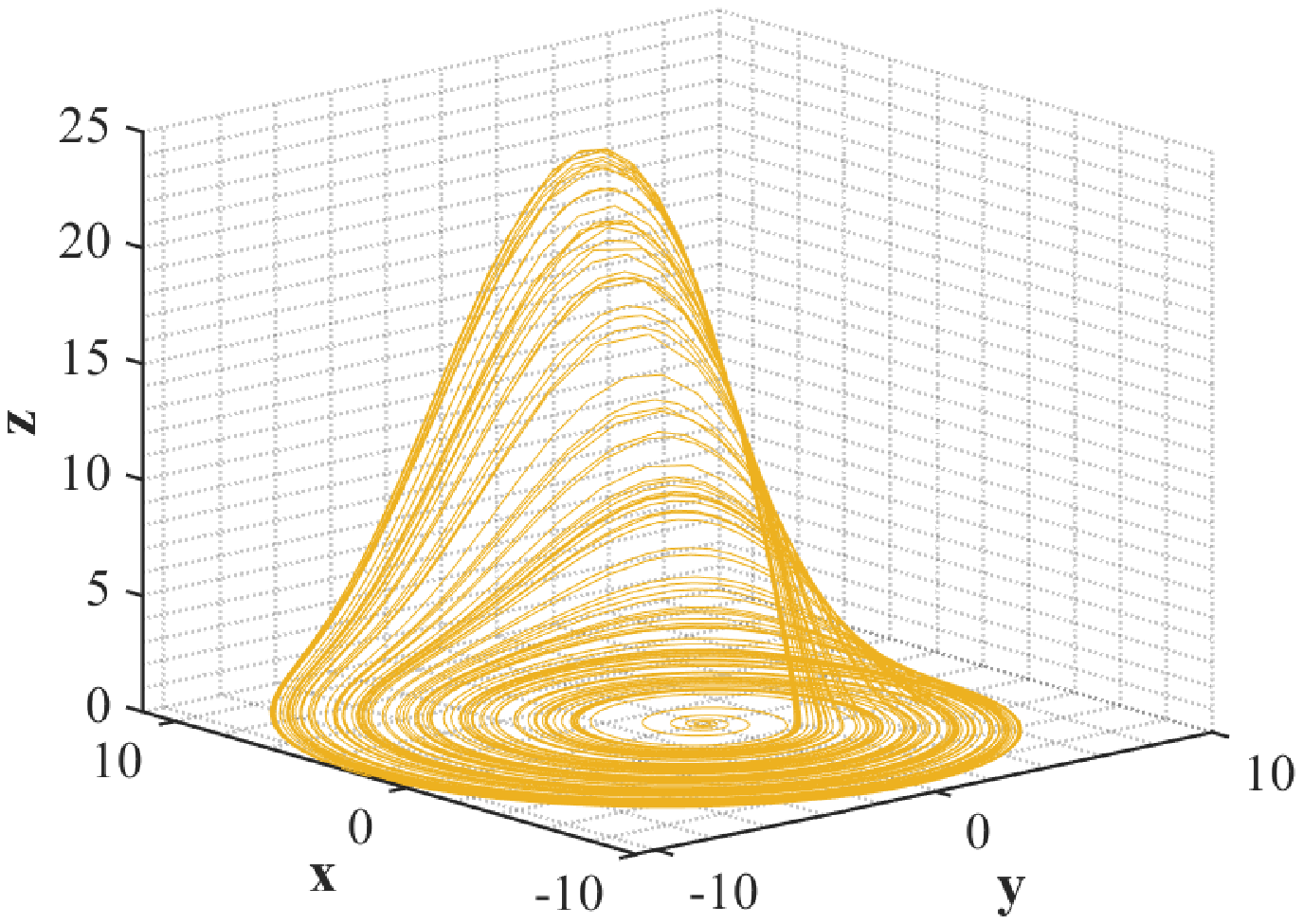}
\caption{R$\rm\ddot{o}$ssler attractor.}
\label{Fig.3}
\end{minipage}
\begin{minipage}{0.45\linewidth}
\includegraphics[scale=.42]{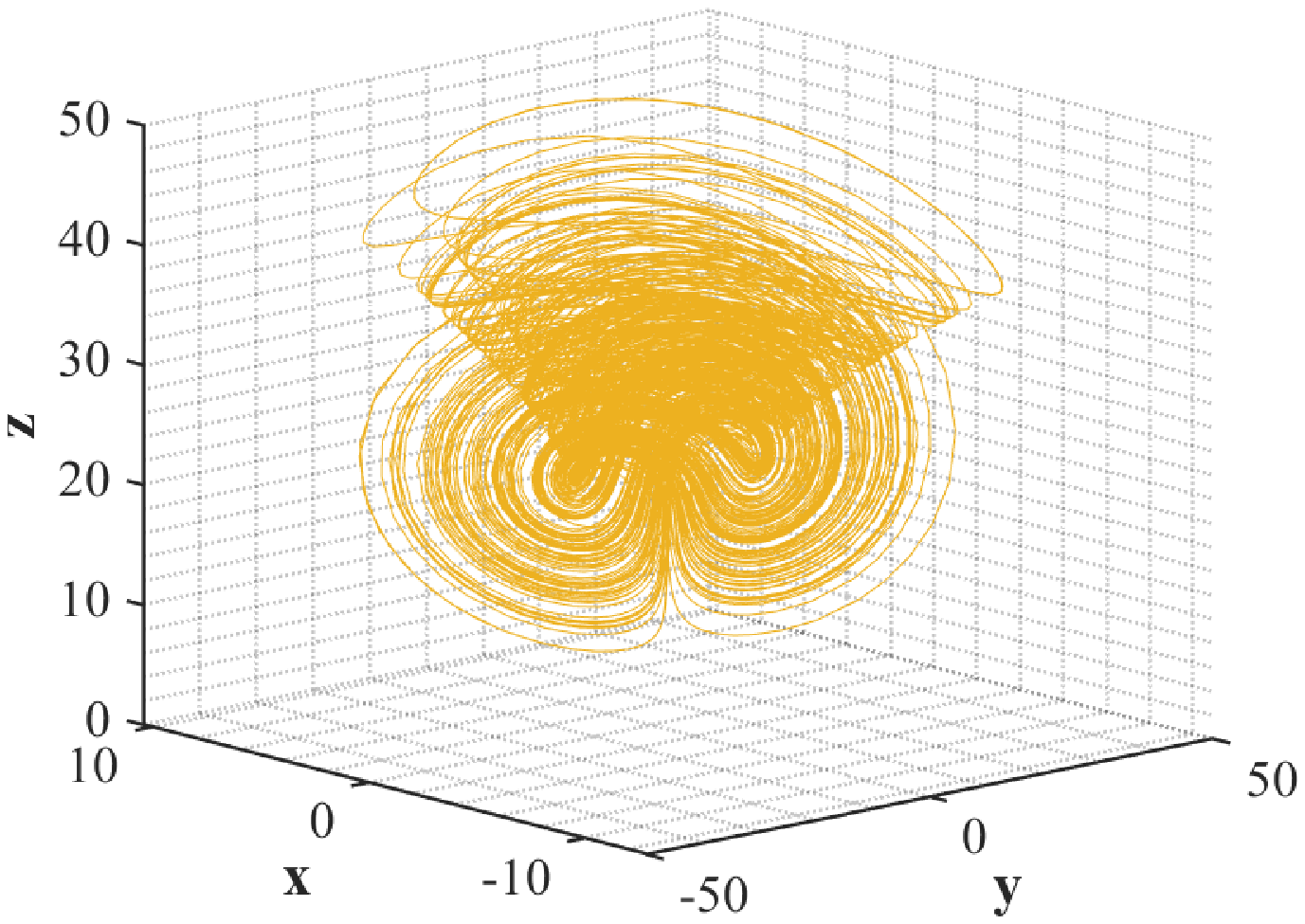}
\caption{Chen attractor.}
\label{Fig.4}
\end{minipage}
\end{figure}

R$\rm\ddot{o}$ssler attractor behaves similar to Lorenz attractor, and it is the most simple chaotic attractor from the topological point of view. This attractor is applied to modelling equilibrium in chemical reactions which is a chaotic solution to the system of three differential equations
\begin{equation}\label{4}
\dot{x}=-y-z,~~
\dot{y}=x+ay,~~
\dot{z}=b+z(x-c),
\end{equation}
where $x$, $y$, $z$ denote the system states and three parameters $a$, $b$ and $c$ are assumed to be positive \cite{Rossler1976}. We select numerical values of parameters as $a=0.2$, $b=0.2$, $c=5.7$ and give a typical orbit of R$\rm\ddot{o}$ssler attractor, which admits chaotic behavior, as shown in Fig. \ref{Fig.3}.

Chen attractor is found in the pursuit of chaotification, being similar but topologically not equivalent to Lorenz attractor \cite{Chen1999}. Despite Chen attractor with simple structure is the dual to Lorenz system, it is considered displaying even more sophisticated dynamical behaviors \cite{Ueta2000}. The three-dimensional autonomous system of ordinary differential equations with quadratic nonlinearities that describe Chen system are
\begin{equation}\label{5}
\dot{x}=a(y-x),~~
\dot{y}=(c-a)x-xz+cy, ~~
\dot{z}=xy-bz,
\end{equation}
where $x$, $y$, $z$ are the system states and $a$, $b$, $c$ are real parameters. For parameters values $a=40$, $b=3$, $c=28$, we obtain a Lorenz-based wing attractor as shown in Fig. \ref{Fig.4}.

\begin{figure}[htbp]
\centering
\subfloat[$u=0$.]{\includegraphics[scale=.3]{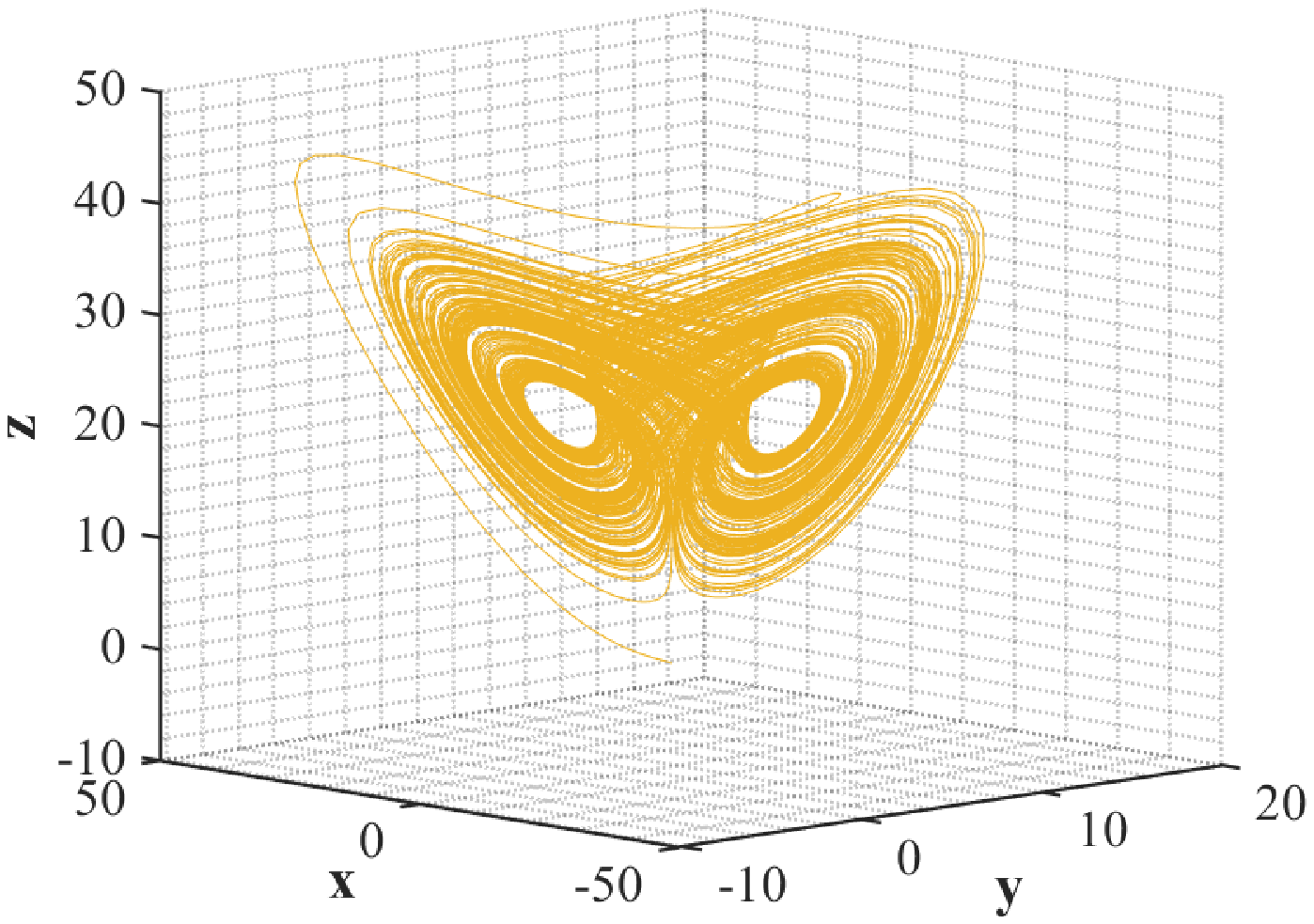}
\label{Fig.5a}}
\hfil
\subfloat[$u=-12$.]{\includegraphics[scale=.3]{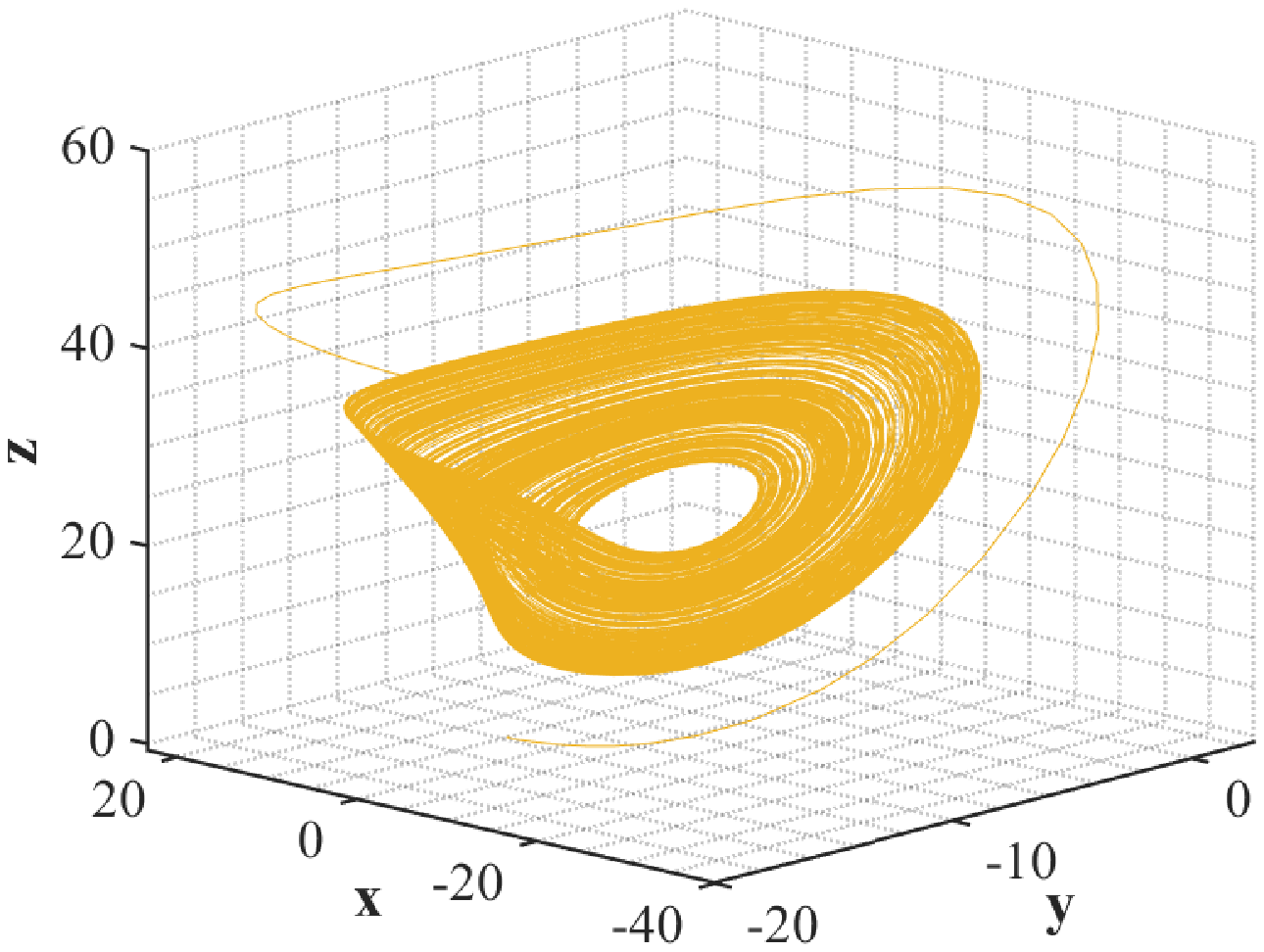}
\label{Fig.5b}}
\hfil
\subfloat[$u=12$.]{\includegraphics[scale=.3]{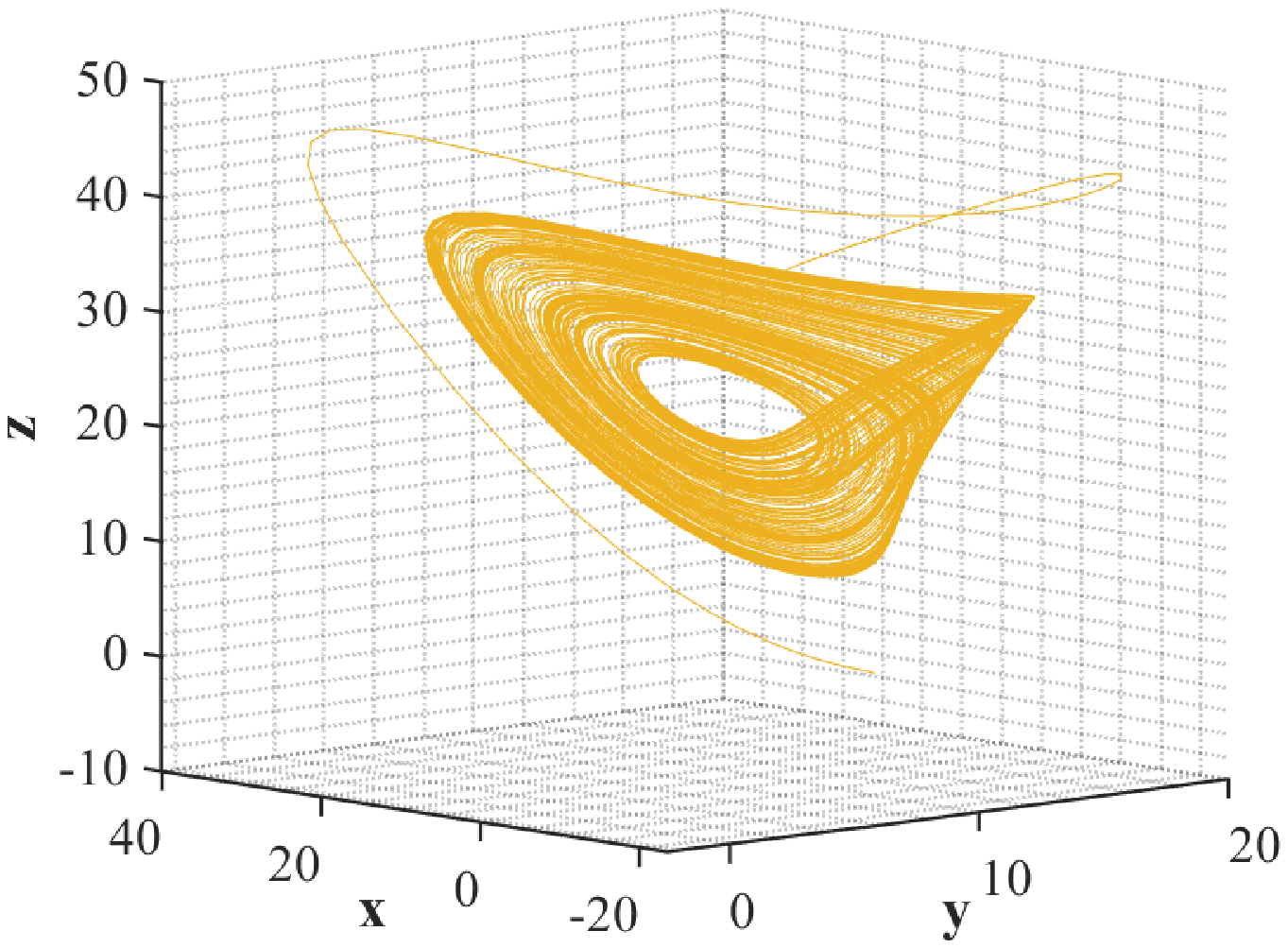}
\label{Fig.5c}}
\caption{L$\rm\ddot{u}$ attractor.}
\label{Fig.5}
\end{figure}

L$\rm\ddot{u}$ attractor is another example that captures the paradigms of chaotic system, which connects Lorenz attractor and Chen attractor and represents the transition from one to the other \cite{Lu2002(1),Lu2002(2)}. In order to reveal the topological structure of this chaotic attractor, consider its controlled system which is obtained by adding a constant $u$ to the second equation of L$\rm\ddot{u}$ system
\begin{equation}\label{6}
\dot{x}=a(y-x),~~
\dot{y}=-xz+cy+u,~~
\dot{z}=xy-bz,
\end{equation}
where $x$, $y$, $z$ denote the system states, $a$, $b$, $c$ are the system parameters. By varying the parameter $u$ considered as ``controller" of the controlled system, one can observe different dynamical behaviors, contributing to a better understanding of all similar and closely related chaotic system \cite{Lu2002(3)}. When $a=36$, $b=3$, $c=20$, all the simulation figures are summarized in Fig. \ref{Fig.5}.

The high sensitivity of chaotic systems to small perturbations of the initial states, together with the complex dynamical behavior characterized by rapidly changing solutions, make the research on chaotic dynamical systems challenging. The purpose of this paper focuses on finding similarity between the orbits of chaotic attractors via the general optimality principle, which will be discussed in next section.

\section{Main results}
\label{sec3}
We are now in the position to demonstrate that similarity of orbits derived by discrete dynamical systems can be found through a strict mathematical principle. This allows us to better understanding the motion trajectory and predicting process trend when looking at the rich behavior of complex physical processes.

\subsection{Simple Dynamical Systems}

Consider the following two discrete dynamical system
\begin{equation}\label{7}
x_{k+1}=f(k,x_k),
\end{equation}
\begin{equation} \label{8}
y_{k+1}=g(k,y_k),
\end{equation}
where the mappings $f,g:\mathbb{N_+}\times\mathbb{R}^n\rightarrow\mathbb{R}^n$ are of $\mathbb{C}^1$.
Starting from initial states $x_0$ and $y_0$, the solutions of (\ref{7}) and (\ref{8}) are denoted as
\begin{equation}\label{9}
x_k=x(k,x_0),
\end{equation}
\begin{equation}\label{10}
y_k=y(k,y_0),
\end{equation}
respectively. We introduce a new concept of similarity transformation matrix to deal with the problem of drawing a relation of similarity between (\ref{7}) and (\ref{8}).

\begin{definition}\label{def1}
Let $[n]$ denote the set $\{1,2,\ldots,n\}$. There exists an $n-$order matrix $A=(a_{ij})\in\Omega$ satisfies:

$\bullet$ $a_{ij}$ denotes some matrix element, where $i\in [n]$ and $j\in[n]$ are the $i$th row and $j$th column of $A$, respectively;

$\bullet$ $\Omega$ is a bounded closed convex set of~~$\mathbb{M}^{n\times n}$ whose interior $\Omega^\circ\neq\emptyset$;

$\bullet$ $y_0=Ax_0$.

We say that $A$ is a similarity transformation matrix if solutions between two discrete dynamical systems (\ref{7}) and (\ref{8}) satisfy $y_k=Ax_k$ at $k$th step. Otherwise, we say that they are not similar if such A does not exist.
\end{definition}

From the definition given above, it follows that the way to estimate similarity transformation matrix $A$ is to minimize the cost functional
\begin{align}\label{11}
J(A)=\min\limits_A\sum\limits_{k=0}^N\|Ax_k-y_k\|_2^2.
\end{align}
Motivated by first-order optimality condition which is the foundation for many of optimization algorithms, we know that the optimal solution of (\ref{11}) is equivalent to
\begin{align}\label{12}
\frac{\partial J(A)}{\partial a_{ij}}
=\frac{\partial}{\partial a_{ij}}(x_k^{\rm T}A^{\rm T}Ax_k)
 -\frac{\partial}{\partial a_{ij}}(2x_k^{\rm T}A^{\rm T}y_k)
 +\frac{\partial}{\partial a_{ij}}(y_k^{\rm T}A^{\rm T}yx_k)
=0.
\end{align}

For the first term of the right-hand side in (\ref{12}), we have
\begin{align}\label{13}
\dfrac{\partial}{\partial a_{ij}}(x_k^{\rm T}A^{\rm T}Ax_k)
&=\dfrac{\partial}{\partial a_{ij}}(x_k^{\rm T}A^{\rm T})Ax_k
    +x_k^{\rm T}A^{\rm T}\dfrac{\partial}{\partial a_{ij}}(Ax_k) \nonumber\\
&=x_k^{\rm T}\left(\dfrac{\partial}{\partial a_{ij}}A^{\rm T}\right)Ax_k
    +x_k^{\rm T}A^{\rm T}\left(\dfrac{\partial}{\partial a_{ij}}A\right)x_k \nonumber\\
&=(0,\ldots,0,x_{kj},0,\ldots,0)Ax_k+x_k^{\rm T}A^{\rm T}\left(
        \begin{array}{lllllll}
             ~0  \\
             ~\vdots\\
             ~0 \\
             x_{kj}\\
             ~0 \\
             ~\vdots \\
             ~0 \\
        \end{array}
    \right) \nonumber\\
&=(a_{i1}x_{kj},\ldots,a_{in}x_{kj})
    \left(
        \begin{array}{lll}
             x_{k1}  \\
             ~~\vdots \\
             x_{kn} \\
        \end{array}
    \right)
   +(x_{k1},\ldots,x_{kn})
    \left(
        \begin{array}{lll}
             a_{i1}x_{kj} \\
             ~~~\vdots\\
             a_{in}x_{kj} \\
        \end{array}
    \right) \nonumber\\
&=\sum\limits_{r=1}^n a_{ir}x_{kj}x_{kr}+x_{kj}\sum\limits_{r=1}^n a_{ir}x_{kr} \nonumber\\
&=2x_{kj}\sum\limits_{r=1}^n a_{ir}x_{kr},
\end{align}
where $i,j\in[n]$, $k=0,1,2,\ldots,n$.

In the following, we give the estimate of variational equation (\ref{10}), which plays an important role in analysis of the optimal principle. The result can be shown by induction. For the case $k=0$, from the fact that $y_1=g(0,y_0)=g(0,Ax_0)$, we have
\begin{align}\label{14}
\dfrac{\partial y_1}{\partial a_{ij}}
&=\dfrac{\partial}{\partial {a_{ij}}}g(0,Ax_0)\
=\dfrac{\partial y_1}{\partial y_0^{\rm T}}\dfrac{\partial y_0}{\partial a_{ij}}\nonumber\\
&=\left(\begin{array}{lll}
\dfrac{\partial y_{11}}{\partial y_{01}}\cdots\dfrac{\partial y_{11}}{\partial y_{0n}}\\
~~\vdots~~~~\ddots ~~\vdots\\
\dfrac{\partial y_{1n}}{\partial y_{01}}\cdots\dfrac{\partial y_{1n}}{\partial y_{0n}}
           \end{array}
     \right)
     \left(\begin{array}{lllllll}
       ~0\\
       ~\vdots\\
       ~0\\
       x_{0j}\\
       ~0\\
       ~\vdots\\
       ~0
     \end{array}
     \right)
=\left(\begin{array}{lll}
     \dfrac{\partial y_{11}}{\partial y_{0i}}x_{0j}\\
     ~~~~~\vdots\\
     \dfrac{\partial y_{1n}}{\partial y_{0i}}x_{0j}
     \end{array}
     \right)
=\left(\begin{array}{lll}
     \dfrac{\partial y_{11}}{\partial y_{0i}}\\
     ~~~\vdots\\
     \dfrac{\partial y_{1n}}{\partial y_{0i}}
     \end{array}
     \right) x_{0j}.
\end{align}
When $k=1$, it follows from the iterative formula $y_2=g(1,y_1)=g(1,g(0,Ax_0))$ and (\ref{14}) that
\begin{align*}
\dfrac{\partial y_2}{\partial a_{ij}}
&=\dfrac{\partial}{\partial {a_{ij}}}g(1,g(0,Ax_0))
=\dfrac{\partial y_2}{\partial y_1^{\rm T}}
    \dfrac{\partial y_1}{\partial a_{ij}} \\
&=\left(\begin{array}{lll}
\dfrac{\partial y_{21}}{\partial y_{11}}\cdots\dfrac{\partial y_{21}}{\partial y_{1n}}\\
~~\vdots~~~~\ddots ~~\vdots\\
\dfrac{\partial y_{2n}}{\partial y_{11}}\cdots\dfrac{\partial y_{2n}}{\partial y_{1n}}
           \end{array}
     \right)
\left(\begin{array}{lll}
     \dfrac{\partial y_{11}}{\partial y_{0i}}\\
     ~~~\vdots\\
     \dfrac{\partial y_{1n}}{\partial y_{0i}}
     \end{array}
     \right) x_{0j} \\
&=\left(\begin{array}{lll}
\dfrac{\partial y_{21}}{\partial y_{11}}\dfrac{\partial y_{11}}{\partial y_{0i}}+
\cdots+
\dfrac{\partial y_{21}}{\partial y_{1n}}\dfrac{\partial y_{1n}}{\partial y_{0i}}\\
~~~~~~~~~~~~~~~~~\vdots\\
\dfrac{\partial y_{2n}}{\partial y_{11}}\dfrac{\partial y_{11}}{\partial y_{0i}}+
\cdots+
\dfrac{\partial y_{2n}}{\partial y_{1n}}\dfrac{\partial y_{1n}}{\partial y_{0i}}
           \end{array}
     \right)x_{0j}
=\left(\begin{array}{lll}
\sum\limits_{t=1}^n
\dfrac{\partial y_{21}}{\partial y_{1t}}\dfrac{\partial y_{1t}}{\partial y_{0i}}\\
~~~~~~~\vdots\\
\sum\limits_{t=1}^n
\dfrac{\partial y_{2n}}{\partial y_{1t}}\dfrac{\partial y_{1t}}{\partial y_{0i}}
           \end{array}
     \right)x_{0j}.
\end{align*}
Let us make the induction hypothesis, assuming that the following expression is true for $k-1$, that is
\begin{align}\label{15}
\dfrac{\partial y_k}{\partial a_{ij}}
=\dfrac{\partial}{\partial {a_{ij}}}g(k-1,y_{k-1})
=\left(\begin{array}{lll}
\sum\limits_{t=1}^n
\dfrac{\partial y_{k1}}{\partial y_{(k-1)t}}\dfrac{\partial y_{(k-1)t}}{\partial y_{(k-2)t}} \cdots \dfrac{\partial y_{1t}}{\partial y_{0i}}\\
~~~~~~~~~~~~~~~~\vdots\\
\sum\limits_{t=1}^n
\dfrac{\partial y_{kn}}{\partial y_{(k-1)t}}\dfrac{\partial y_{(k-1)t}}{\partial y_{(k-2)t}} \cdots \dfrac{\partial y_{1t}}{\partial y_{0i}}\\
           \end{array}
     \right)x_{0j}.
\end{align}
We now show that it continues to hold for $k$. By combining (\ref{10}) with (\ref{15}), we obtain
\begin{align*}
\dfrac{\partial y_{k+1}}{\partial a_{ij}}
&=\dfrac{\partial}{\partial {a_{ij}}}g(k,y_k)
=\dfrac{\partial y_{k+1}}{\partial y_k^{\rm T}}
    \dfrac{\partial}{\partial a_{ij}}g(k-1,y_{k-1})\\
&=\left(\begin{array}{lll}
\dfrac{\partial y_{(k+1)1}}{\partial y_{k1}}\cdots\dfrac{\partial y_{(k+1)1}}{\partial y_{kn}}\\
~~~~~\vdots~~~~~\ddots ~~~~~~\vdots\\
\dfrac{\partial y_{(k+1)n}}{\partial y_{k1}}\cdots\dfrac{\partial y_{(k+1)n}}{\partial y_{kn}}
           \end{array}
     \right)
     \left(\begin{array}{lll}
\sum\limits_{t=1}^n
\dfrac{\partial y_{k1}}{\partial y_{(k-1)t}}\dfrac{\partial y_{(k-1)t}}{\partial y_{(k-2)t}} \cdots \dfrac{\partial y_{1t}}{\partial y_{0i}}\\
~~~~~~~~~~~~~~~~\vdots\\
\sum\limits_{t=1}^n
\dfrac{\partial y_{kn}}{\partial y_{(k-1)t}}\dfrac{\partial y_{(k-1)t}}{\partial y_{(k-2)t}} \cdots \dfrac{\partial y_{1t}}{\partial y_{0i}}\\
           \end{array}
     \right)x_{0j} \\
&=\left(\begin{array}{lll}
\sum\limits_{t=1}^n
\dfrac{\partial y_{(k+1)1}}{\partial y_{kt}}\dfrac{\partial y_{kt}}{\partial y_{(k-1)t}} \cdots \dfrac{\partial y_{1t}}{\partial y_{0i}}\\
~~~~~~~~~~~~~~~~\vdots\\
\sum\limits_{t=1}^n
\dfrac{\partial y_{(k+1)n}}{\partial y_{kt}}\dfrac{\partial y_{kt}}{\partial y_{(k-1)t}} \cdots \dfrac{\partial y_{1t}}{\partial y_{0i}}\\
           \end{array}
     \right)x_{0j}.
\end{align*}
Then, if $\{y_k\}$ is generated by (\ref{10}), we deduce that
\begin{equation}\label{16}
\dfrac{\partial y_{k+1}}{\partial a_{ij}}
=\dfrac{\partial}{\partial a_{ij}}g(k,y_k)
=\left(\begin{array}{lll}
\sum\limits_{t=1}^n
\dfrac{\partial y_{(k+1)1}}{\partial y_{kt}}\dfrac{\partial y_{kt}}{\partial y_{(k-1)t}} \cdots \dfrac{\partial y_{1t}}{\partial y_{0i}}\\
~~~~~~~~~~~~~~~~\vdots\\
\sum\limits_{t=1}^n
\dfrac{\partial y_{(k+1)n}}{\partial y_{kt}}\dfrac{\partial y_{kt}}{\partial y_{(k-1)t}} \cdots \dfrac{\partial y_{1t}}{\partial y_{0i}}\\
           \end{array}
     \right)x_{0j}.
\end{equation}

For the second term of the right-hand side in (\ref{12}), we obtain the following equation by means of derivative rule of compound function,
\begin{equation}\label{17}
\dfrac{\partial}{\partial a_{ij}}(x_k^{\rm T}A^{\rm T}y_k)
=x_k^{\rm T}\left(\dfrac{\partial}{\partial {a_{ij}}}A^{\rm T}\right)y_k
 +x_k^{\rm T}A^{\rm T}\left(\dfrac{\partial}{\partial {a_{ij}}}y_k\right).
\end{equation}
It is obvious that
\begin{equation}\label{18}
x_k^{\rm T}\left(\dfrac{\partial}{\partial a_{ij}}A^{\rm T}\right)y_k
=(0,\ldots,0,x_{kj},0,\ldots,0)y_k
=x_{kj}y_{ki}.
\end{equation}
Together with (\ref{16}), it yields that
\begin{align}\label{19}
x_k^{\rm T}A^{\rm T}\left(\dfrac{\partial}{\partial {a_{ij}}}y_k\right)
=&x_k^{\rm T}A^{\rm T} \left(\begin{array}{lll}
\sum\limits_{t=1}^n
\dfrac{\partial y_{k1}}{\partial y_{(k-1)t}}\dfrac{\partial y_{(k-1)t}}{\partial y_{(k-2)t}}
\cdots \dfrac{\partial y_{1t}}{\partial y_{0i}}\\
~~~~~~~~~~~~~~~~~\vdots\\
\sum\limits_{t=1}^n
\dfrac{\partial y_{kn}}{\partial y_{(k-1)t}}\dfrac{\partial y_{(k-1)t}}{\partial y_{(k-2)t}}
\cdots \dfrac{\partial y_{1t}}{\partial y_{0i}}
           \end{array}
     \right)x_{0j} \nonumber\\
=& x_k^{\rm T}
\left(\begin{array}{lll}
a_{11}\sum\limits_{t=1}^n
\dfrac{\partial y_{k1}}{\partial y_{(k-1)t}}\cdots \dfrac{\partial y_{1t}}{\partial y_{0i}}
+\cdots
+a_{n1}\sum\limits_{t=1}^n\dfrac{\partial y_{kn}}{\partial y_{(k-1)t}}
\cdots \dfrac{\partial y_{1t}}{\partial y_{0i}}\\
~~~~~~~~~~~~~~~~~~~~~~~~~~~~~~~~~~~~~~~~~~~~~~~\vdots\\
a_{1n}\sum\limits_{t=1}^n
\dfrac{\partial y_{k1}}{\partial y_{(k-1)t}}
\cdots \dfrac{\partial y_{1t}}{\partial y_{0i}}
+\cdots+a_{nn}\sum\limits_{t=1}^n\dfrac{\partial y_{kn}}{\partial y_{(k-1)t}}
\cdots \dfrac{\partial y_{1t}}{\partial y_{0i}}
                       \end{array}
                       \right) x_{0j} \nonumber\\
=& x_{k1}\left(a_{11}\sum\limits_{t=1}^n
\dfrac{\partial y_{k1}}{\partial y_{(k-1)t}}
\cdots \dfrac{\partial y_{1t}}{\partial y_{0i}}
+\cdots+a_{n1}\sum\limits_{t=1}^n\dfrac{\partial y_{kn}}{\partial y_{(k-1)t}}
\cdots \dfrac{\partial y_{1t}}{\partial y_{0i}} \right)x_{0j}  \nonumber\\
& +\cdots \nonumber\\
& +x_{kn}\left(a_{1n}\sum\limits_{t=1}^n
\dfrac{\partial y_{k1}}{\partial y_{(k-1)t}}
\cdots \dfrac{\partial y_{1t}}{\partial y_{0i}}
+\cdots+a_{nn}\sum\limits_{t=1}^n\dfrac{\partial y_{kn}}{\partial y_{(k-1)t}}
\cdots \dfrac{\partial y_{1t}}{\partial y_{0i}} \right)x_{0j} \nonumber\\
=&x_{0j}\sum\limits_{r=1}^n
           \sum\limits_{s=1}^n
           \sum\limits_{t=1}^n x_{kr}a_{sr}
\left(\dfrac{\partial y_{ks}}{\partial y_{(k-1)t}}
\dfrac{\partial y_{(k-1)t}}{\partial y_{(k-2)t}}
\cdots\dfrac{\partial y_{1t}}{\partial y_{0i}}\right).
\end{align}

Using (\ref{18}) and (\ref{19}) in (\ref{17}), we can easily get
\begin{align}\label{20}
\dfrac{\partial}{\partial a_{ij}}\left(x_k^{\rm T}A^{\rm T}y_k\right)
=x_{kj}y_{ki}+x_{0j}\sum\limits_{r=1}^n
           \sum\limits_{s=1}^n
           \sum\limits_{t=1}^n x_{kr}a_{sr}
\left(\dfrac{\partial y_{ks}}{\partial y_{(k-1)t}}
\dfrac{\partial y_{(k-1)t}}{\partial y_{(k-2)t}}
\cdots\dfrac{\partial y_{1t}}{\partial y_{0i}}\right),
\end{align}
where $i,j\in[n]$, $k=0,1,2,\ldots,n$.

For the last term of the right-hand side in (\ref{12}), it is enough to compute that
\begin{align}\label{21}
\dfrac{\partial}{\partial a_{ij}}(y_k^{\rm T}y_k)
=&\left(\dfrac{\partial}{\partial a_{ij}}y_k^{\rm T}\right)y_k
    +y_k^{\rm T}\left(\dfrac{\partial}{\partial a_{ij}}y_k\right) \nonumber\\
=&\left(
\sum\limits_{t=1}^n
\dfrac{\partial y_{k1}}{\partial y_{(k-1)t}} \cdots \dfrac{\partial y_{1t}}{\partial y_{0i}},
\cdots,
\sum\limits_{t=1}^n
\dfrac{\partial y_{kn}}{\partial y_{(k-1)t}} \cdots \dfrac{\partial y_{1t}}{\partial y_{0i}}
     \right)x_{0j}y_k \nonumber\\
&+y_k^{\rm T}\left(\begin{array}{lll}
\sum\limits_{t=1}^n
\dfrac{\partial y_{k1}}{\partial y_{(k-1)t}}\dfrac{\partial y_{(k-1)t}}{\partial y_{(k-2)t}} \cdots \dfrac{\partial y_{1t}}{\partial y_{0i}}\\
~~~~~~~~~~~~~~~~\vdots\\
\sum\limits_{t=1}^n
\dfrac{\partial y_{kn}}{\partial y_{(k-1)t}}\dfrac{\partial y_{(k-1)t}}{\partial y_{(k-2)t}} \cdots \dfrac{\partial y_{1t}}{\partial y_{0i}}\\
           \end{array}
     \right)x_{0j} \nonumber\\
=& 2x_{0j}\left(
y_{k1}\sum\limits_{t=1}^n
\dfrac{\partial y_{k1}}{\partial y_{(k-1)t}} \cdots \dfrac{\partial y_{1t}}{\partial y_{0i}}+
\cdots+
y_{kn}\sum\limits_{t=1}^n
\dfrac{\partial y_{kn}}{\partial y_{(k-1)t}} \cdots \dfrac{\partial y_{1t}}{\partial y_{0i}}
     \right) \nonumber\\
=& 2x_{0j}\sum\limits_{s=1}^n\sum\limits_{t=1}^n y_{ks}
\left(\dfrac{\partial y_{ks}}{\partial y_{(k-1)t}}\dfrac{\partial y_{(k-1)t}}{\partial y_{(k-2)t}} \cdots \dfrac{\partial y_{1t}}{\partial y_{0i}}\right),
\end{align}
where $i,j\in[n]$, $k=0,1,2,\ldots,n$.

In this way, one of the main results in this paper is built up by substituting (\ref{13}), (\ref{20}) and (\ref{21}) into (\ref{12}).

\begin{proposition} \label{pro:1}
Let $\{x_k\}$ and $\{y_k\}$ be generated by (\ref{7}) and (\ref{8}), respectively. If matrix $A$ is the optimal solutions of (\ref{11}), then there exists the optimal principle via first-order optimality condition:
\begin{align}\label{22}
x_{kj}\sum\limits_{r=1}^n a_{ir}x_{kr}
-x_{kj}y_{ki}
-x_{0j}\sum\limits_{r=1}^n\sum\limits_{s=1}^n \sum\limits_{t=1}^n x_{kr}a_{sr}
          \left( \dfrac{\partial y_{ks}}{\partial y_{(k-1)t}}
      \dfrac{\partial y_{(k-1)t}}{\partial y_{(k-2)t}}
      \cdots
      \dfrac{\partial y_{1t}}{\partial y_{0i}}\right) \nonumber\\
+x_{0j}\sum\limits_{s=1}^n \sum\limits_{t=1}^n y_{ks}
   \left(\dfrac{\partial y_{ks}}{\partial y_{(k-1)t}}
      \dfrac{\partial y_{(k-1)t}}{\partial y_{(k-2)t}}
      \cdots
      \dfrac{\partial y_{1t}}{\partial y_{0i}}
  \right)
=0,~~~~~~~~~~~~~~~~~~
\end{align}
where $x_{kj}$ denotes the $j$th component of column vector $x_k$ at step $k$, and other similar representations have the same meanings.
\end{proposition}

\subsection{Homotopy Dynamical Systems}
In light of above analysis, a natural extension of that is a more complex similarity study for general dynamical systems.

Consider hybrid systems formed by two chaotic attractors via homotopy method, which will show richer and more interesting dynamical behavior. Two simple systems as shown in (\ref{7}) and (\ref{8}) can be connected by constructing such a homotopy
\begin{equation}\label{23}
H(k,y_k,\lambda)=(1-\lambda)f(k,y_k)+\lambda g(k,y_k),
\end{equation}
where $\lambda\in[0,1]$ is an embedding parameter. Note that the homotopy $H$ is exactly a path connecting $f$ and $g$ such that $H(k,y_k,0)=f(k,y_k)$ and $H(k,y_k,1)=g(k,y_k)$. As the parameter $\lambda$ increases from 0 to 1, the homotopy $H$ varies continuously from one system to another.

Take into account two dynamical systems with the following forms
\begin{equation}\label{24}
x_{k+1}=H_1(k,x_k,\lambda_1),
\end{equation}
\begin{equation}\label{25}
y_{k+1}=H_2(k,y_k,\lambda_2),
\end{equation}
where homotopies (\ref{24}) and (\ref{25}) stand for general hybrid dynamical systems with the embedding parameters $\lambda_1,\lambda_2\in[0,1]$ changing per iteration. The solutions of (\ref{24}) and (\ref{25}) are denoted as
\begin{equation}\label{26}
x_k=H_1(k,x_0,\lambda_1),
\end{equation}
\begin{equation}\label{27}
y_k=H_2(k,y_0,\lambda_2),
\end{equation}
respectively. At this point, the similarity transformation matrix, becoming related to the parameter $\lambda=(\lambda_1,\lambda_2)^{\rm T}$, is re-expressed as $A=A(\lambda)$. Let the cost functional be written in the form
\begin{equation}\label{28}
J(A,\lambda)=\min\limits_{A,\lambda}\sum\limits_{k=0}^N\|Ax_k-y_k\|_2^2.
\end{equation}

According to Karush-Kuhn-Tucker (KKT for short) optimality conditions that are often checked for investigating whether a solution of nonlinear programming problem is optimal, $(A,\lambda)$ is a stationary point of (\ref{28}) if and only if
\begin{equation}\label{29}
\frac{\partial J(A,\lambda)}{\partial a_{ij}}
=\frac{\partial}{\partial a_{ij}}(x_k^{\rm T}A^{\rm T}Ax_k)
 -\frac{\partial}{\partial a_{ij}}(2x_k^{\rm T}A^{\rm T}y_k)
 +\frac{\partial}{\partial a_{ij}}(y_k^{\rm T}A^{\rm T}yx_k)
=0,
\end{equation}
and
\begin{equation}\label{30}
\frac{\partial J(A,\lambda)}{\partial \lambda}
=\frac{\partial}{\partial \lambda}(x_k^{\rm T}A^{\rm T}Ax_k)
 -\frac{\partial}{\partial \lambda}(2x_k^{\rm T}A^{\rm T}y_k)
 +\frac{\partial}{\partial \lambda}(y_k^{\rm T}A^{\rm T}yx_k)
=0,
\end{equation}
simultaneously. The derivation of (\ref{29}) is the same as (\ref{12}), and hence the details are omitted here.

For the first term of the right-hand side in (\ref{30}), by direct calculation, we get
\begin{align}\label{31}
\dfrac{\partial}{\partial \lambda}(x_k^{\rm T}A^{\rm T}Ax_k)
=&\left(\dfrac{\partial}{\partial \lambda}x_k^{\rm T}\right)A^{\rm T}Ax_k
   +x_k^{\rm T}A^{\rm T}A\left(\dfrac{\partial}{\partial \lambda}x_k\right) \nonumber\\
=&\left(\dfrac{\partial x_{k1}^{\rm T}}{\partial \lambda},\cdots,
        \dfrac{\partial x_{kn}^{\rm T}}{\partial \lambda}\right)A^{\rm T}Ax_k
   +x_k^{\rm T}A^{\rm T}A
   \left(
        \begin{array}{lll}
             \dfrac{\partial{x_{k1}}}{\partial{\lambda}} \\
             ~~~\vdots\\
             \dfrac{\partial{x_{kn}}}{\partial{\lambda}} \\
        \end{array}
    \right) \nonumber\\
=&\left(a_{11}\dfrac{\partial{x_{k1}}}{\partial{\lambda}}+\cdots
        +a_{1n}\dfrac{\partial{x_{kn}}}{\partial{\lambda}},
       \cdots,
       a_{n1}\dfrac{\partial{x_{k1}}}{\partial{\lambda}}+\cdots
       +a_{nn}\dfrac{\partial{x_{kn}}}{\partial{\lambda}}
  \right)
  \left(
        \begin{array}{lll}
             \sum\limits_{s=1}^na_{1s}x_{ks}\\
             ~~~~~~~\vdots\\
             \sum\limits_{s=1}^na_{ns}x_{ks} \\
        \end{array}
\right)  \nonumber\\
& +\left(\sum\limits_{s=1}^na_{1s}x_{ks},\cdots,\sum\limits_{s=1}^na_{ns}x_{ks}\right)
   \left(
        \begin{array}{lll}
           a_{11}\dfrac{\partial{x_{k1}}}{\partial{\lambda}}+\cdots
          +a_{1n}\dfrac{\partial{x_{kn}}}{\partial{\lambda}}\\
             ~~~~~~~~~~~~~~~~~~~\vdots\\
           a_{n1}\dfrac{\partial{x_{k1}}}{\partial{\lambda}}+\cdots
          +a_{nn}\dfrac{\partial{x_{kn}}}{\partial{\lambda}}\\
        \end{array}
        \right)  \nonumber\\
=& 2\sum\limits_{s=1}^na_{1s}x_{ks}
  \left(a_{11}\dfrac{\partial{x_{k1}}}{\partial{\lambda}}+\cdots
        +a_{1n}\dfrac{\partial{x_{kn}}}{\partial{\lambda}}\right)+\cdots \nonumber\\
&+2\sum\limits_{s=1}^na_{ns}x_{ks}
  \left(a_{n1}\dfrac{\partial{x_{k1}}}{\partial{\lambda}}+\cdots
       +a_{nn}\dfrac{\partial{x_{kn}}}{\partial{\lambda}}\right).
\end{align}

For the second term of the right-hand side in (\ref{30}), it follows from derivative rule of compound function that
\begin{align}\label{32}
\dfrac{\partial}{\partial \lambda}(x_k^{\rm T}A^{\rm T}y_k)
=\left(\dfrac{\partial}{\partial \lambda}x_k^{\rm T}\right)A^{\rm T}y_k
 +x_k^{\rm T}A^{\rm T}\left(\dfrac{\partial}{\partial \lambda}y_k\right).
\end{align}
On the one hand,
\begin{align}\label{33}
\left(\dfrac{\partial}{\partial \lambda}x_k^{\rm T}\right)A^{\rm T}y_k
=&\left(\dfrac{\partial{x_{k1}}}{\partial{\lambda}},\cdots,
       \dfrac{\partial{x_{kn}}}{\partial{\lambda}}\right)A^{\rm T}y_k \nonumber\\
=&\left(a_{11}\dfrac{\partial{x_{k1}}}{\partial{\lambda}}+\cdots
       +a_{1n}\dfrac{\partial{x_{kn}}}{\partial{\lambda}},
       \cdots,
       a_{n1}\dfrac{\partial{x_{k1}}}{\partial{\lambda}}+\cdots
       +a_{nn}\dfrac{\partial{x_{kn}}}{\partial{\lambda}}\right)
   \left(\begin{array}{lll}
             y_{k1}\\
             ~~\vdots\\
             y_{kn} \\
        \end{array}
   \right) \nonumber\\
=&\left(a_{11}\dfrac{\partial{x_{k1}}}{\partial{\lambda}}
       +\cdots
       +a_{1n}\dfrac{\partial{x_{kn}}}{\partial{\lambda}}\right)y_{k1}
 +\cdots
 +\left(a_{n1}\dfrac{\partial{x_{k1}}}{\partial{\lambda}}
       +\cdots
       +a_{nn}\dfrac{\partial{x_{kn}}}{\partial{\lambda}}\right)y_{kn} \nonumber\\
=&\sum\limits_{t=1}^na_{t1}y_{kt}\dfrac{\partial{x_{k1}}}{\partial{\lambda}}
 +\cdots
 +\sum\limits_{t=1}^na_{tn}y_{kt}\dfrac{\partial{x_{kn}}}{\partial{\lambda}}.
\end{align}
On the other hand,
\begin{align}\label{34}
x_k^{\rm T}A^{\rm T}\left(\dfrac{\partial}{\partial \lambda}y_k\right)
=&x_k^{\rm T}A^{\rm T}
   \left(\begin{array}{lll}
             \dfrac{\partial{y_{k1}}}{\partial{\lambda}}\\
             ~~~\vdots\\
             \dfrac{\partial{y_{kn}}}{\partial{\lambda}} \\
        \end{array}
   \right) \nonumber\\
=&\left(\sum\limits_{s=1}^na_{1s}x_{ks},\cdots,\sum\limits_{s=1}^na_{ns}x_{ks}\right)
    \left(\begin{array}{lll}
             \dfrac{\partial{y_{k1}}}{\partial{\lambda}}\\
             ~~~\vdots\\
             \dfrac{\partial{y_{kn}}}{\partial{\lambda}} \\
        \end{array}
   \right) \nonumber\\
=&\sum\limits_{s=1}^na_{1s}x_{ks}\dfrac{\partial{y_{k1}}}{\partial{\lambda}}
  +\cdots
  +\sum\limits_{s=1}^na_{ns}x_{ks}\dfrac{\partial{y_{kn}}}{\partial{\lambda}}.
\end{align}
Combining (\ref{33}) and (\ref{34}) with (\ref{32}), we obtain
\begin{align}\label{35}
\dfrac{\partial}{\partial \lambda}(x_k^{\rm T}A^{\rm T}y_k)
=\sum\limits_{t=1}^n a_{t1}y_{kt}\dfrac{\partial{x_{k1}}}{\partial \lambda}
+\cdots
+\sum\limits_{t=1}^n a_{tn}y_{kt}\dfrac{\partial{x_{kn}}}{\partial \lambda}
+\sum\limits_{s=1}^n a_{1s}x_{ks}\dfrac{\partial{y_{k1}}}{\partial \lambda}
+\cdots
+\sum\limits_{s=1}^n a_{ns}x_{ks}\dfrac{\partial{y_{kn}}}{\partial \lambda}.
\end{align}

For the last term of the right-hand side in (\ref{30}), it is obvious that
\begin{align}\label{36}
\dfrac{\partial}{\partial \lambda}(y_k^{\rm T}y_k)
=&\left(\dfrac{\partial}{\partial \lambda}y_k^{\rm T}\right)y_k
 +y_k^{\rm T}\left(\dfrac{\partial}{\partial \lambda}y_k\right) \nonumber\\
=&\left(\dfrac{\partial{y_{k1}}}{\partial \lambda},
        \cdots,
        \dfrac{\partial{y_{kn}}}{\partial \lambda}\right)
  \left(\begin{array}{lll}
             y_{k1}\\
             ~~\vdots\\
             y_{kn} \\
        \end{array}
   \right)
+(y_{k1},\cdots,y_{kn})
  \left(\begin{array}{lll}
             \dfrac{\partial{y_{k1}}}{\partial{\lambda}}\\
             ~~~\vdots\\
             \dfrac{\partial{y_{kn}}}{\partial{\lambda}} \\
        \end{array}
   \right) \nonumber\\
=&2y_{k1}\dfrac{\partial{y_{k1}}}{\partial \lambda}
 +\cdots
 +2y_{kn}\dfrac{\partial{y_{kn}}}{\partial \lambda}.
\end{align}

According to all derivations above, we deduce the other main result of this paper, which deals with how similar between orbits of general dynamical systems. By substituting (\ref{31}), (\ref{35}) and (\ref{36}) into (\ref{30}), we give the existence of the solution of (\ref{28}).

\begin{proposition} \label{pro:2}
Let $\{x_k\}$ and $\{y_k\}$ be generated by (\ref{24}) and (\ref{25}), respectively.
If matrix $A=A(\lambda)$ and parameter $\lambda$ are the optimal solutions of (\ref{28}), then there exist the following general optimal principle:
\begin{align}\label{37}
x_{kj}\sum\limits_{r=1}^n \bar{a}_{ir}x_{kr}
-x_{kj}y_{ki}
-x_{0j}\sum\limits_{r=1}^n\sum\limits_{s=1}^n \sum\limits_{t=1}^n x_{kr}\bar{a}_{sr}
          \left(\dfrac{\partial y_{ks}}{\partial y_{(k-1)t}}
      \dfrac{\partial y_{(k-1)t}}{\partial y_{(k-2)t}}
      \cdots
      \dfrac{\partial y_{1t}}{\partial y_{0i}}\right)\nonumber\\
~~~~~~~~~~~~~~~~+x_{0j}\sum\limits_{s=1}^n \sum\limits_{t=1}^n y_{ks}
   \left(\dfrac{\partial y_{ks}}{\partial y_{(k-1)t}}
      \dfrac{\partial y_{(k-1)t}}{\partial y_{(k-2)t}}
      \cdots
      \dfrac{\partial y_{1t}}{\partial y_{0i}}
  \right)
=0,  ~~~~~~~~~~~~~~~~~
\end{align}
and
\begin{align}\label{38}
\left(\sum\limits_{r=1}^n\sum\limits_{s=1}^n\bar{a}_{r1}\bar{a}_{rs}x_{ks}
      -\sum\limits_{t=1}^n\bar{a}_{t1}y_{kt}\right)\frac{\partial x_{k1}}{\partial\lambda}
+\cdots
+\left(\sum\limits_{r=1}^n\sum\limits_{s=1}^n\bar{a}_{rn}\bar{a}_{rs}x_{ks}
      -\sum\limits_{t=1}^n\bar{a}_{tn}y_{kt}\right)\frac{\partial x_{kn}}{\partial\lambda} \nonumber\\
-\left(\sum\limits_{s=1}^n\bar{a}_{1s}x_{ks}-y_{k1}\right)\frac{\partial y_{k1}}{\partial\lambda}
-\cdots
-\left(\sum\limits_{s=1}^n\bar{a}_{ns}x_{ks}-y_{kn}\right)\frac{\partial y_{kn}}{\partial\lambda}=0, ~~~~~~~~~~~~
\end{align}
where the similarity transformation matrix becomes related to the parameter $\lambda=(\lambda_1,\lambda_2)^{\rm T}$, and $\bar{a}_{ij}$ stands for some matrix element of $A=A(\lambda)$, other representations of components have the same meanings as before.
\end{proposition}

\begin{remark}
Although (\ref{37}) is formally consistent with (\ref{22}), in fact, each component of similarity transformation matrix $A$ is related to parameter $\lambda$.
\end{remark}

\begin{remark}
The existing literature on the variation of functional with respect to matrix rather than just vector or even scalar is rare, as considered in (\ref{11}) and (\ref{28}). Taking variation of a matrix can be converted to the partial derivative of each matrix elements, which is exactly the tedious calculations of derivation, especially when the number of iterations is large.
\end{remark}

For most cases, it is almost impossible to find similarity transformation matrix that makes solutions between two discrete dynamical systems exactly similar. By virtue of the characterizations of function
\begin{equation*}
h(\omega)=\frac{\log(1+\omega)}{\omega},~\omega\in\mathbb{R},
\end{equation*}
we define a similarity function by setting $\omega=\dfrac{1}{N}\sum\limits_{k=1}^N\|Ax_k-y_k\|_2^2$. When no similarity transformation matrix can be found to make two orbits completely similar, we also ask to what extent they are similar.

\begin{definition}\label{def2}
The similarity degree $\rho(A)$ of solutions between two discrete dynamical systems is defined as
\begin{equation}\label{39}
\rho(A)=\left\{\begin{array}{ll}
                ~~~~~ 1 ,& ~~{\rm{if}} \,\omega=0,\\
                \dfrac{\log(1+\omega)}{\omega}, & ~~{\rm{otherwise.}}\,\end{array}\right.
\end{equation}
It is easy to see that $\rho(A)$$\in$$(0,1]$, according to the above definition.
\end{definition}

We close this section by summing up that the solutions of two systems are said to be completely similar if $\rho(A)=1$, otherwise, some are similar.

\section{Experimental results}
\label{sec4}
In this section, some examples are given to show the utility of general optimal principle proposed in this paper. All codes are written in MATLAB R2021a and run on PC with 1.80 GHz CPU processor and 8.00 GB RAM memory. Unless otherwise specified, the numerical results are accurate to four decimal places throughout this paper.

To cope with over-fitting, L2-norm regularization is introduced naturally. In reinforcement learning and neural networks, this happens frequently when samples are limited and computation is expensive. The cost functional in (\ref{11}) or (\ref{28}) is augmented to include a L2-norm penalty of matrix $A$ with the following form
\begin{equation}\label{40}
\tilde{J}(A)=\min\limits_A\sum\limits_{k=0}^N\|Ax_k-y_k\|_2^2+\tau\|A\|_2^2,
\end{equation}
where $\tau$ is a positive constant called regularization parameter that balances the two objective terms.

In the optimal control theory, the most fundamental but crucial two optimization methods are Pontryagin's maximum principle and Bellman's dynamic programming. Taking into account the proposed optimal principles, we conduct the numerical simulations by following Pontryagin's maximum principle and Bellman's dynamic programming respectively for similarity of orbits between various chaotic systems.

The chaotic systems chosen in this section are all solved numerically by means of the widely used fourth-order Runge-Kutta method with time step size equal 0.01. The following numerical experiments mainly include two parts.

\subsection{Pontryagin's maximum principle}
As a necessary condition to solve optimal control problems, Pontryagin's maximum principle was proposed by Pontryagin and his group in the 1960s \cite{Pontryagin1962}. Outstanding feature of Pontryagin's maximum principle lies in that the optimal control signal transfering dynamical system from one state to another can be found under the condition that the state or input is fixed. The following examples concern to verify Pontryagin's maximum principle formulated in terms of the proposed optimal principle (\ref{22}) when studying the similarity of orbits between two chaotic attractors.

{\bf Example~4.1.}
Similarity of orbits between Lorenz attractor and Chua's circuit.

Let $\{x_k\}$ and $\{y_k\}$ be the numerical solutions derived from Lorenz and Chua systems for 2000 time steps (namely, $N$$=$$2000$) from the same initial states
\begin{equation*}
x_0=y_0=(0.1,0.1,0.1)^{\rm T}.
\end{equation*}
We divide the sequences into multi-stage decisions consisting of 10 steps for each (denoted by $N_1$$-$$N_{200}$), with the final state of previous stage as the initial condition of current stage. The optimal similarity transformation matrix of each stage is found by optimal principle (\ref{22}) whose accuracy performance is assessed by similarity degree (\ref{39}).

\begin{figure}[htbp]
\centering
\subfloat[Without~penalty.]{\includegraphics[scale=.4]{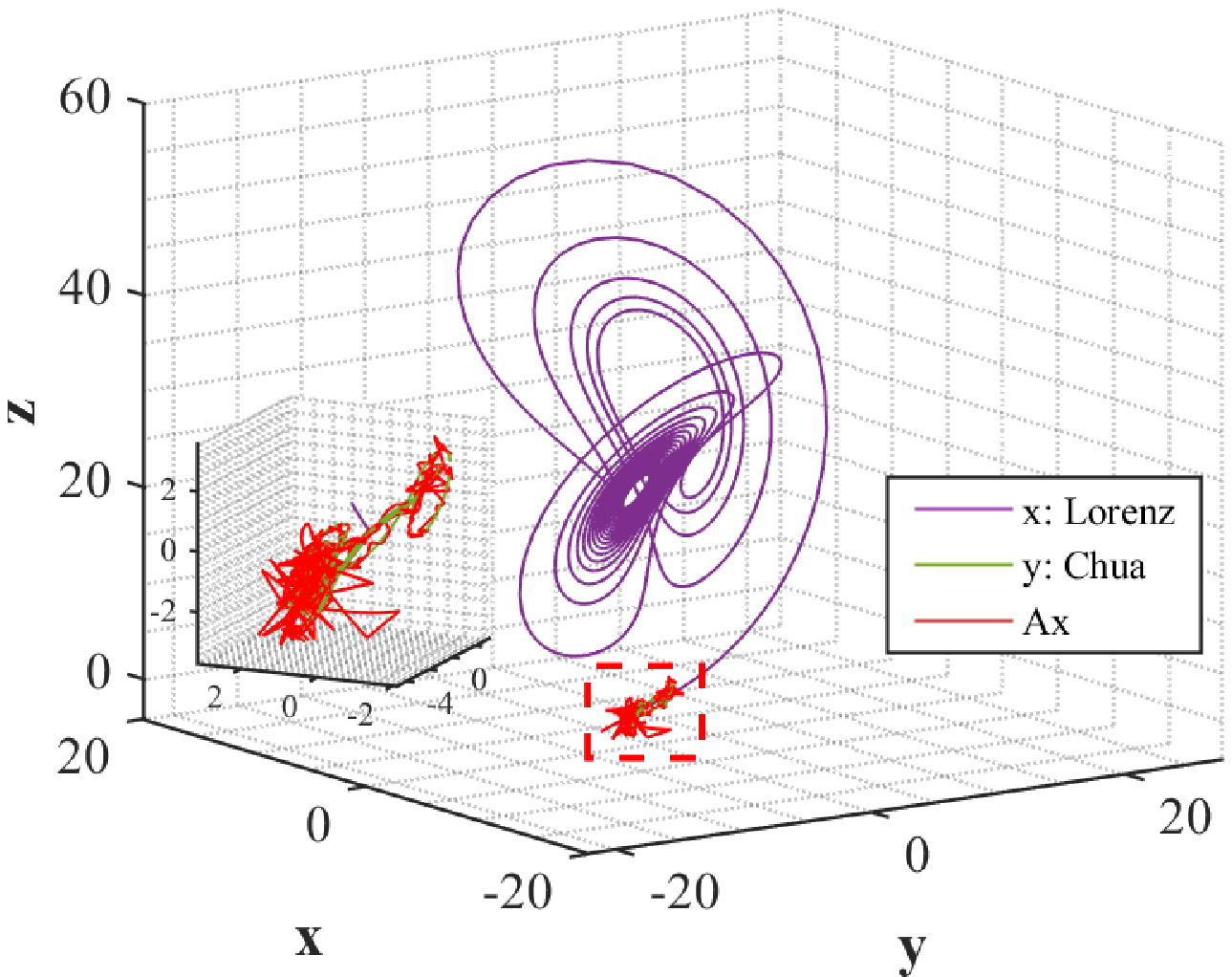}
\label{Fig.6a}}
\hfil
\subfloat[Based on L2-norm penalty.]{\includegraphics[scale=.4]{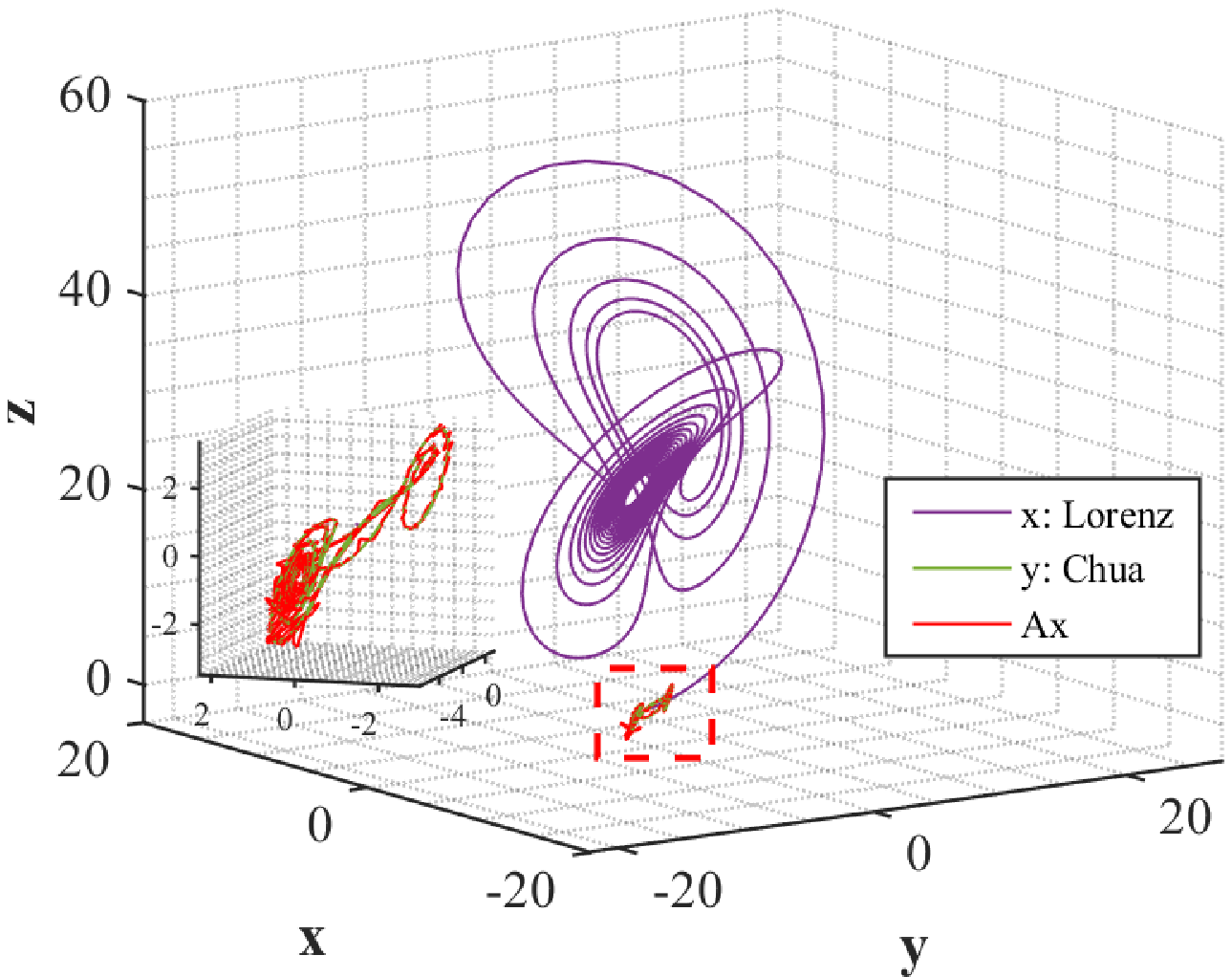}
\label{Fig.6b}}
\caption{Three dimensional stereograms of Example~4.1.}
\label{Fig.6}
\end{figure}

\begin{figure}[htbp]
\centering
\subfloat[Without~penalty.]{\includegraphics[scale=.28]{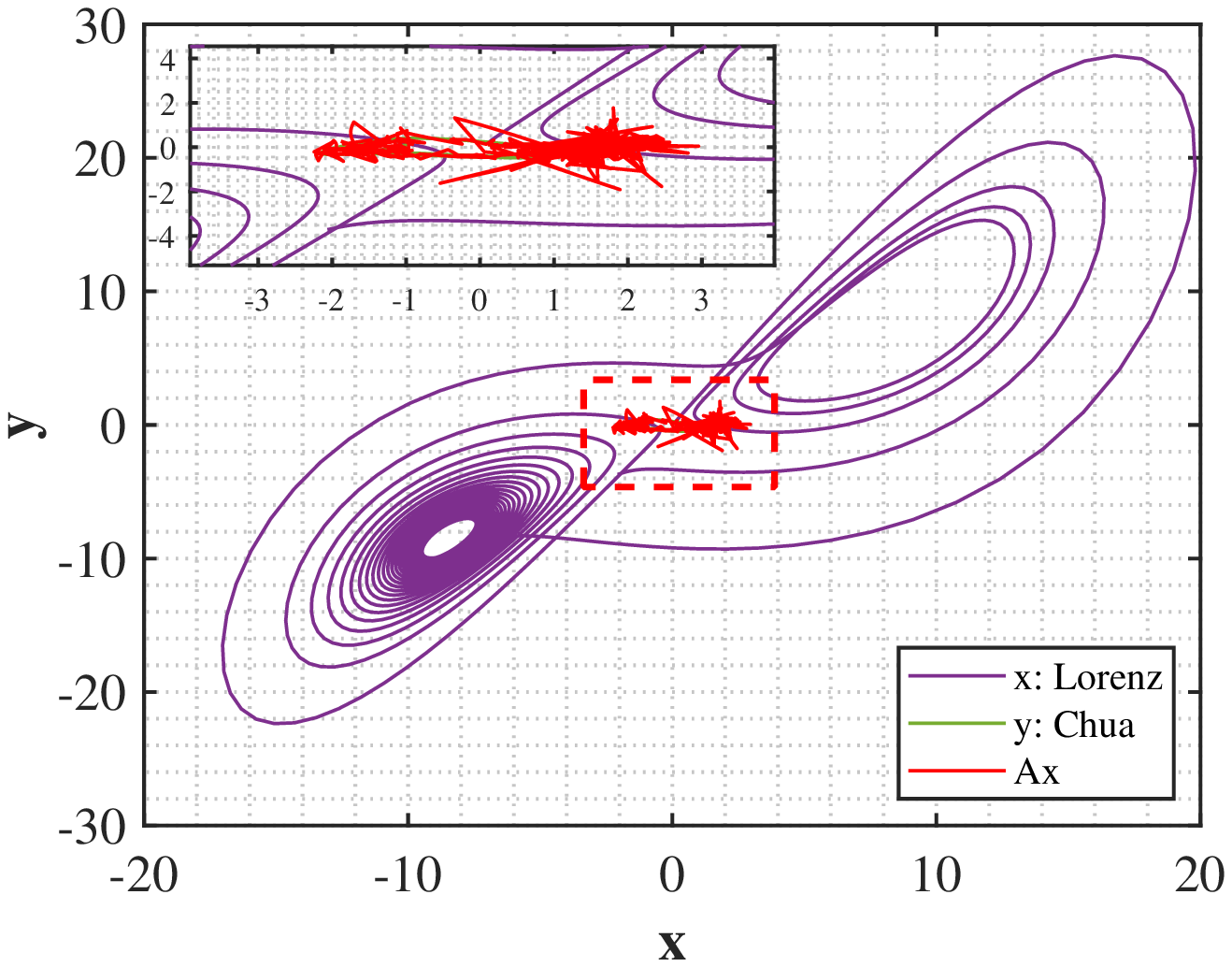}
\hfil
\includegraphics[scale=.28]{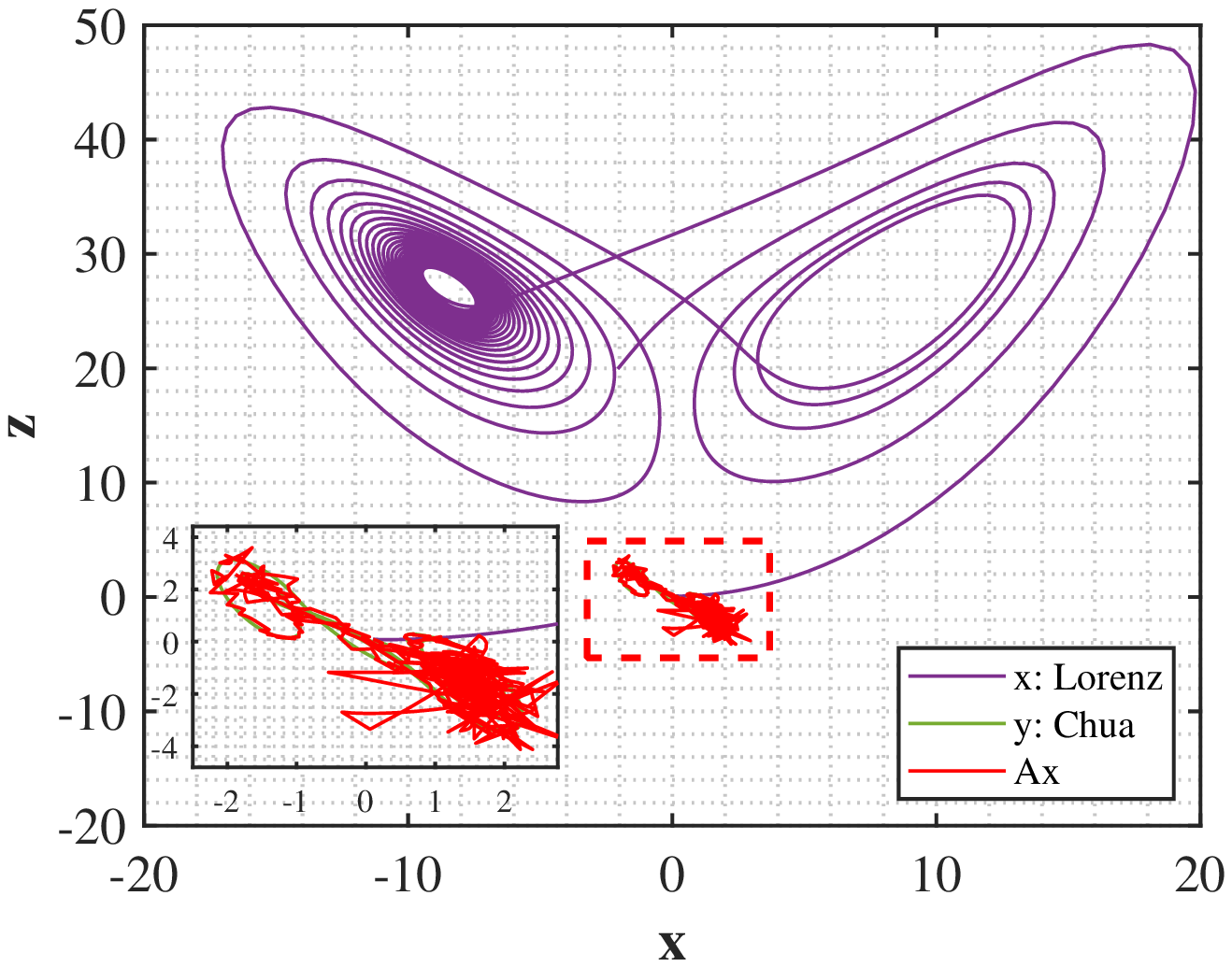}
\hfil
\includegraphics[scale=.28]{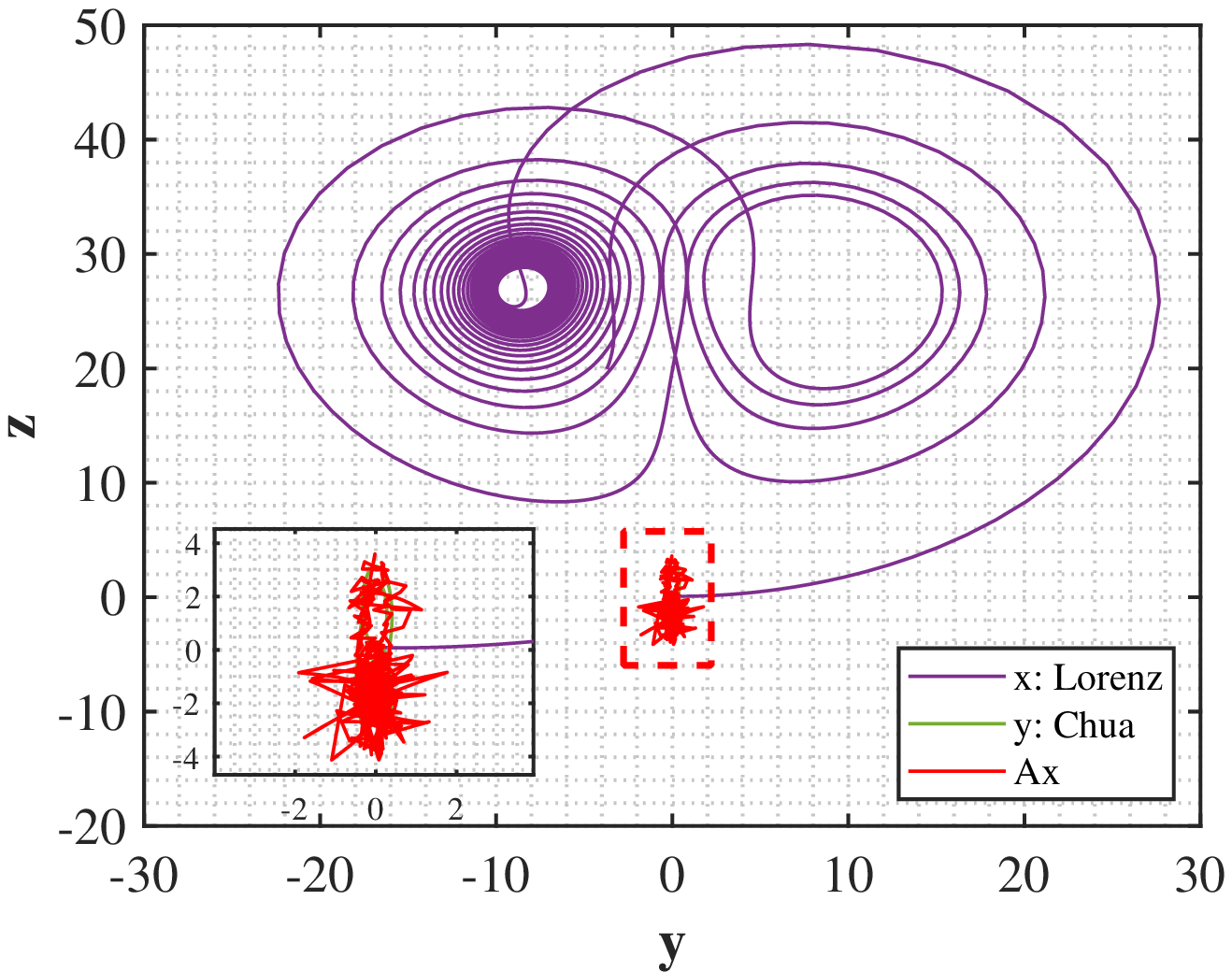}
\label{Fig.7a}}
\hfil
\subfloat[Based on L2-norm penalty.]{\includegraphics[scale=.28]{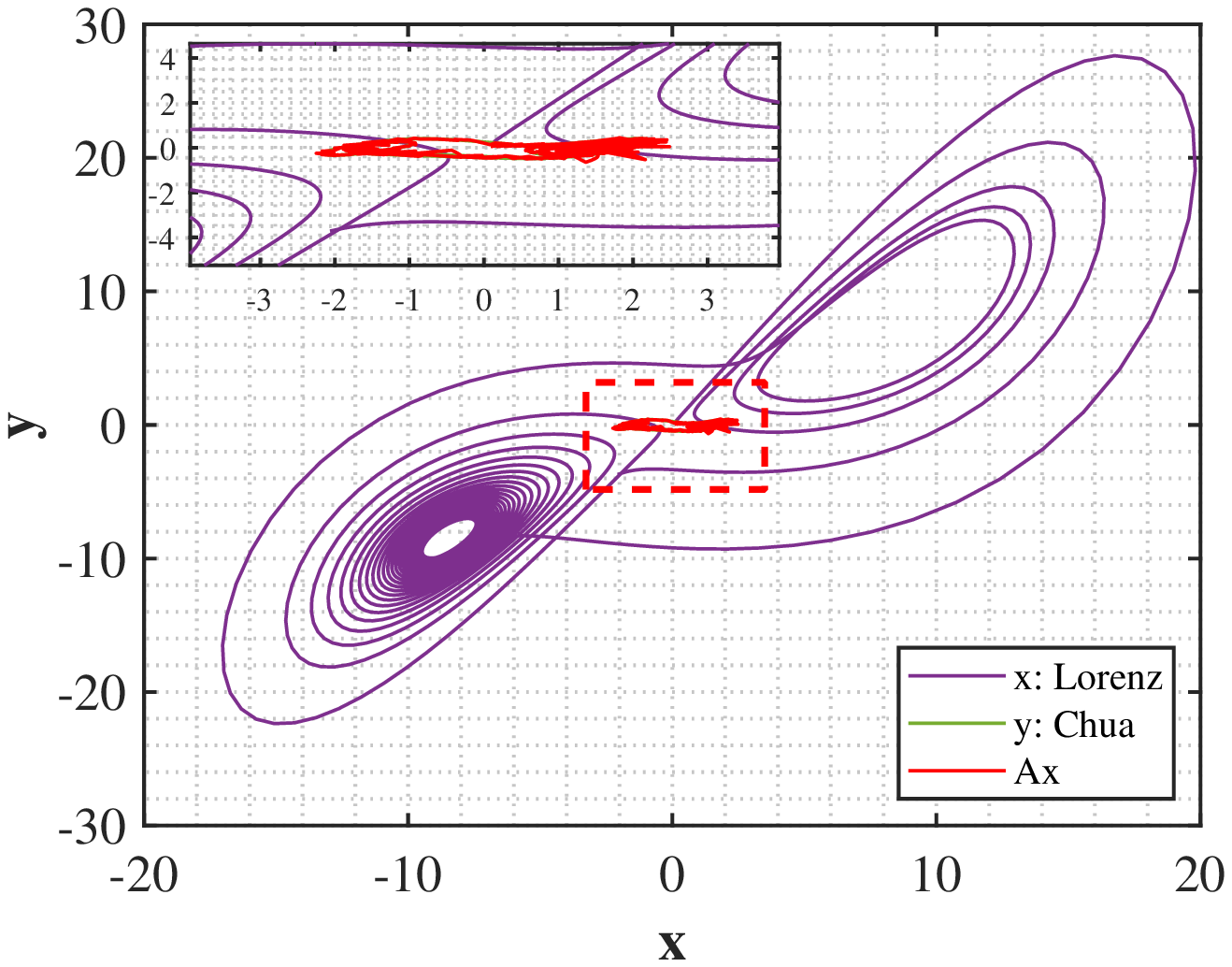}
\hfil
\includegraphics[scale=.28]{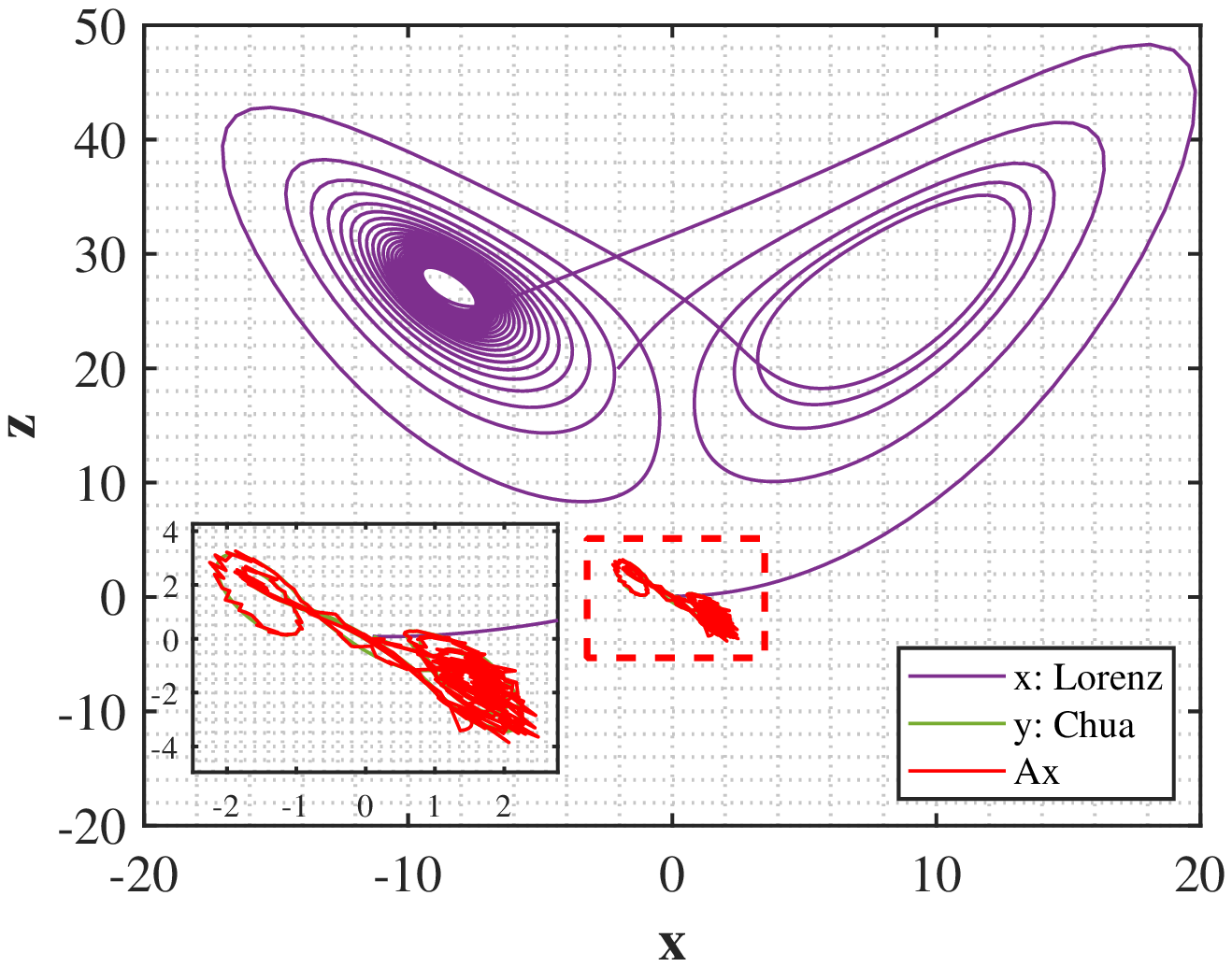}
\hfil
\includegraphics[scale=.28]{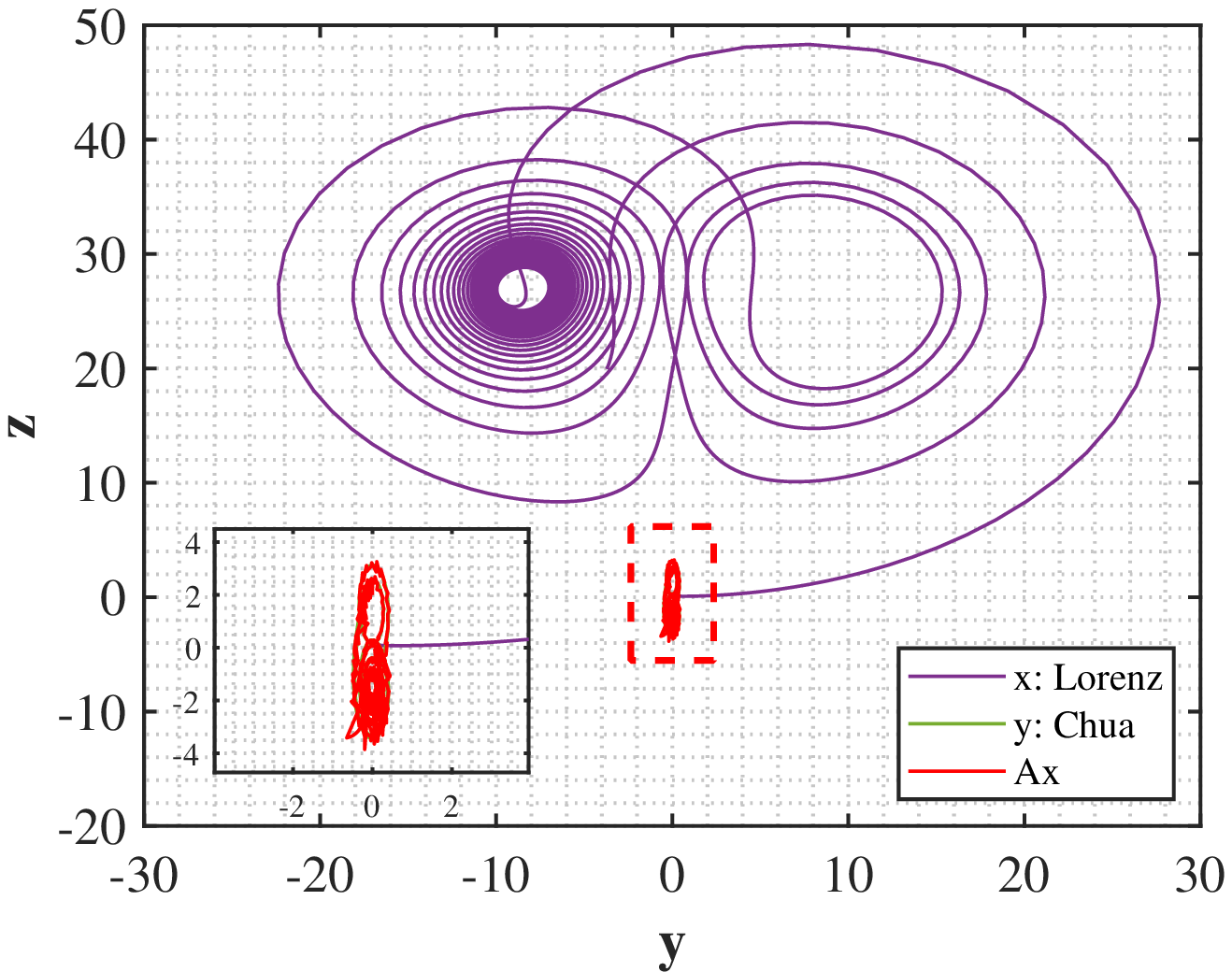}
\label{Fig.7b}}
\caption{Two dimensional plans of Example~4.1.}
\label{Fig.7}
\end{figure}

Figs. \ref{Fig.6}-\ref{Fig.7} describe three dimensional stereograms and two dimensional plans of Lorenz attractor, Chua's circuit and the trajectory acted by the optimal similarity transformation matrix for each stage. Enlarge the trajectories marked by dotted box so as to find similarity of orbits between $\{Ax_k\}$ and $\{y_k\}$ more clearly. We can observe that Lorenz attractor and Chua's circuit with different orbits become mainly similar after the optimal similarity transformation matrix is employed according to Figs. \ref{Fig.6a} and \ref{Fig.7a}. As can be seen in Figs. \ref{Fig.6b} and \ref{Fig.7b}, if the regularization parameter $\tau$ is selected as $10^{-4}$, the orbit $\{x_k\}$ under action of optimal similarity transformation matrix is surprisingly close to the orbit $\{y_k\}$, supporting the availability of tuning parameter and the stability of L2-norm regularization.

\begin{figure}[htbp]
\centering
\subfloat[Without~penalty.]{\includegraphics[scale=.4]{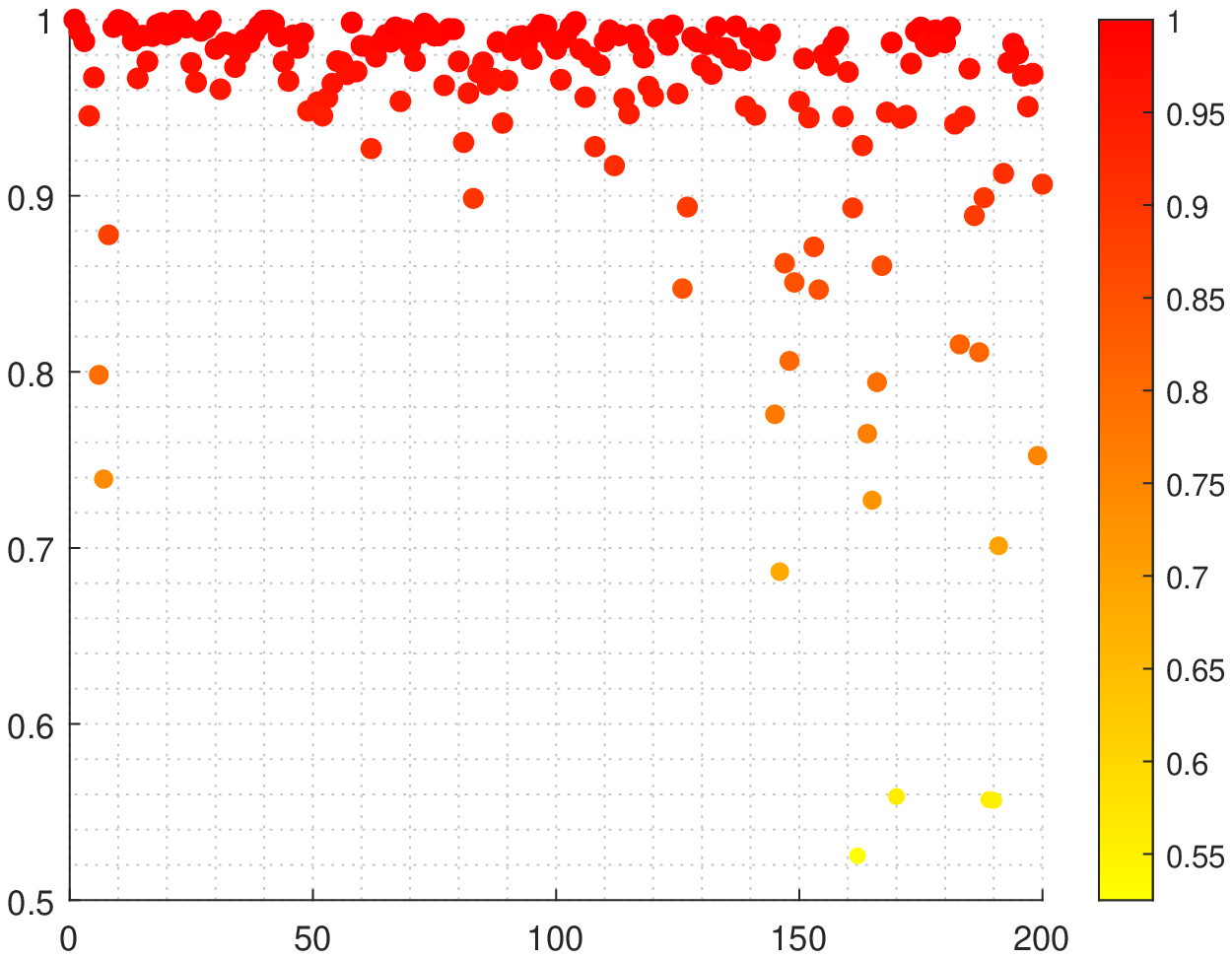}
\label{Fig.8a}}
\hfil
\subfloat[Based on L2-norm penalty.]{\includegraphics[scale=.4]{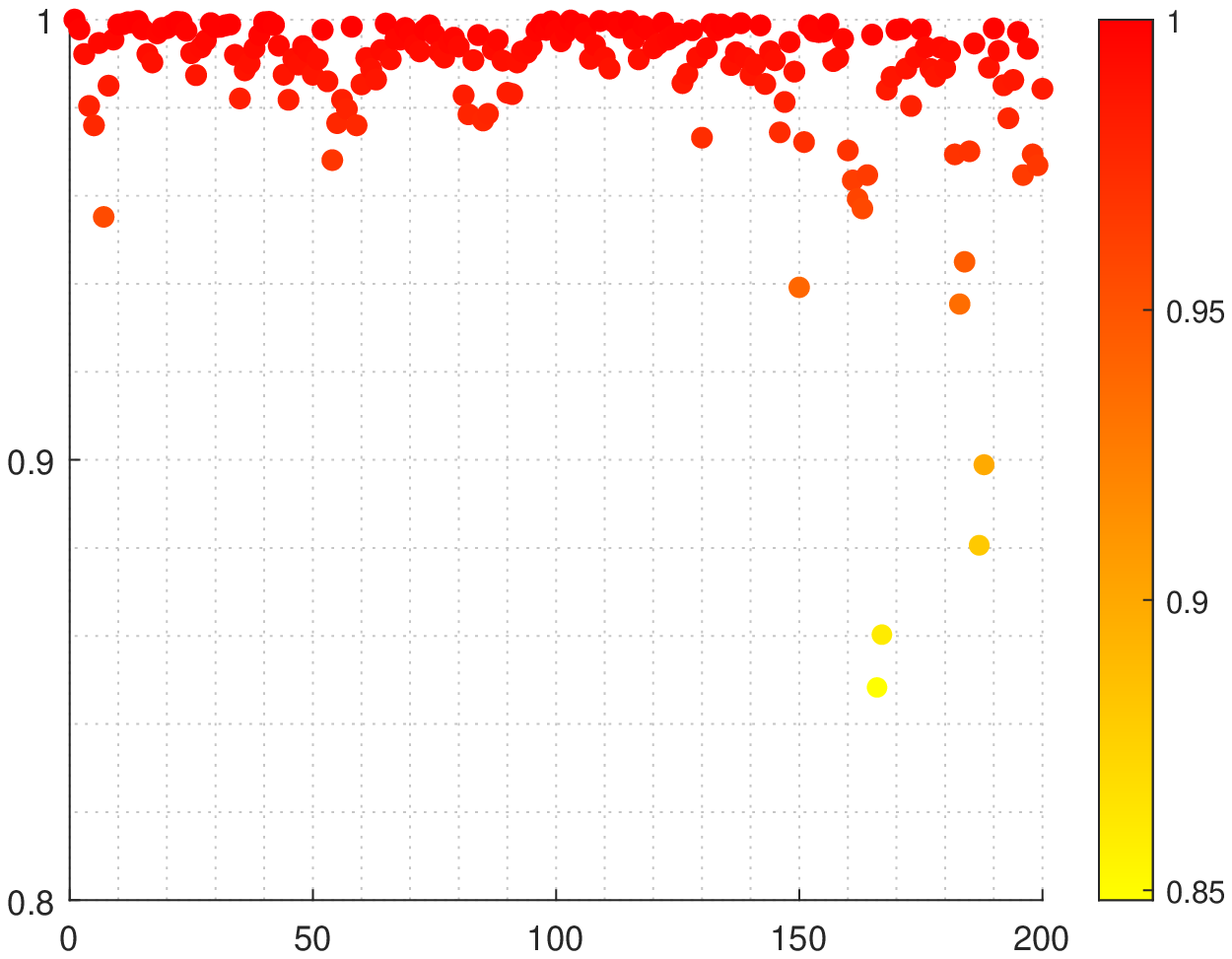}
\label{Fig.8b}}
\caption{Similarity degree of Example~4.1.}
\label{Fig.8}
\end{figure}

Results of similarity degree are shown in Fig. \ref{Fig.8}. From Fig. \ref{Fig.8a}, we can see that the values of similarity degree can achieve over 0.9 in the majority of results. When L2-norm penalty is introduced, only 4 results are less than 0.9. This exactly implies that the stability of solution is improved by introducing L2-norm regularization with suitable regularization parameter.

{\bf Example~4.2.}
Similarity of orbits between Lorenz attractor and R$\rm\ddot{o}ssler$ attractor.

Let $\{x_k\}$ and $\{y_k\}$ be the numerical solutions of Lorenz and R$\rm\ddot{o}ssler$ systems, sharing the same initial states $x_0$$=$$y_0$$=$$(0.1,0.1,0.1)^{\rm T}$. We perform the same experimental procedure as in Example 4.1.

\begin{figure}[htbp]
\centering
\subfloat[Without~penalty.]{\includegraphics[scale=.4]{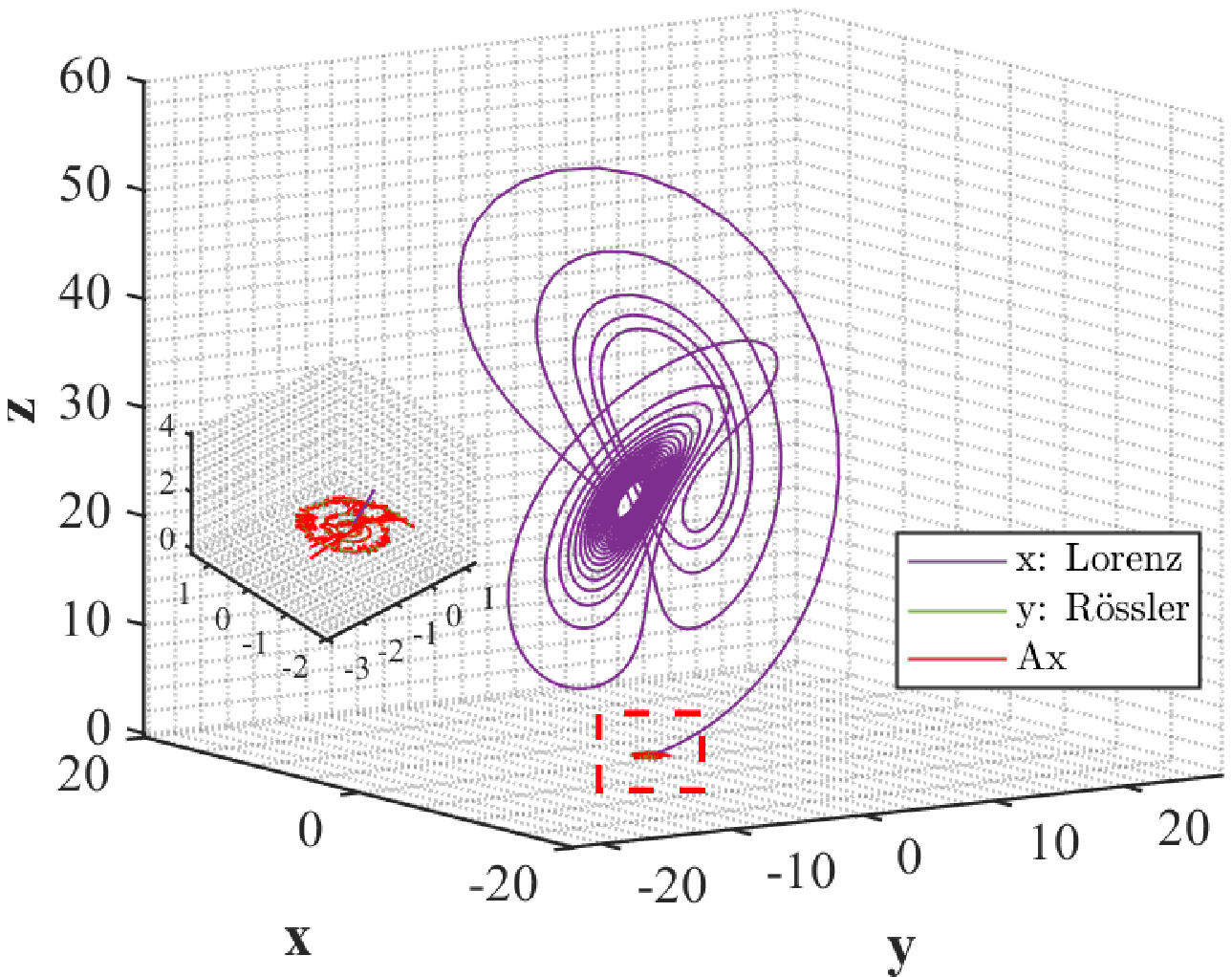}
\label{Fig.9a}}
\hfil
\subfloat[Based on L2-norm penalty.]{\includegraphics[scale=.4]{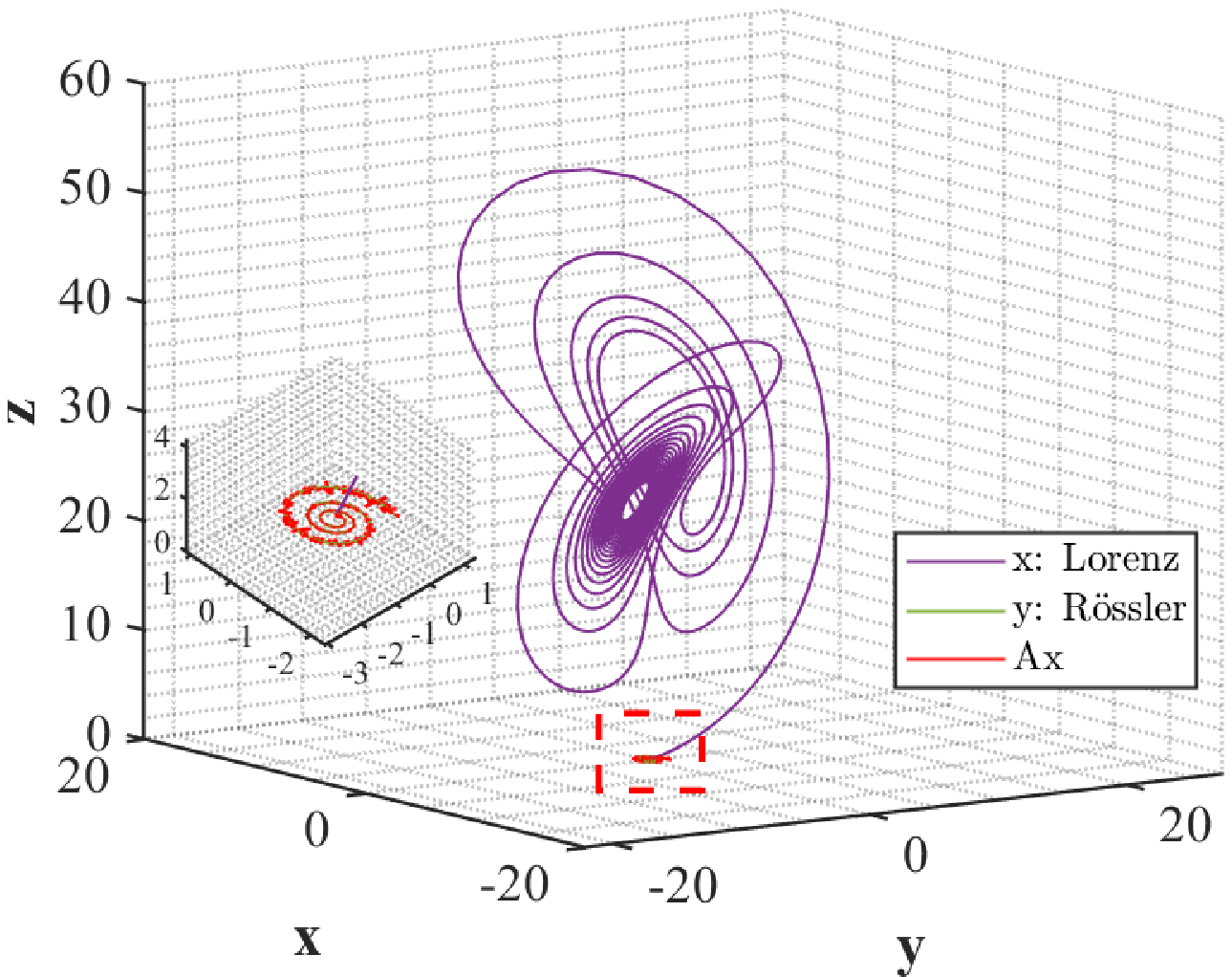}
\label{Fig.9b}}
\caption{Three dimensional stereograms of Example~4.2.}
\label{Fig.9}
\end{figure}

\begin{figure}[htbp]
\centering
\subfloat[Without~penalty.]{\includegraphics[scale=.28]{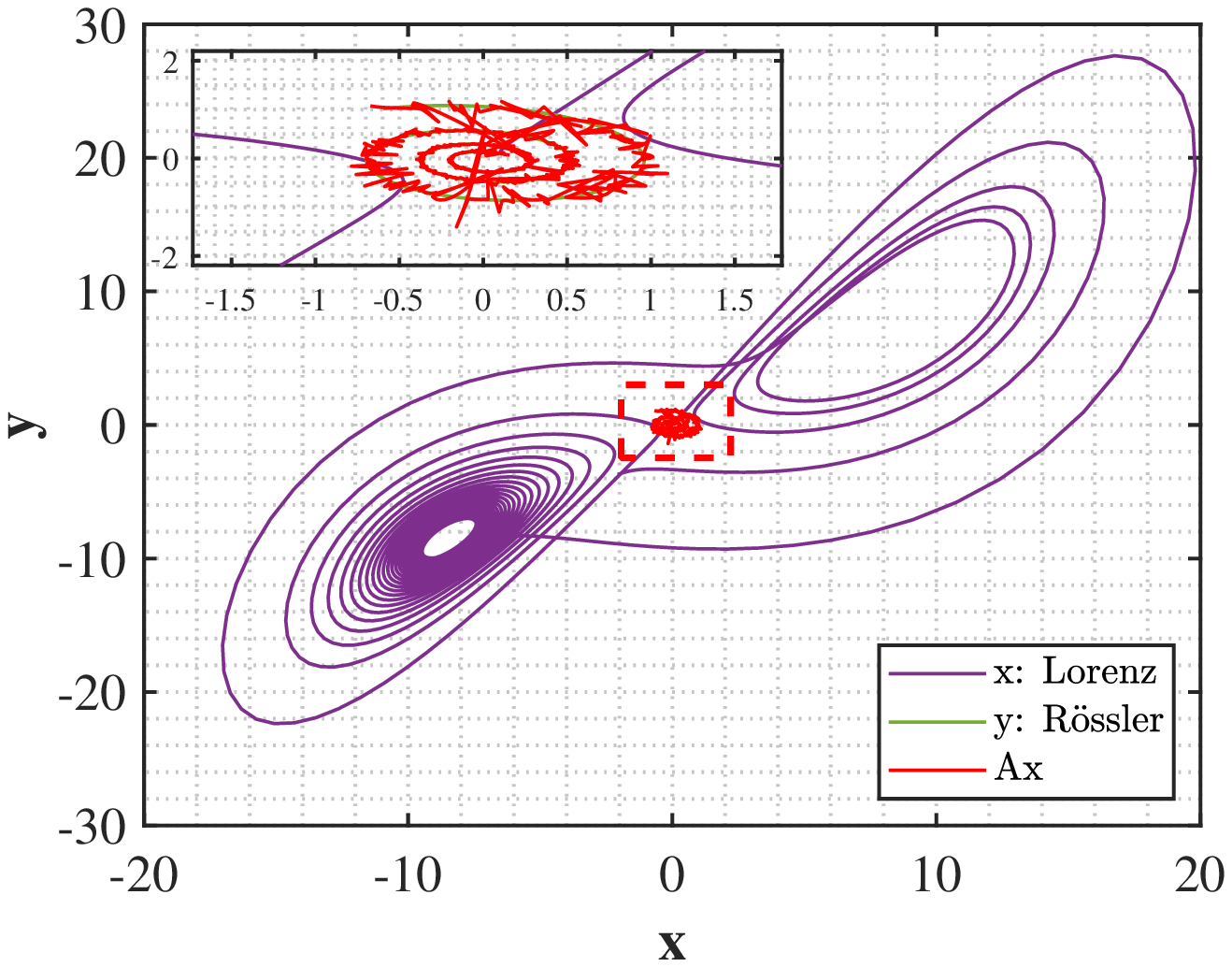}
\hfil
\includegraphics[scale=.28]{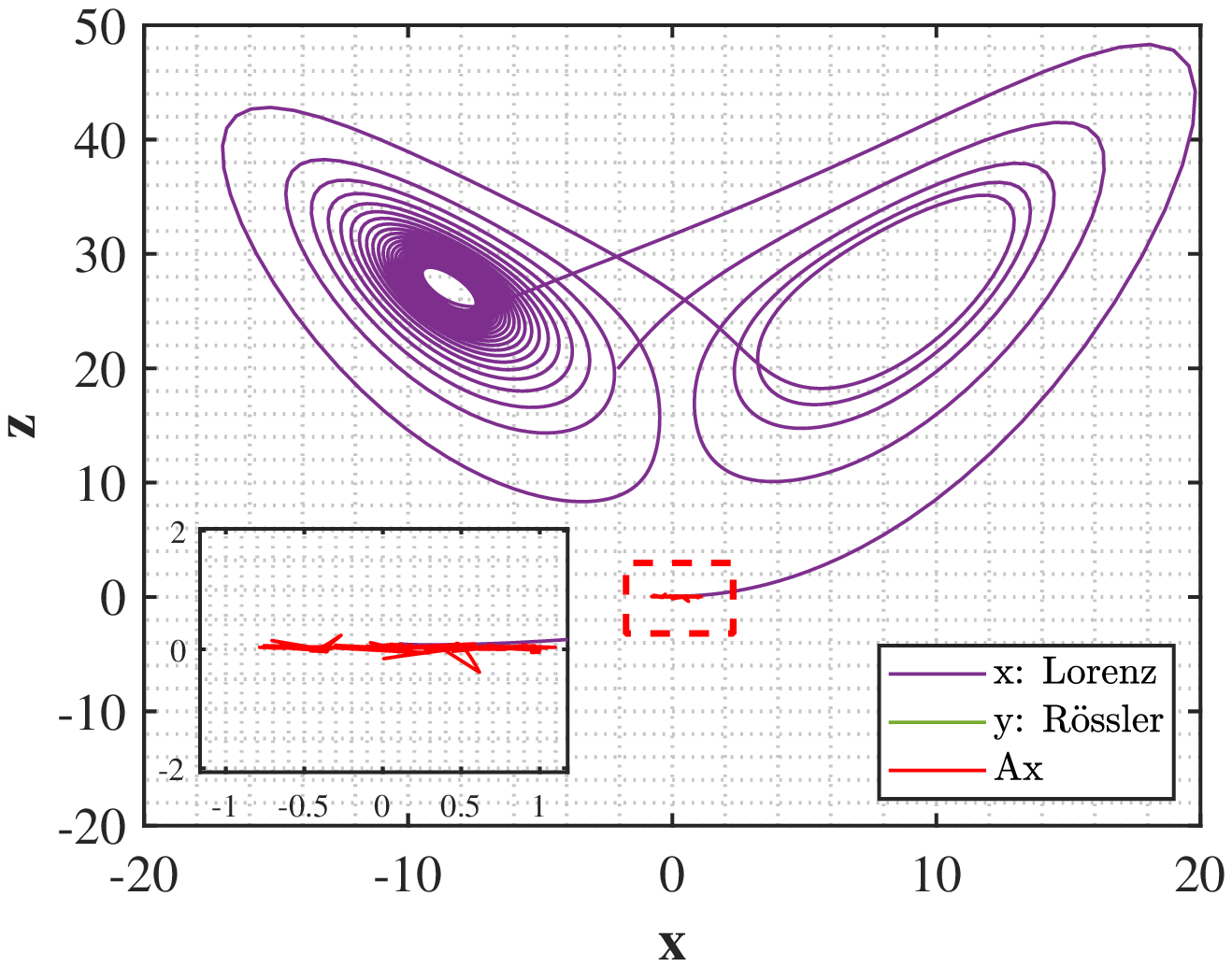}
\hfil
\includegraphics[scale=.28]{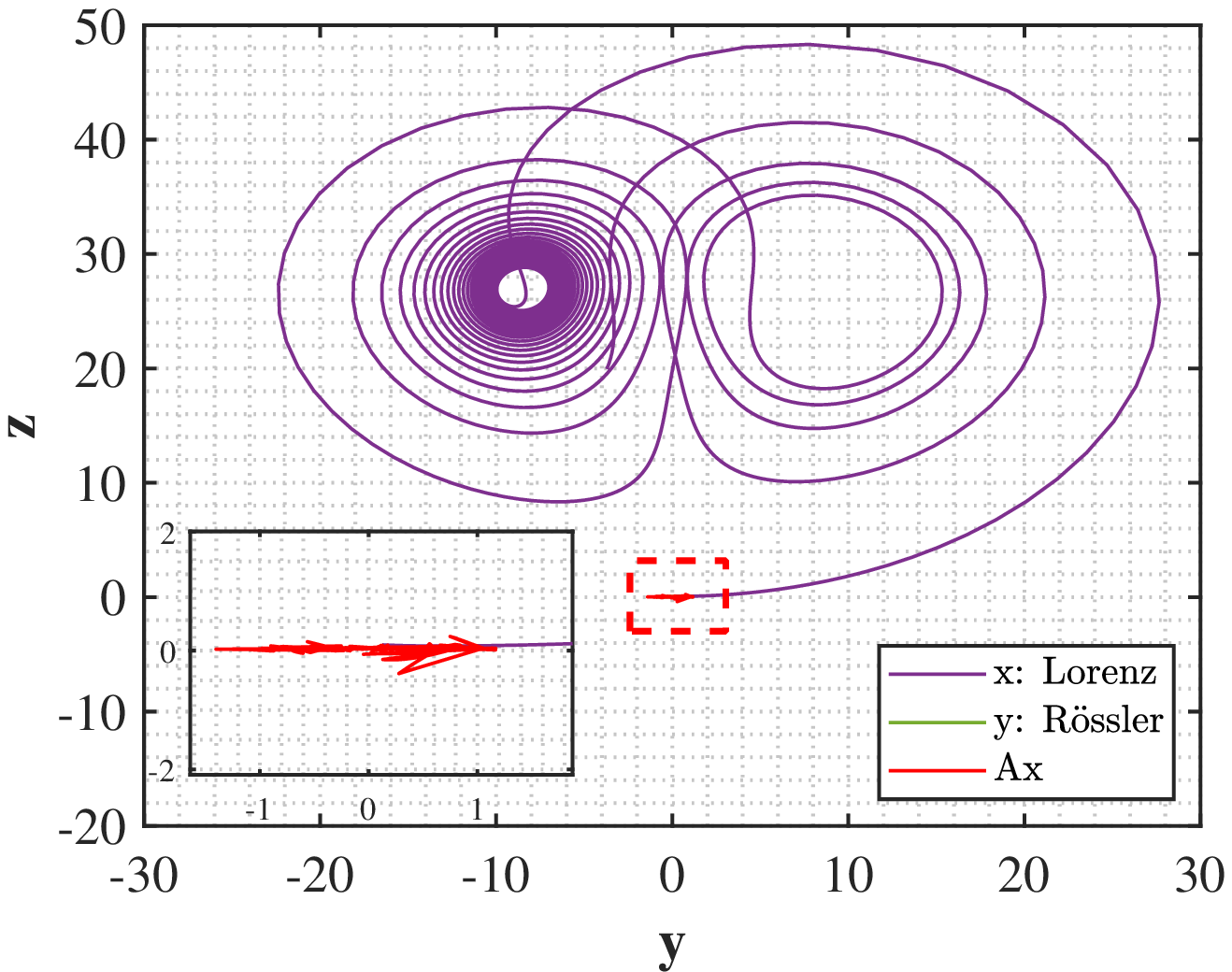}
\label{Fig.10a}}
\hfil
\subfloat[Based on L2-norm penalty.]{\includegraphics[scale=.28]{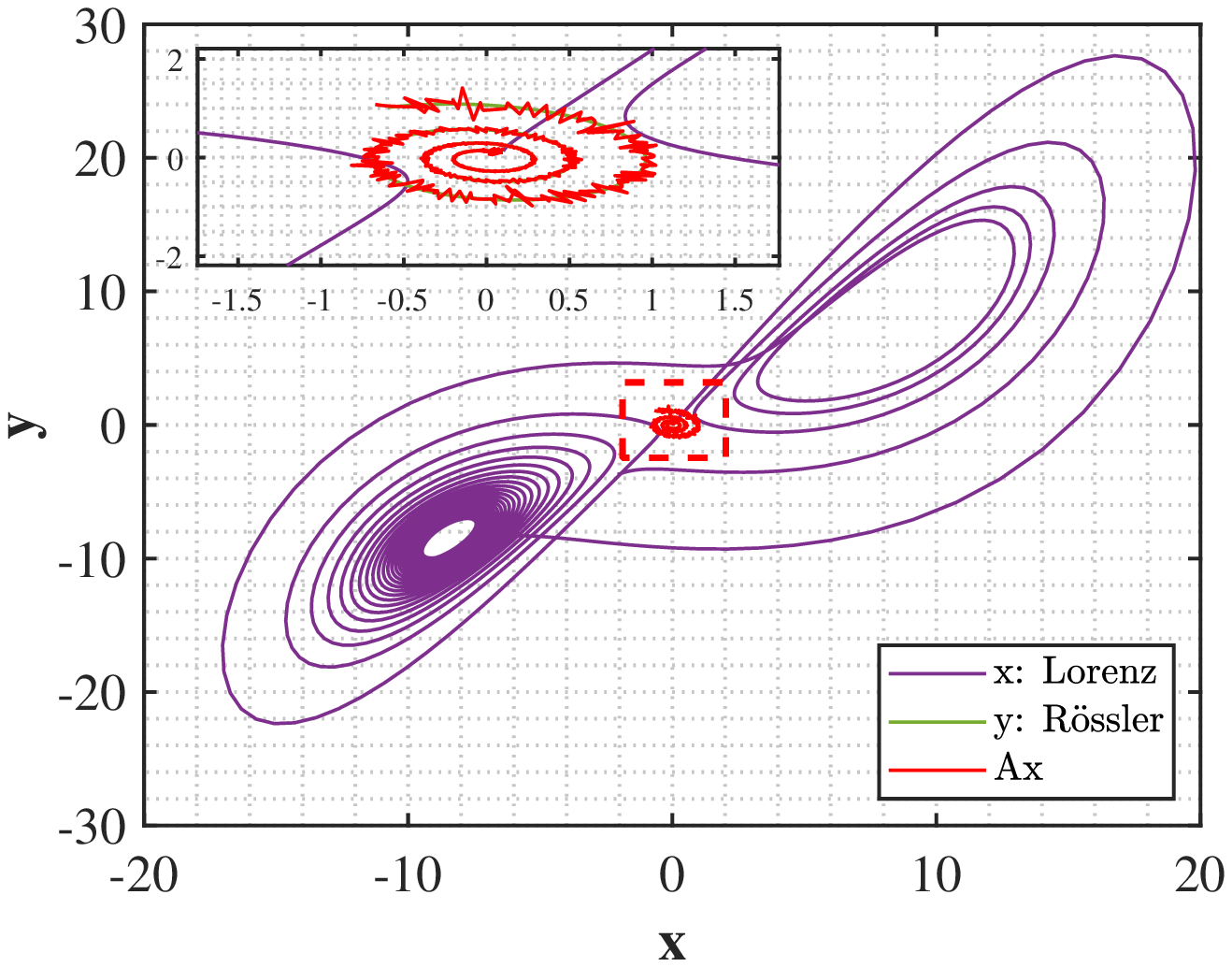}
\hfil
\includegraphics[scale=.28]{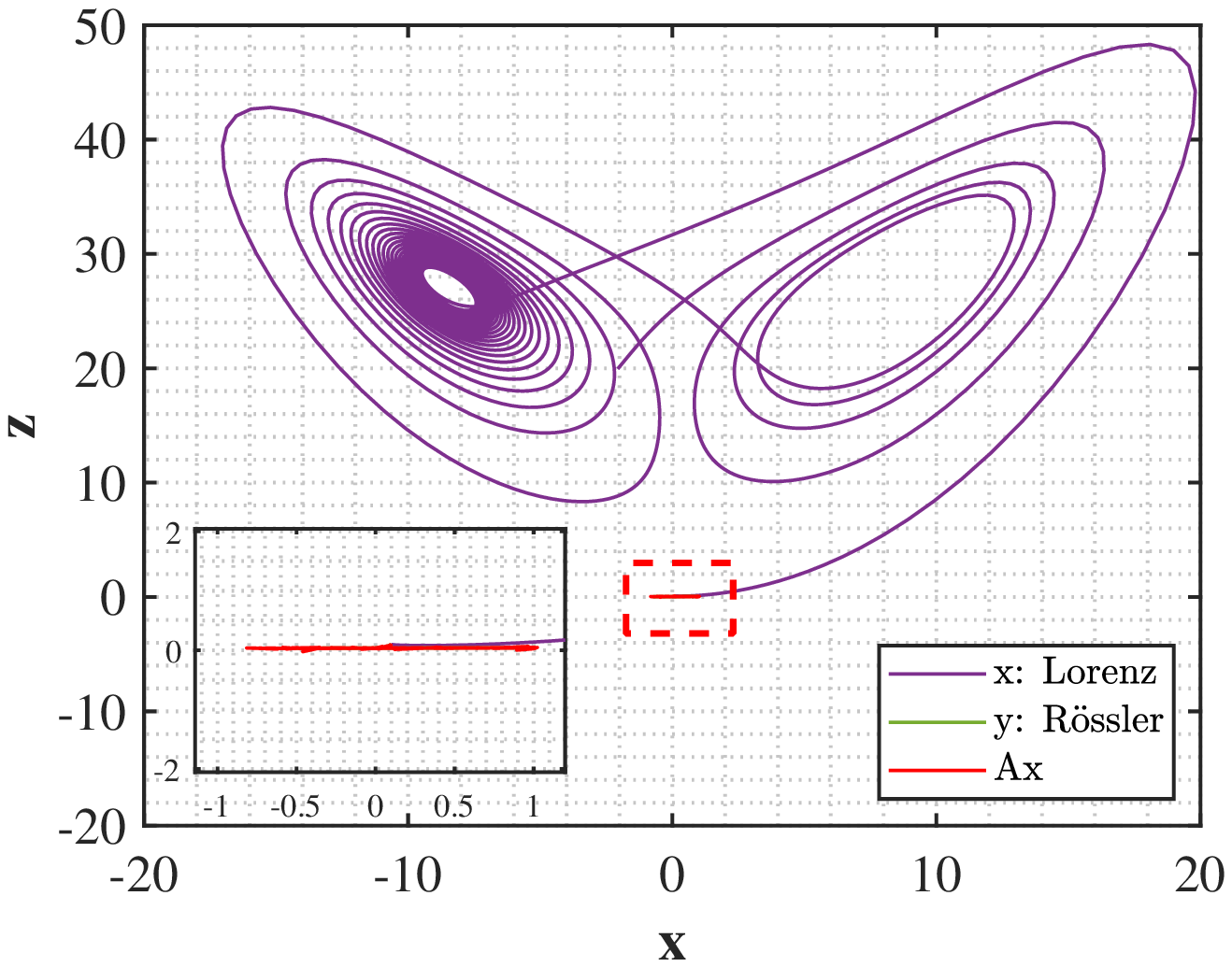}
\hfil
\includegraphics[scale=.28]{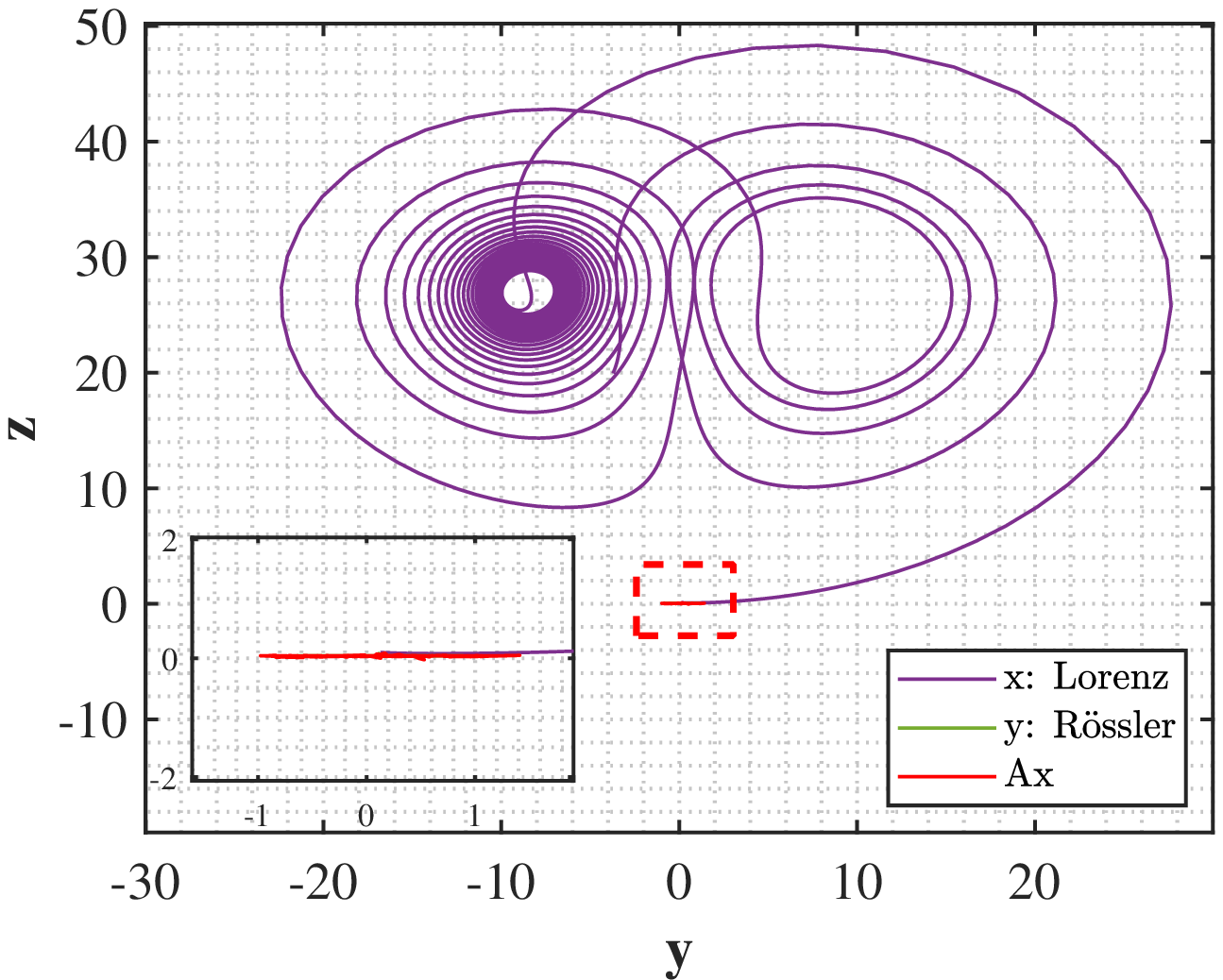}
\label{Fig.10b}}
\caption{Two dimensional plans of Example~4.2.}
\label{Fig.10}
\end{figure}

Figs. \ref{Fig.9} and \ref{Fig.10} shown the stereograms and plans of Lorenz attractor, R$\rm\ddot{o}$ssler attractor and the trajectory acted by optimal similarity transformation matrix for each stage when $N$$=$$2000$. Even if overlapped trajectories marked by dotted box are enlarged, we could still observe that the orbits of R$\rm\ddot{o}$ssler system almost coincides with that of Lorenz system acted by optimal similarity transformation matrix.


Lorenz attractor with butterfly shape and R$\rm \ddot{o}ssler$ attractor with spiral shape that appear to be different geometrically, become remarkably similar within the proper precision under the action of Pontryagin's maximum principle using optimal principle (\ref{22}), which is greatly an amazing finding.

Surprisingly, the numerical results indicated in Fig. 11 show that only six results are less than 0.95 without regularization term. It should be point out that when we carry out tests to find optimal similarity transformation matrix based on L2-norm penalty (\ref{40}), all values of similarity degree are greater than 0.98.
\begin{figure}[htbp]
\centering
\subfloat[Without~penalty.]{\includegraphics[scale=.4]{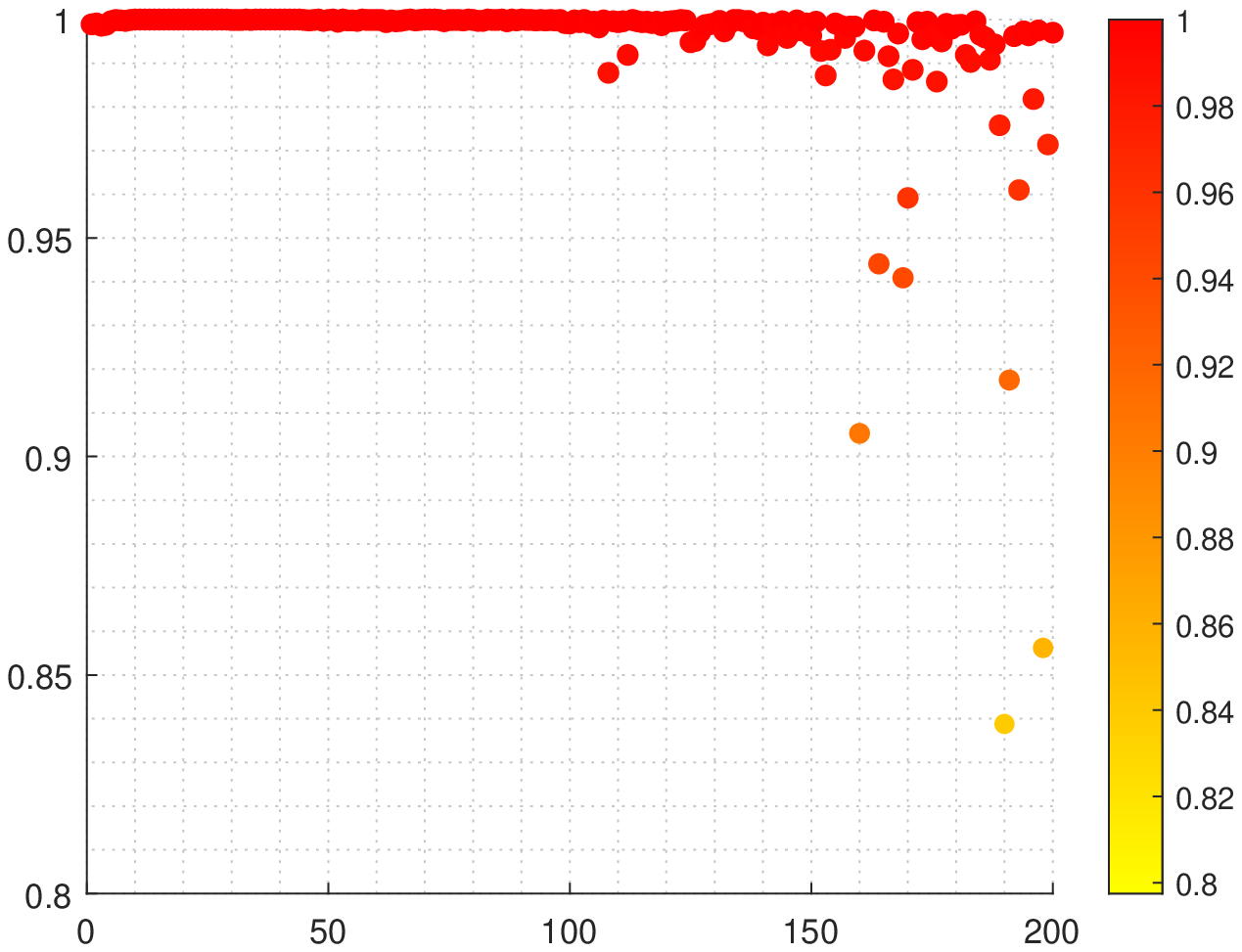}
\label{Fig.11a}}
\hfil
\subfloat[Based on L2-norm penalty.]{\includegraphics[scale=.4]{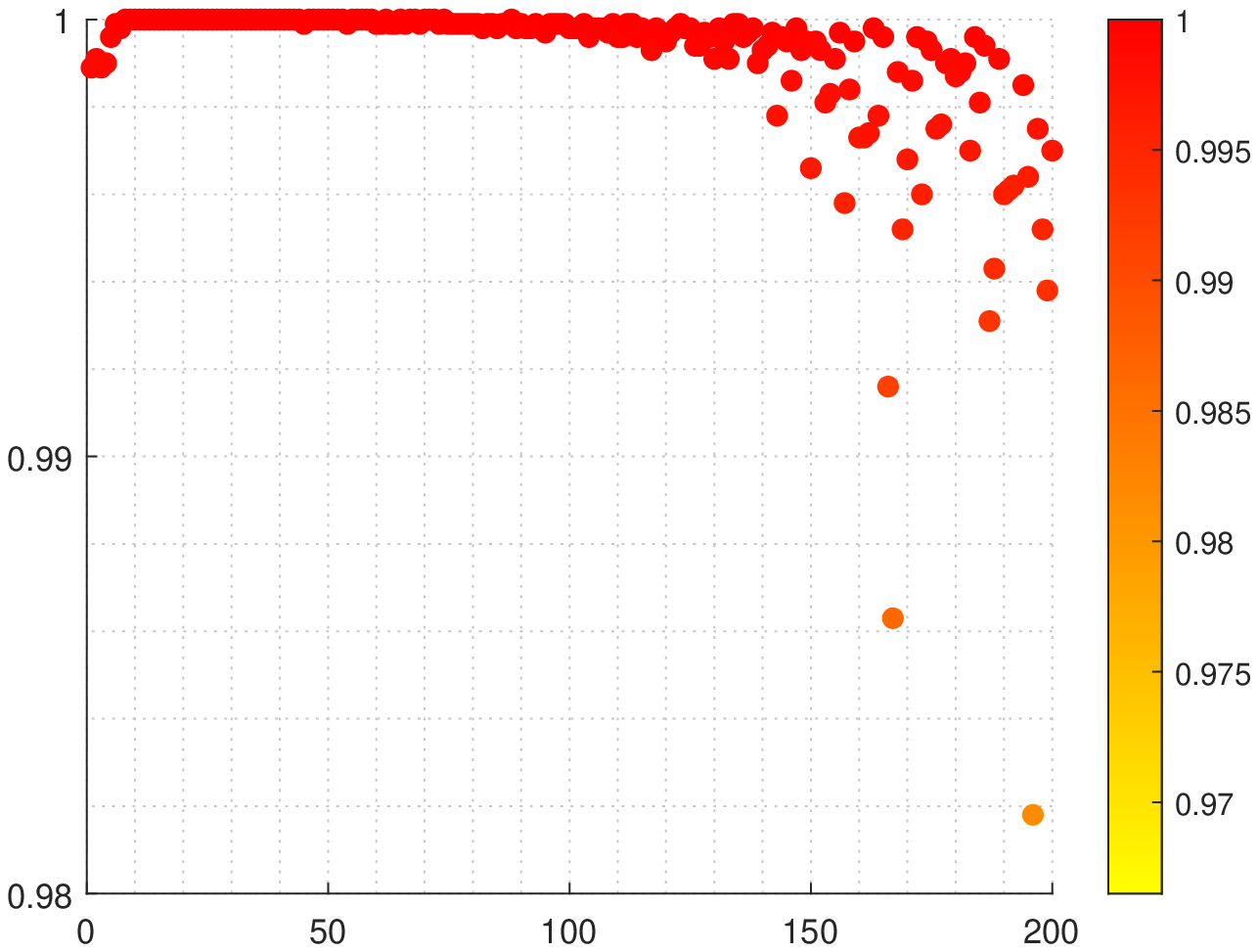}
\label{Fig.11b}}
\caption{Similarity degree of Example~4.2.}
\label{Fig.11}
\end{figure}

To summarize, even without the use of L2-norm regularization, we can still get very satisfactory results to some extent.

\subsection{Bellman's dynamic programming}
Dynamic programming was studied by Bellman to deal with situations where the best decisions are made in stages \cite{Bellman1957}. Whatever the initial state and initial decision are, the decisions that will follow must also constitute an optimal policy for the remaining problems, when the stage and state formed by the first step decision are considered as initial conditions. Applying Bellman's dynamic programming in terms of the proposed optimal principle, we analyze the similarity of orbits between two chaotic attractors, see the following three examples.

{\bf Example~4.3.}
Similarity of orbits between Lorenz attractor and Chen attractor.

Let $\{x_k\}$ and $\{y_k\}$ be the numerical solutions obtained from Lorenz and Chen systems for 2000 time steps. Instead of solving 2000 steps one at a time, we consider multi-stage decision making process by breaking the complex problems into ten simple subproblems denote as $N_1$$-$$N_{200}$ with 10 steps for each.

More precisely, for stage $N_1$, the initial conditions of two systems are determined as
\begin{equation*}
x_0=(0.1,0.1,0.1)^{\rm T}~\mbox{and}~y_0=Ax_0
\end{equation*}
with $A$ unknown. Make use of Runge-Kutta method, we calculate the states of this stage including the values of $x_1$$-$$x_{10}$ and $y_1$$-$$y_{10}$ whose components are represented by expressions containing elements of $A$. By solving a nonlinear equations formulated by (\ref{22}), we obtain an approximate solution with a high precision.

For the following stages, we do the same actions and compare the current similarity transformation matrix $\bar{A}$ with the one obtained by previous stage denoted as $\bar{\bar{A}}$ subject to similarity degree defined in (\ref{39}), then let
\begin{equation*}
A=\max\limits_{\bar{A},\bar{\bar{A}}} \{\rho(\bar{A}), \rho(\bar{\bar{A}})\}
\end{equation*}
be the optimal similarity transformation matrix of this stage. We obtain the approximate solution which makes similarity degree reach 1.0000 for each stage.

\begin{figure}[htbp]
\centering
\subfloat[]{\includegraphics[scale=.27]{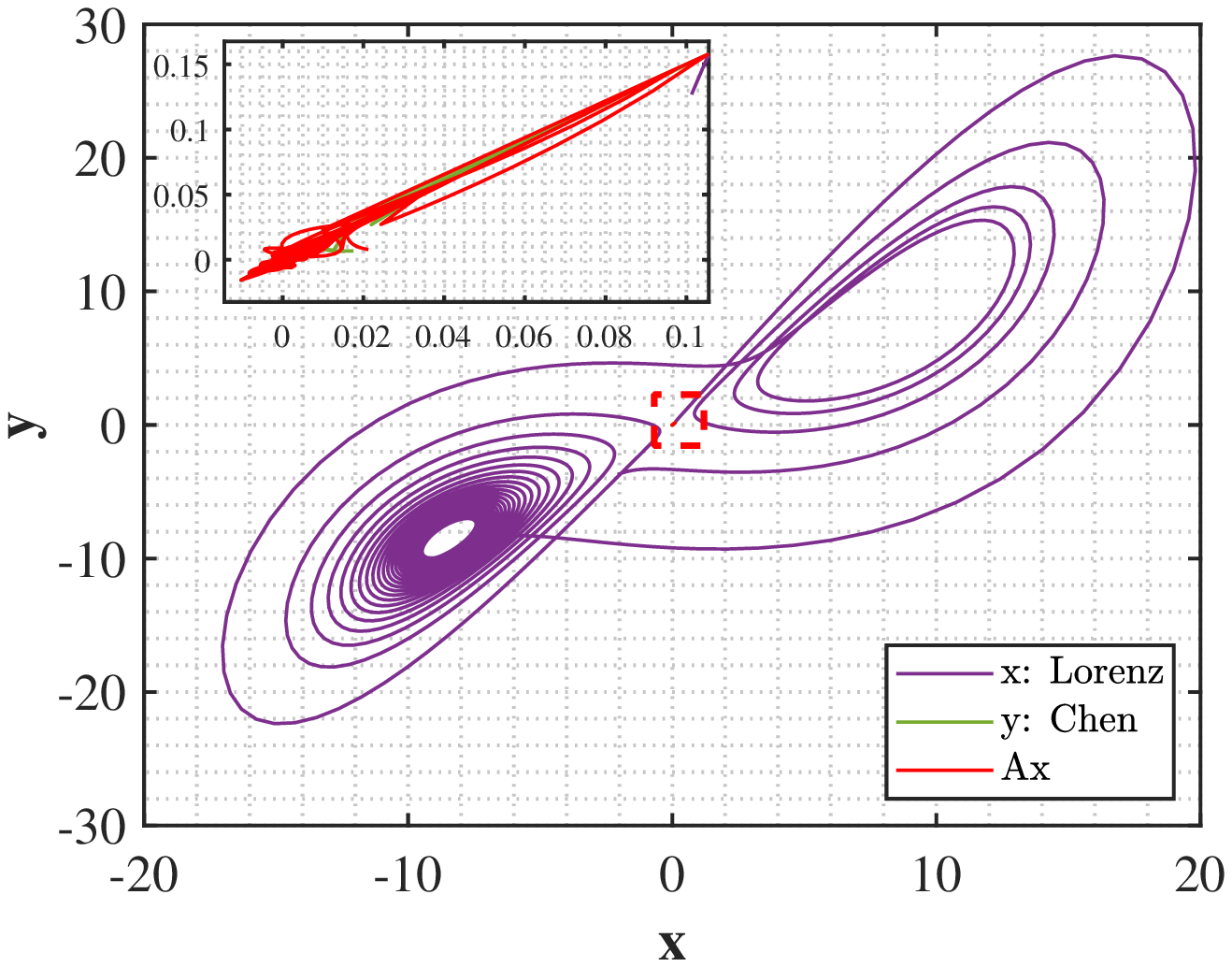}
\label{Fig.12a}}
\hfil
\subfloat[]{\includegraphics[scale=.27]{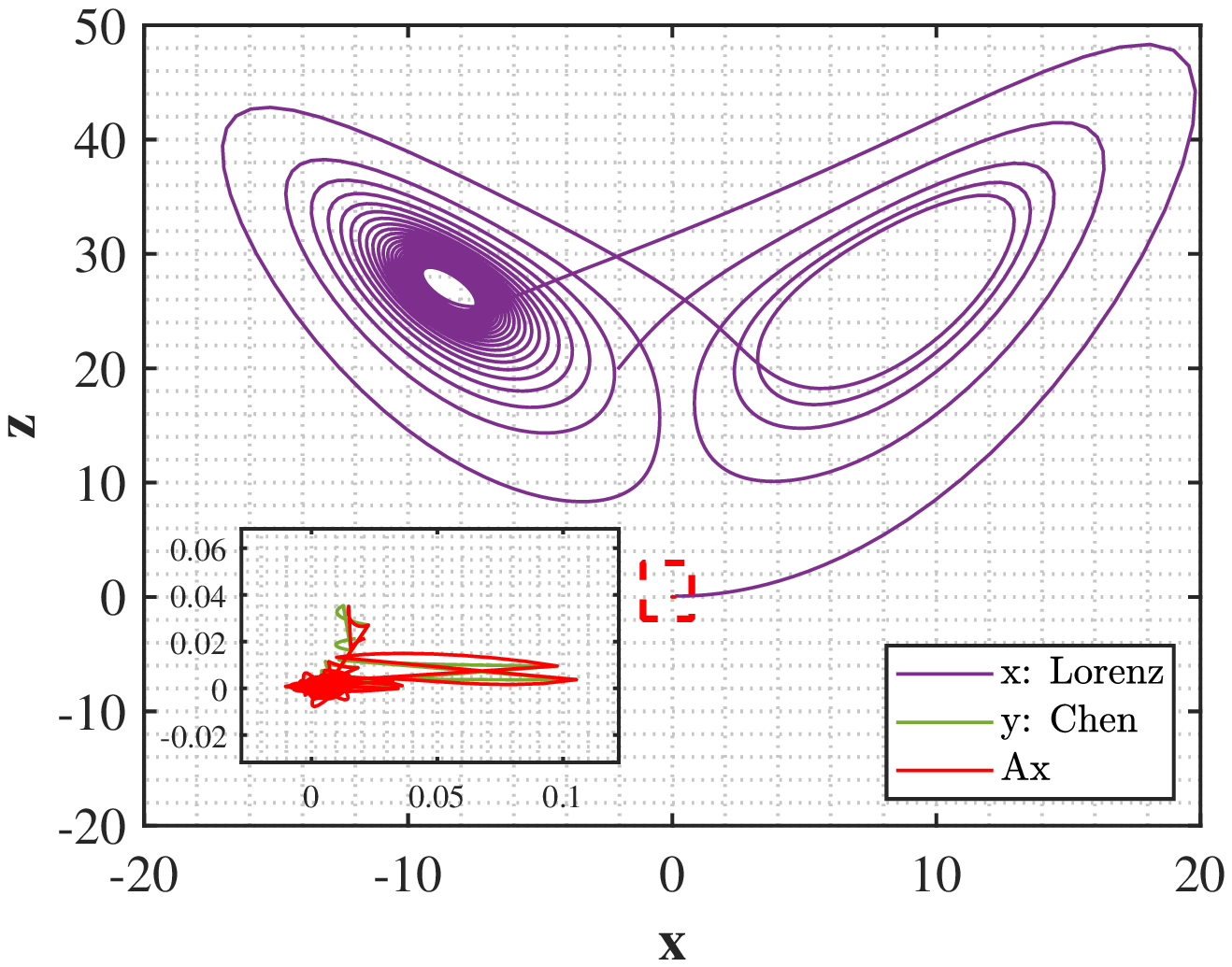}
\label{Fig.12b}}
\hfil
\subfloat[]{\includegraphics[scale=.27]{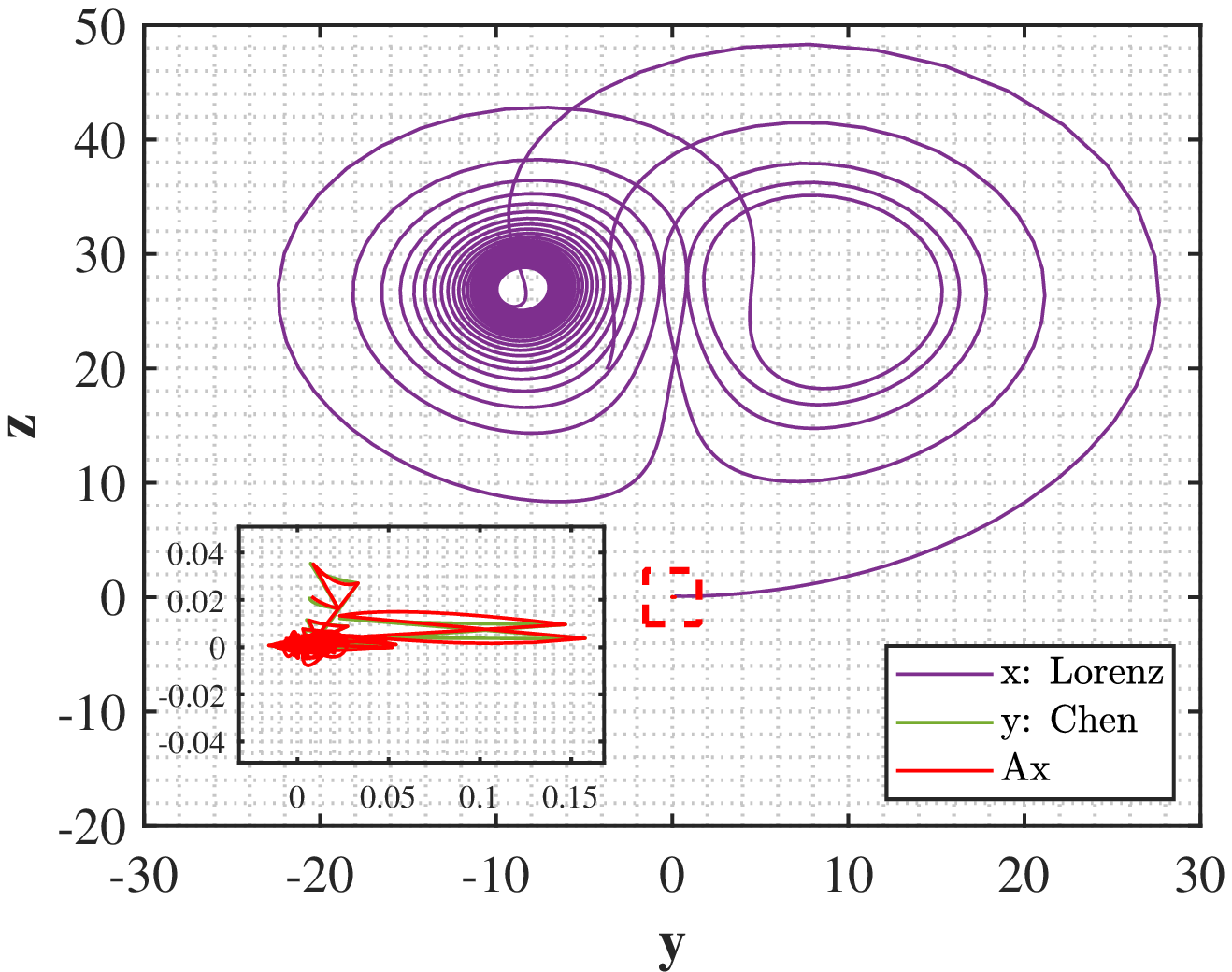}
\label{Fig.12c}}
\caption{Two dimensional plans of Example~4.3.}
\label{Fig.12}
\end{figure}

\begin{figure}[htbp]
\centering
\begin{minipage}{0.45\linewidth}
\includegraphics[scale=.4]{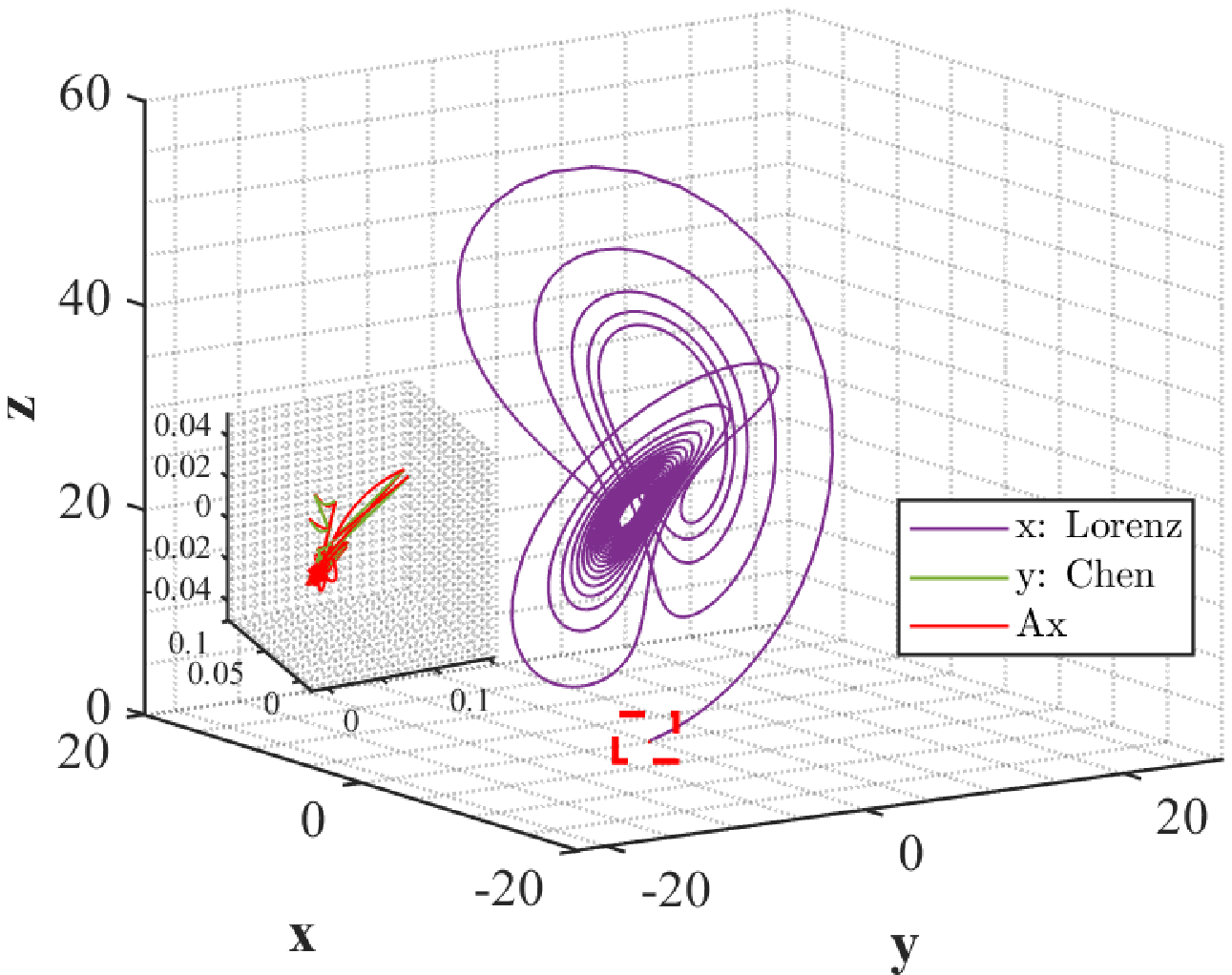}
\caption{Three dimensional stereograms of Example~4.3.}
\label{Fig.13}
\end{minipage}
\begin{minipage}{0.45\linewidth}
\includegraphics[scale=.4]{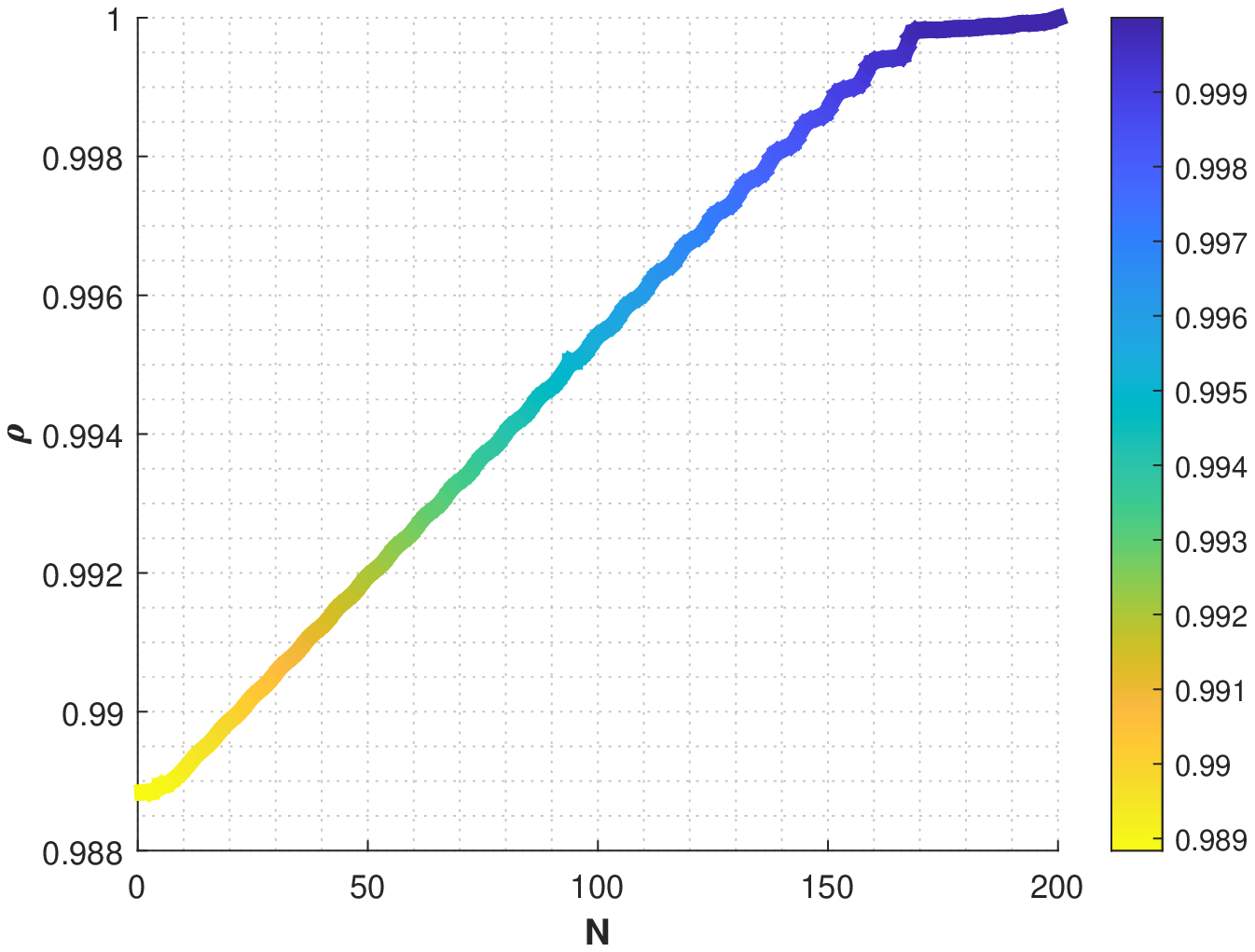}
\caption{Change in similarity degree of Example~4.3.}
\label{Fig.14}
\end{minipage}
\end{figure}

As shown in Figs. \ref{Fig.12}-\ref{Fig.13}, Lorenz and Chen systems with different orbits can become distinct similar through the adjustment of similarity transformation matrix derived by the proposed optimal principle.

Now we elaborate the advantage of multi-stage dynamic programming with the help of numerical results. For $N$$=$$2000$, the optimal similarity transformation matrix $A_0$ can be found according to optimal principle (\ref{22}), which makes similarity degree reach 0.988837. When only the optimal similarity transformation matrix at stage $N_1$ is taken and $A_0$ is still employed in other stages, similarity degree increases to 0.988839. If we adopt the corresponding optimal similarity transformation matrix in both $N_1$ and $N_2$, the similarity degree rises to 0.9888435, and so on. The final similarity degree of solutions between Lorenz attractor and Chen attractor in this example can reach 0.999995 and each stage is optimal at this point, which meets Bellman's principle of optimality. To get a better view of the change in similarity degree, the numerical results of each stage are depicted in Fig. \ref{Fig.14}. We observe that as the number of the optimal similarity transformation matrix in corresponding stage increases, similarity degree is increase progressively.

{\bf Example~4.4.}
Similarity of orbits between Lorenz attractor and L$\rm\ddot{u}$ attractor.

Let $\{x_k\}$ and $\{y_k\}$ be the numerical solutions got from Lorenz system and L$\rm\ddot{u}$ system with $u=0$ for 2000 time steps. Similar to the multi-stage decision in Example 4.3, we also divide the steps into 200 stage. For each stage, only the initial state of sequence $\{x_k\}$ is known. The values of similarity degree can reach 1.0000 for all stages, implying the effectiveness of the optimal similarity transformation matrix.

\begin{figure}[ht]
\centering
\subfloat[]{\includegraphics[scale=.28]{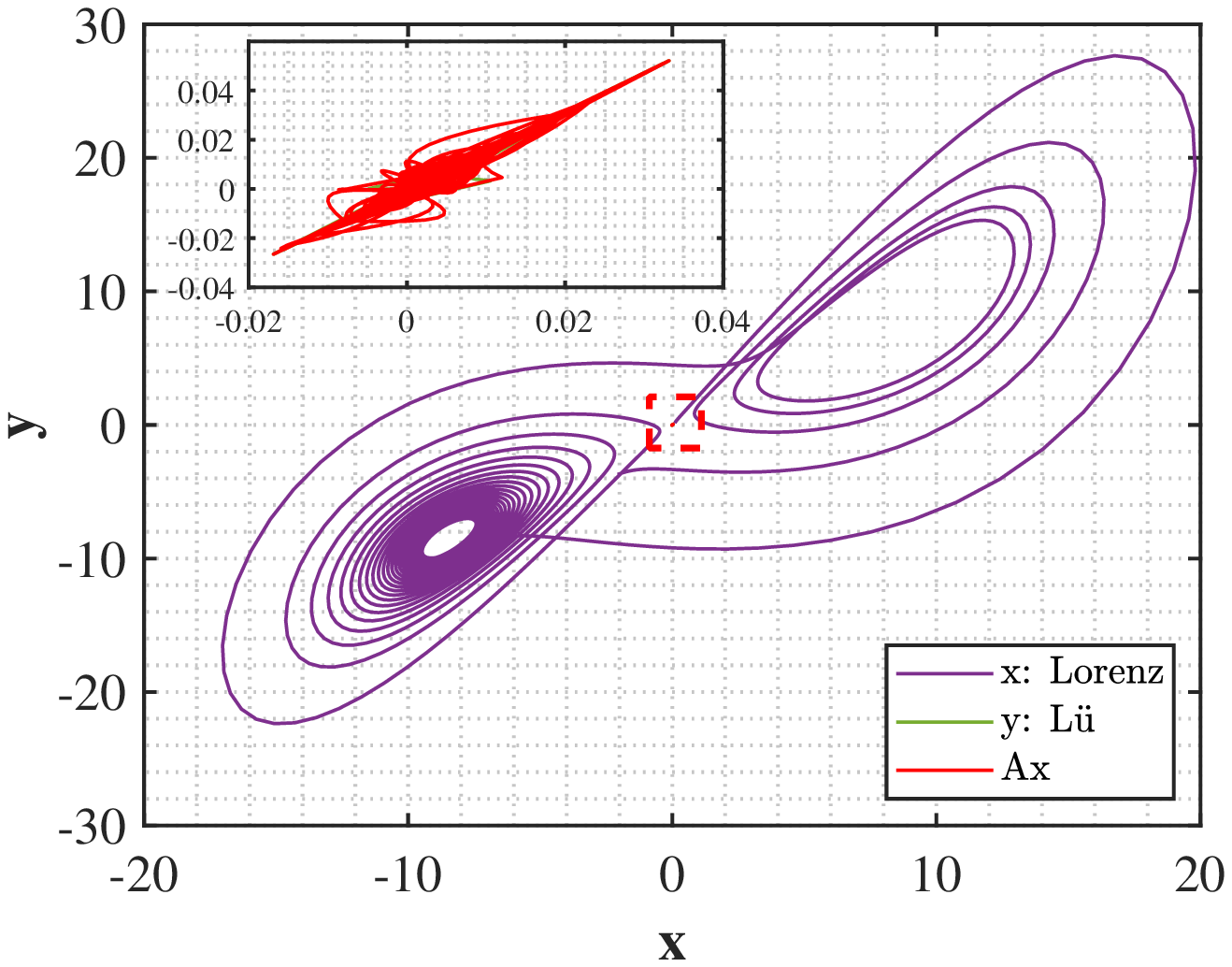}
\label{Fig.15a}}
\hfil
\subfloat[]{\includegraphics[scale=.28]{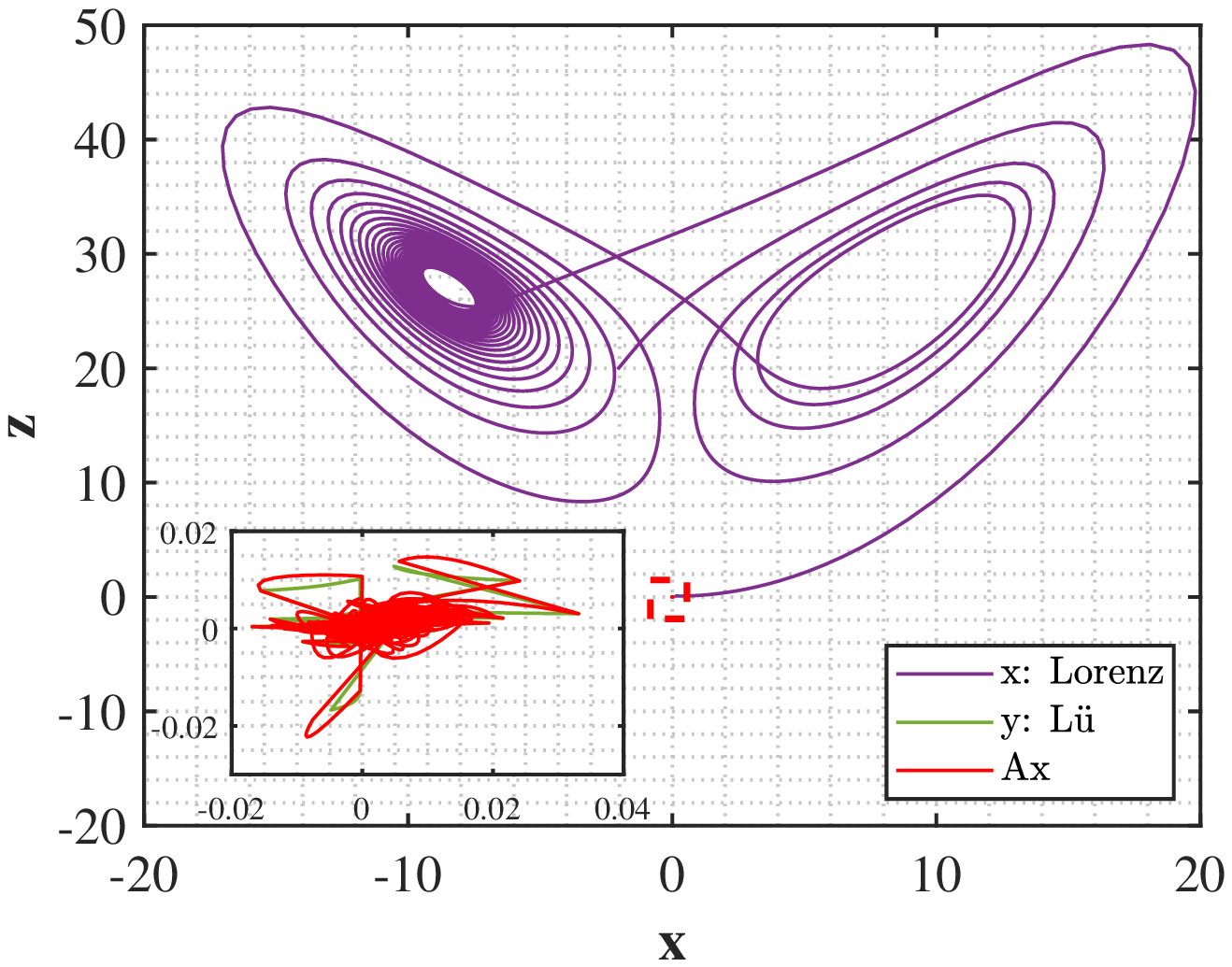}
\label{Fig.15b}}
\hfil
\subfloat[]{\includegraphics[scale=.28]{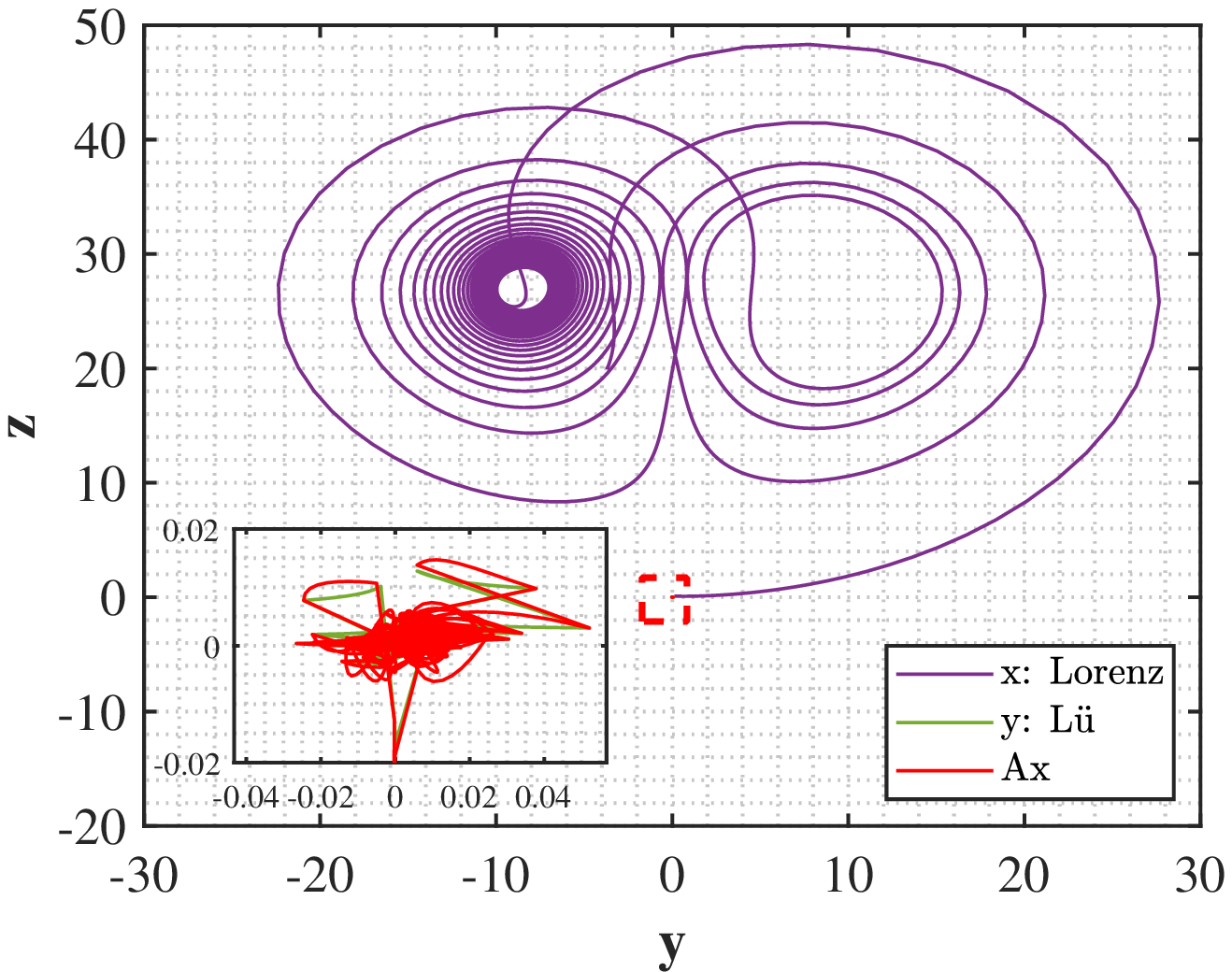}
\label{Fig.15c}}
\caption{Two dimensional plans of Example~4.4.}
\label{Fig.15}
\end{figure}

\begin{figure}[ht]
\centering
\begin{minipage}{0.45\linewidth}
\includegraphics[scale=.4]{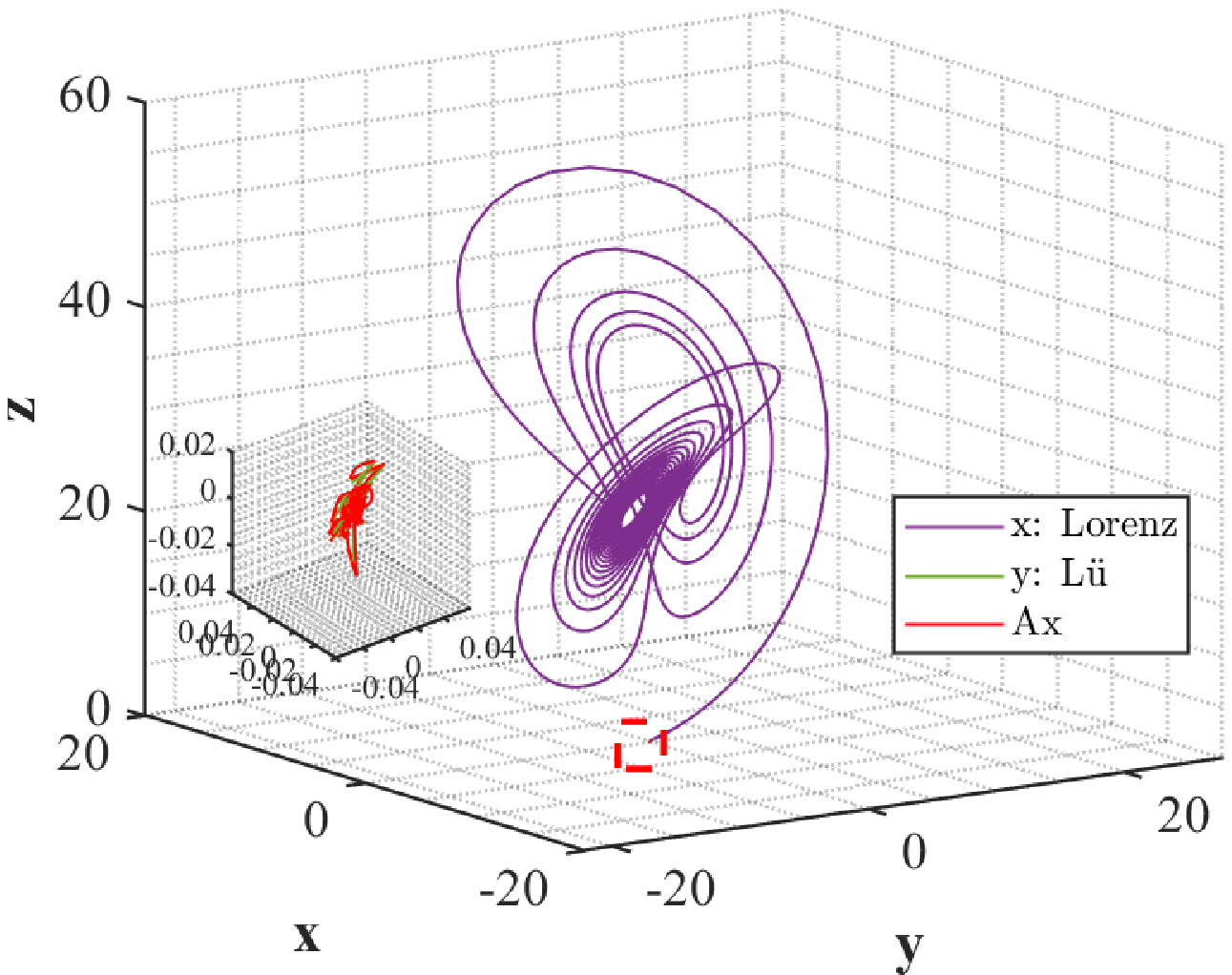}
\caption{Three dimensional stereograms of Example~4.4.}
\label{Fig.16}
\end{minipage}
\begin{minipage}{0.45\linewidth}
\includegraphics[scale=.4]{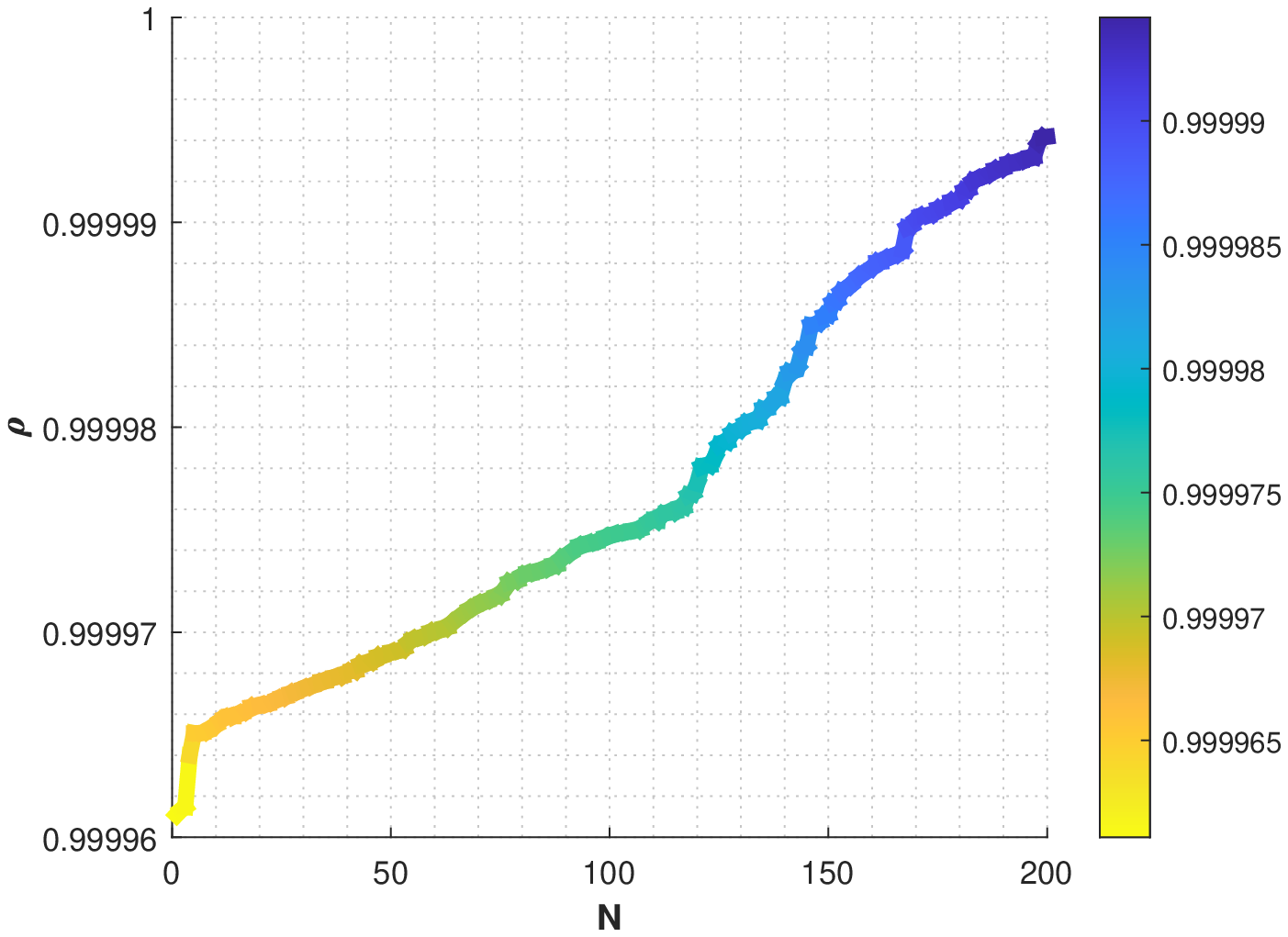}
\caption{Change in similarity degree of Example~4.4.}
\label{Fig.17}
\end{minipage}
\end{figure}

We are surprised by the effectiveness of similarity transformation matrix formulated by the proposed optimal principle, as shown in Figs. \ref{Fig.15}-\ref{Fig.16}. Even if we enlarge the trajectories represented in dotted box to the coordinate diagram with small  horizontal and vertical coordinates, the orbits of L$\rm\ddot{u}$ attractor and Lorenz attractor acted by optimal similarity transformation matrix can still coincide almost exactly. The change of similarity degree gradually increase from $0.999609$ to $0.999994$ with the increase of the number of optimal similarity transformation matrix, satisfying Bellman's principle of optimality, see Fig. \ref{Fig.17}.

{\bf Example~4.5.}
Similarity of orbits between Hybrid attractor and L$\rm\ddot{u}$ attractor.

The last example concerns hybrid Lorenz-Chua chaotic system formed by Lorenz attractor and Chua's circuit using the homotopy approach (\ref{23}), modelling below
\begin{equation}\label{41}
\begin{array}{lll}
\dot{x}=\lambda(-\sigma x+\sigma y)+(1-\lambda)\cdot\alpha[y-x-f(x)],\\
\dot{y}=\lambda(-xz+rz-y)+(1-\lambda)(x-y+z), \\
\dot{z}=\lambda(xy-bz)+(1-\lambda)(-\beta y),
\end{array}
\end{equation}
where the piece-linear function $f(x)$ is defined in (\ref{3}), and all the same parameters as in (\ref{1})-(\ref{3}). The study on hybrid attractor is more challenging due to its more complex topologies and dynamics.

Different dynamical behaviors in L$\rm\ddot{u}$ attractor's controlled system (\ref{6}) can be generated by varying the parameter $u$. For the parameters $u$$=$$-1$, $u$$=$$8$, $u$$=$$-12$ and $u$$=$$12$ that produce complete attractor, partial attractor, left-attractor and right-attractor, we simulate the similarity between orbits of hybrid Lorenz-Chua system and L$\rm\ddot{u}$ attractor respectively.

\begin{figure}[htbp]
\centering
\subfloat[$\rm u$=-1.]{\includegraphics[scale=.28]{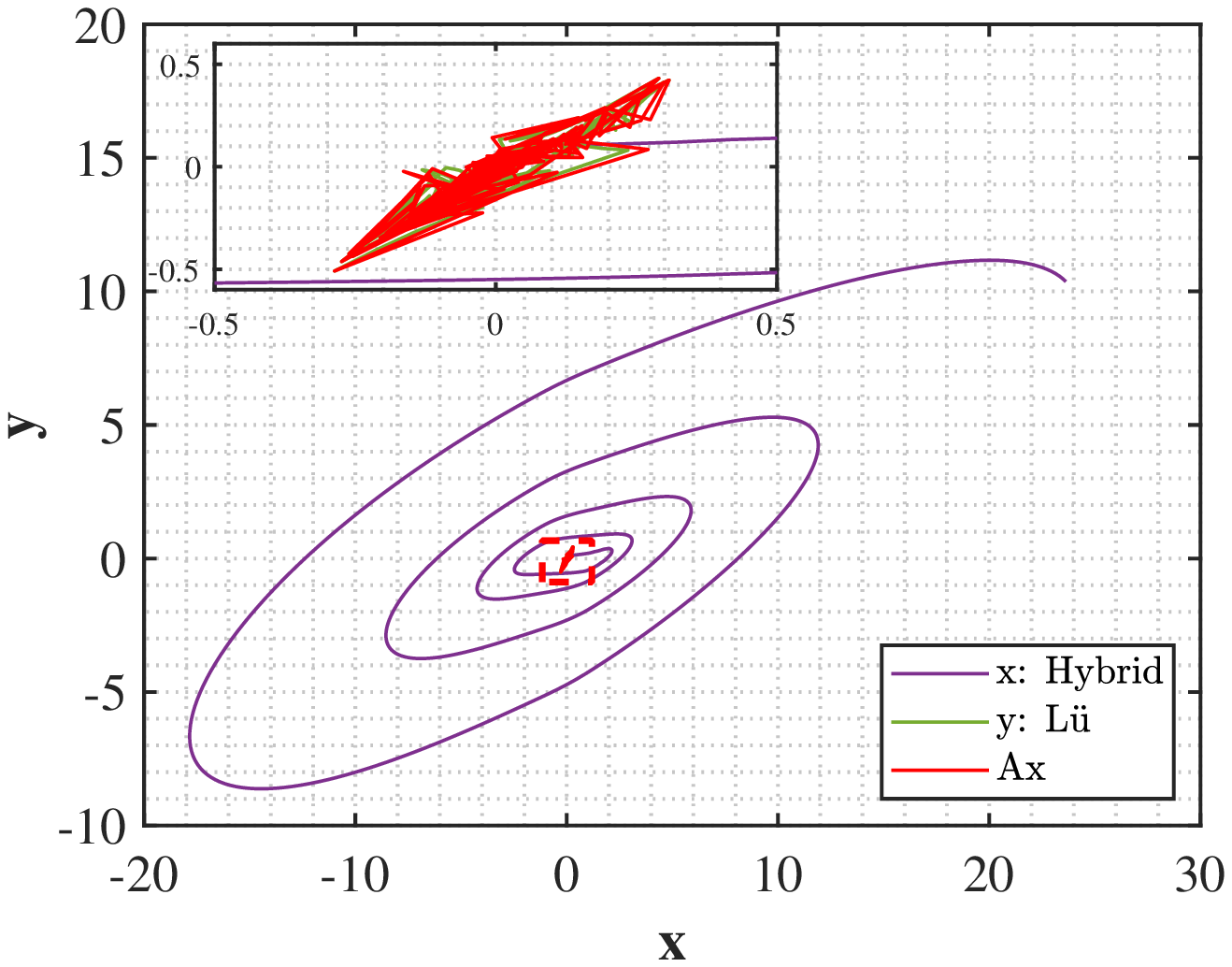}
\hfil
\includegraphics[scale=.28]{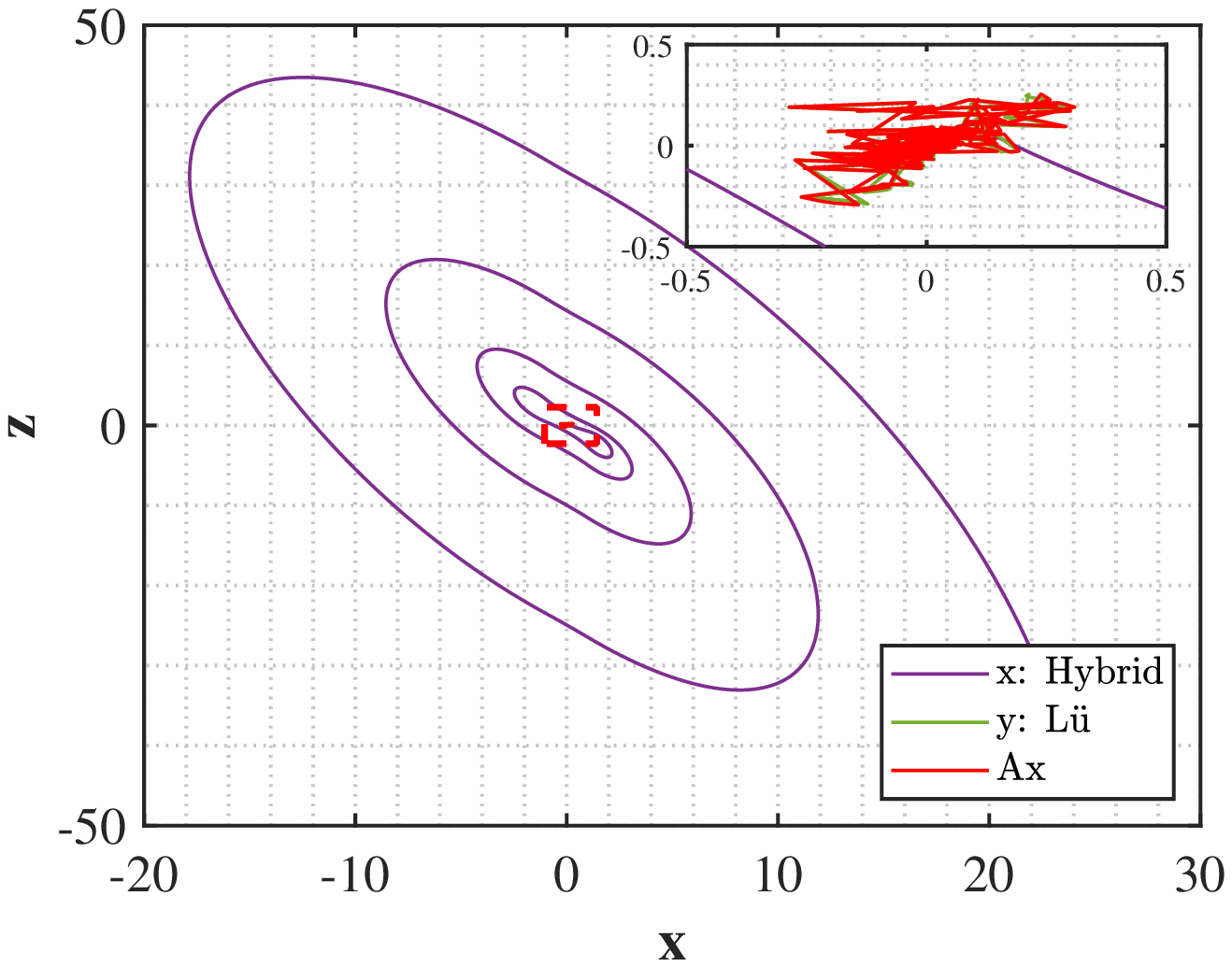}
\hfil
\includegraphics[scale=.28]{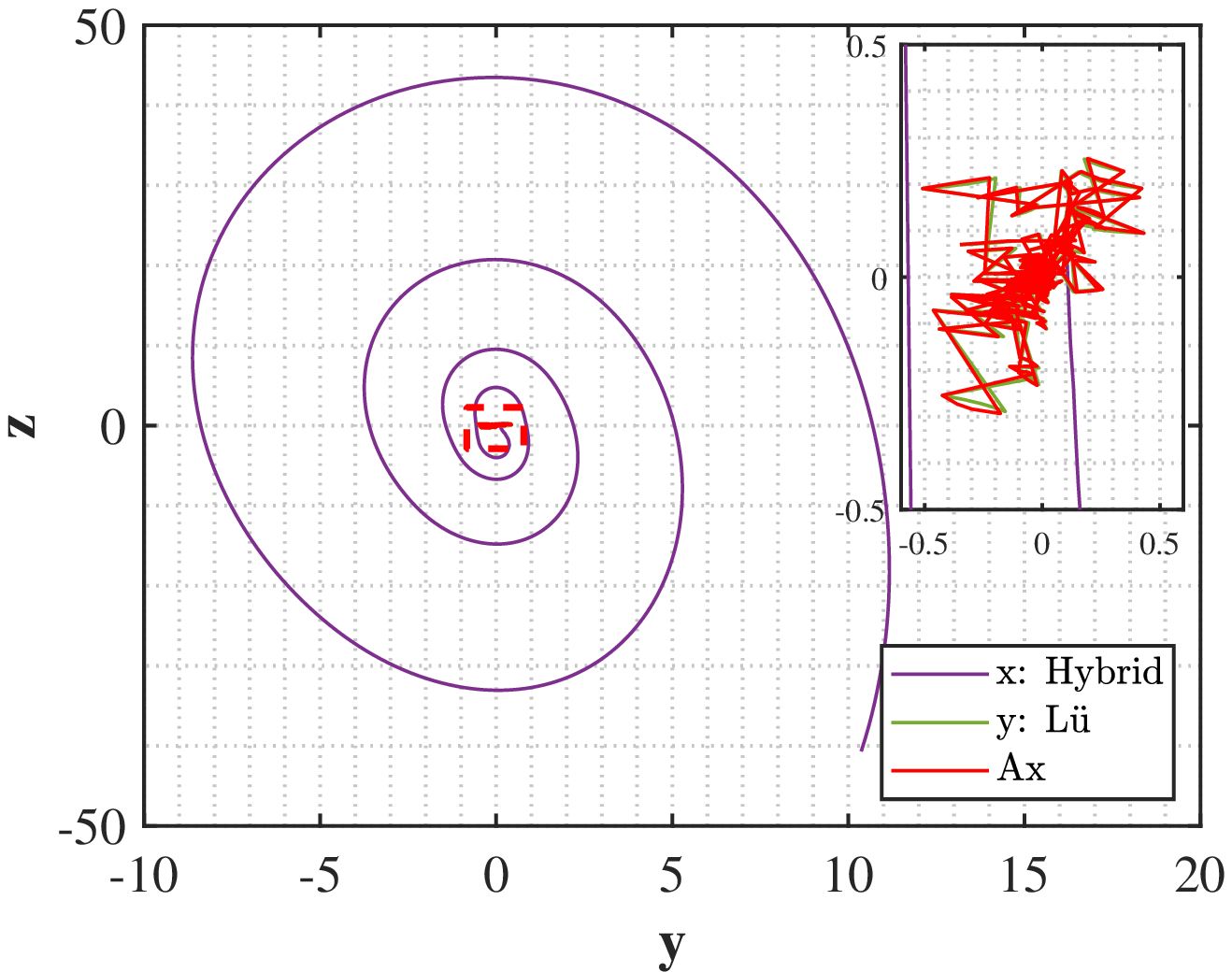}
\label{Fig.18a}}
\hfil
\subfloat[$\rm u$=8.]{\includegraphics[scale=.28]{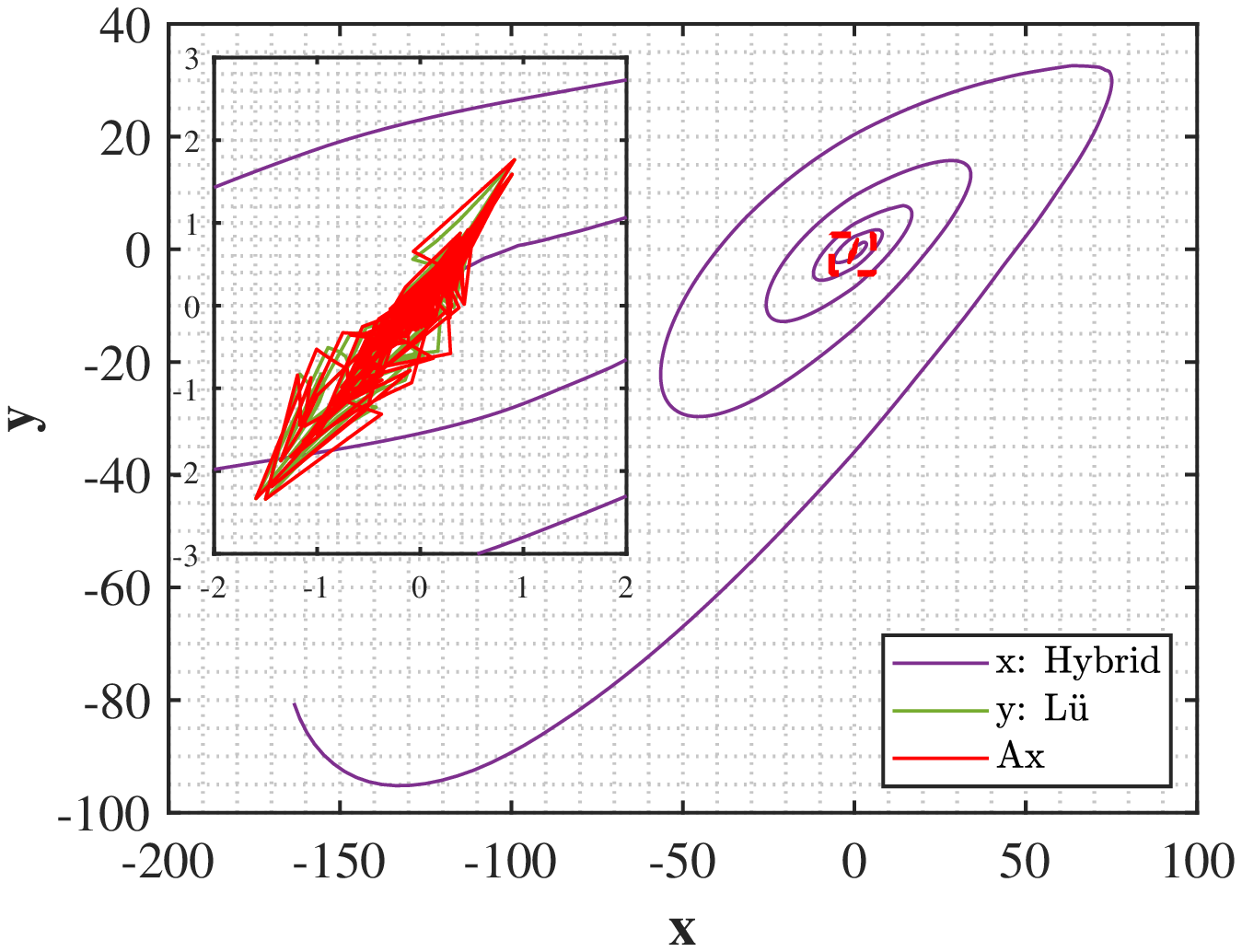}
\hfil
\includegraphics[scale=.28]{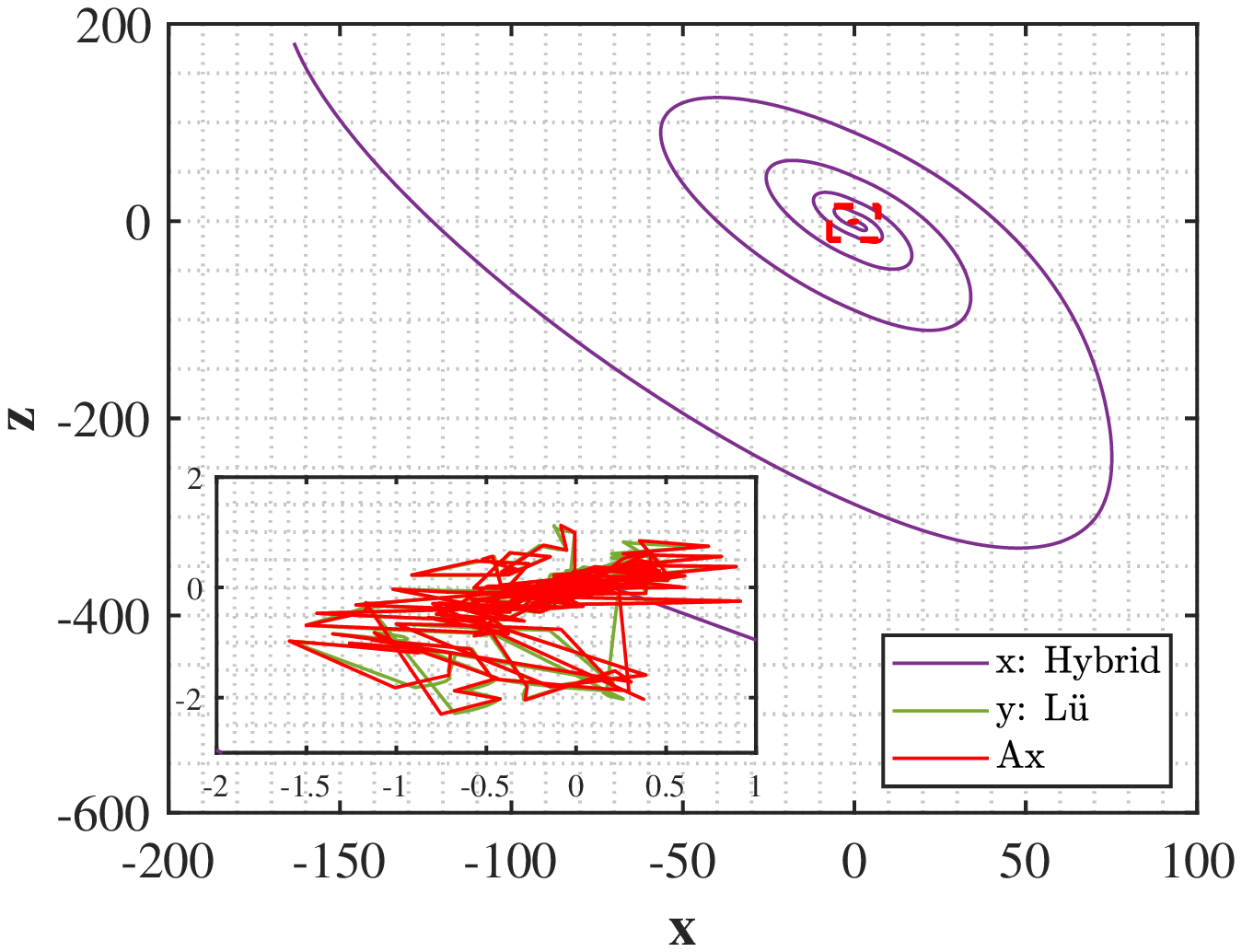}
\hfil
\includegraphics[scale=.28]{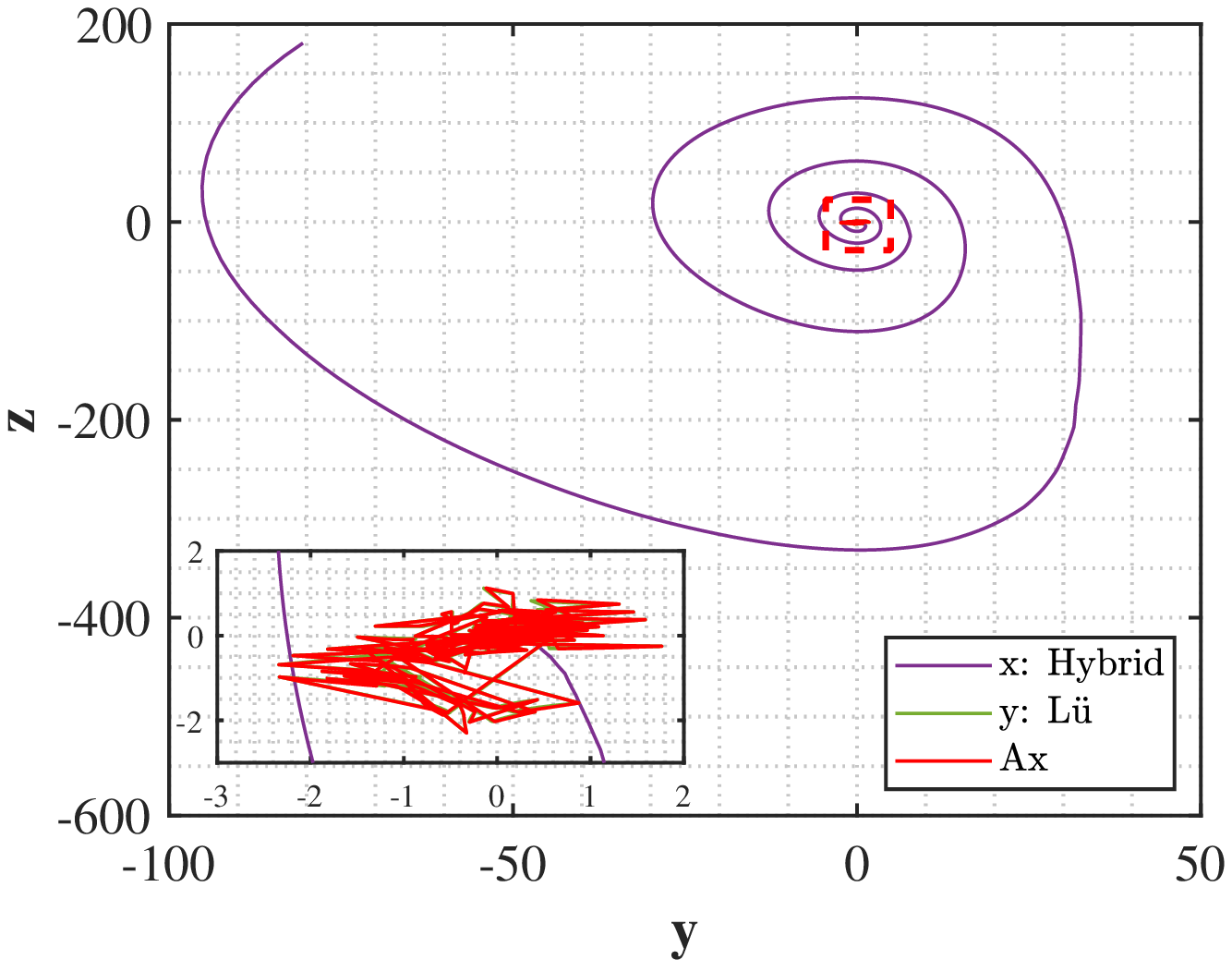}
\label{Fig.18b}}
\hfil
\subfloat[$\rm u$=12.]{\includegraphics[scale=.28]{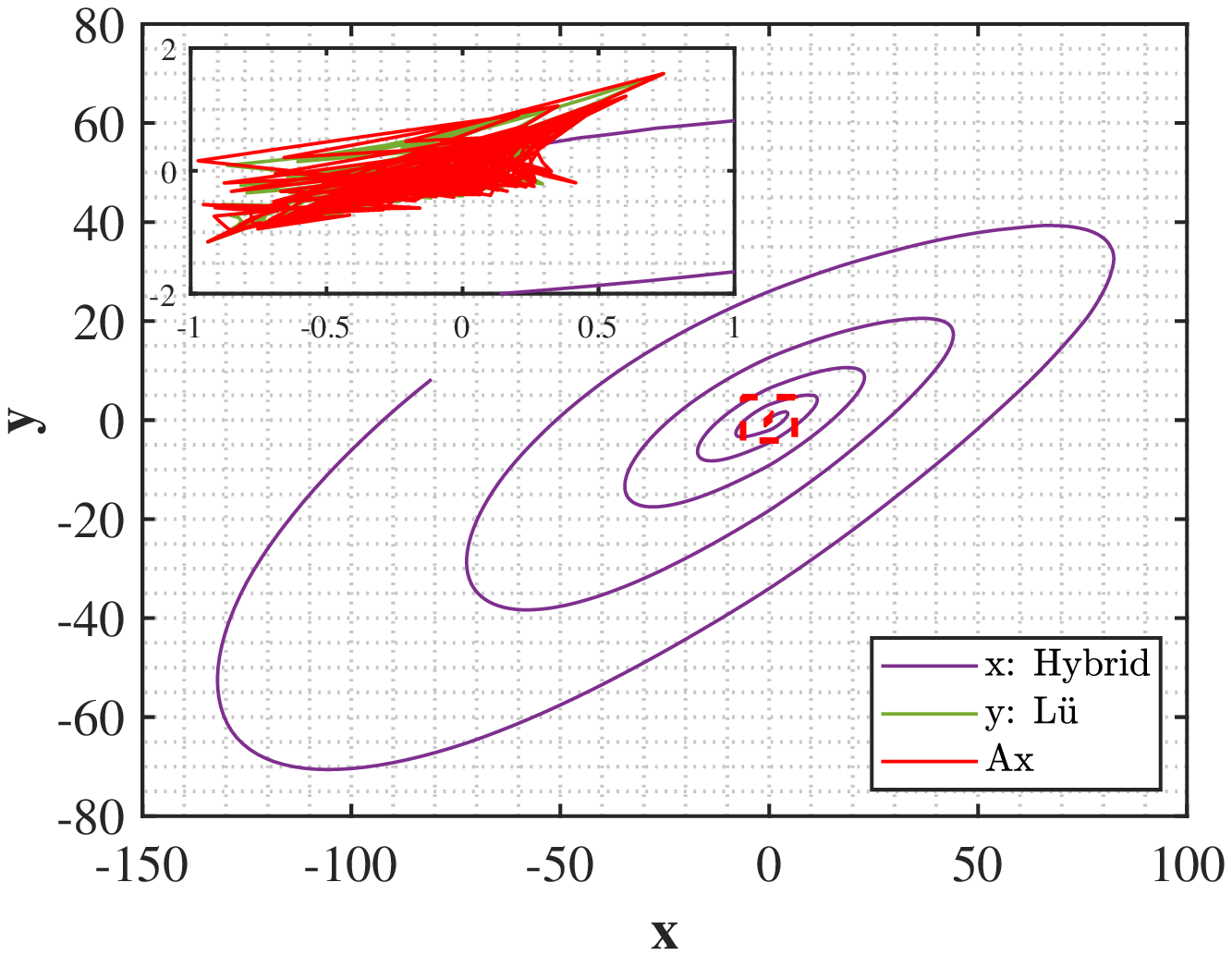}
\hfil
\includegraphics[scale=.28]{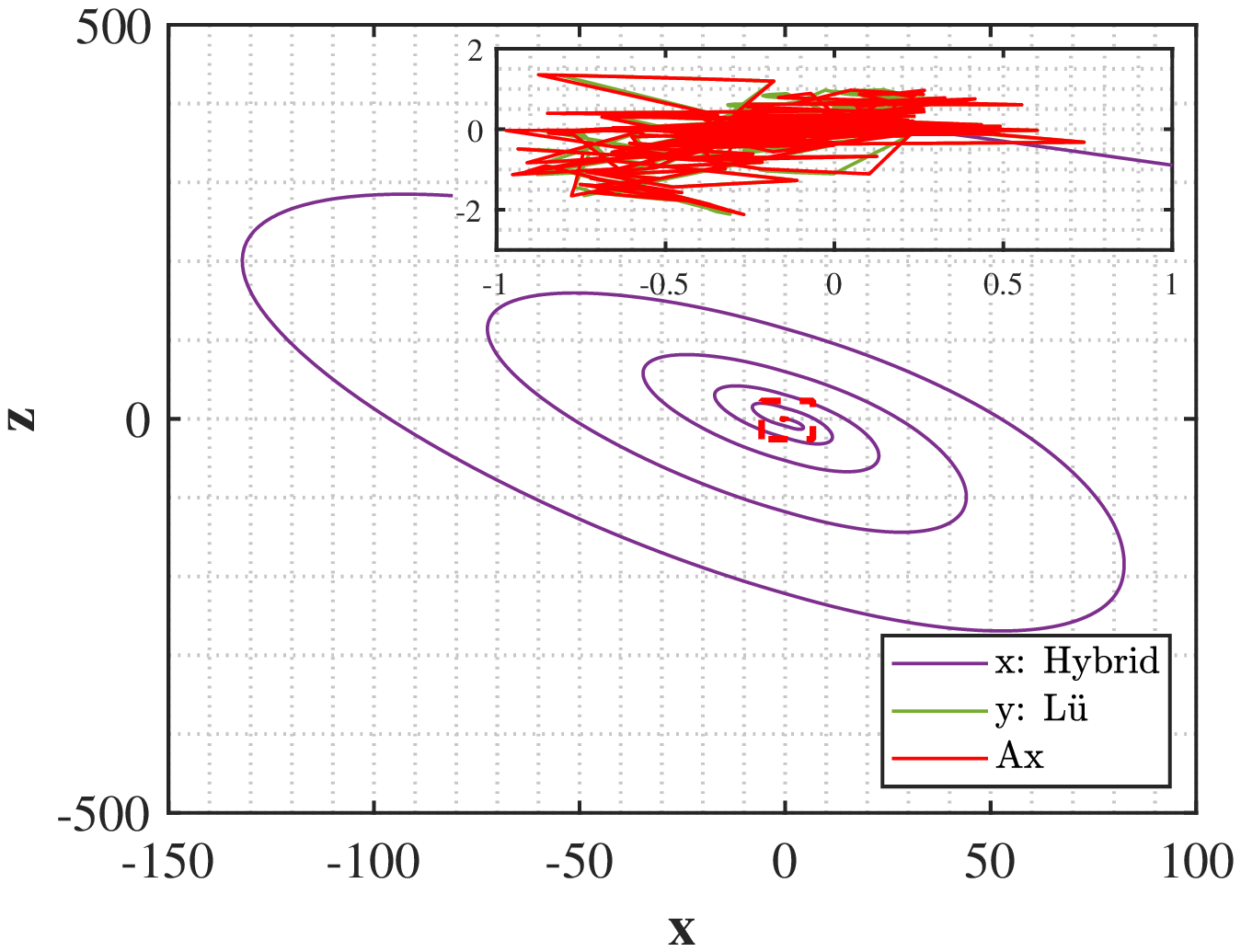}
\hfil
\includegraphics[scale=.28]{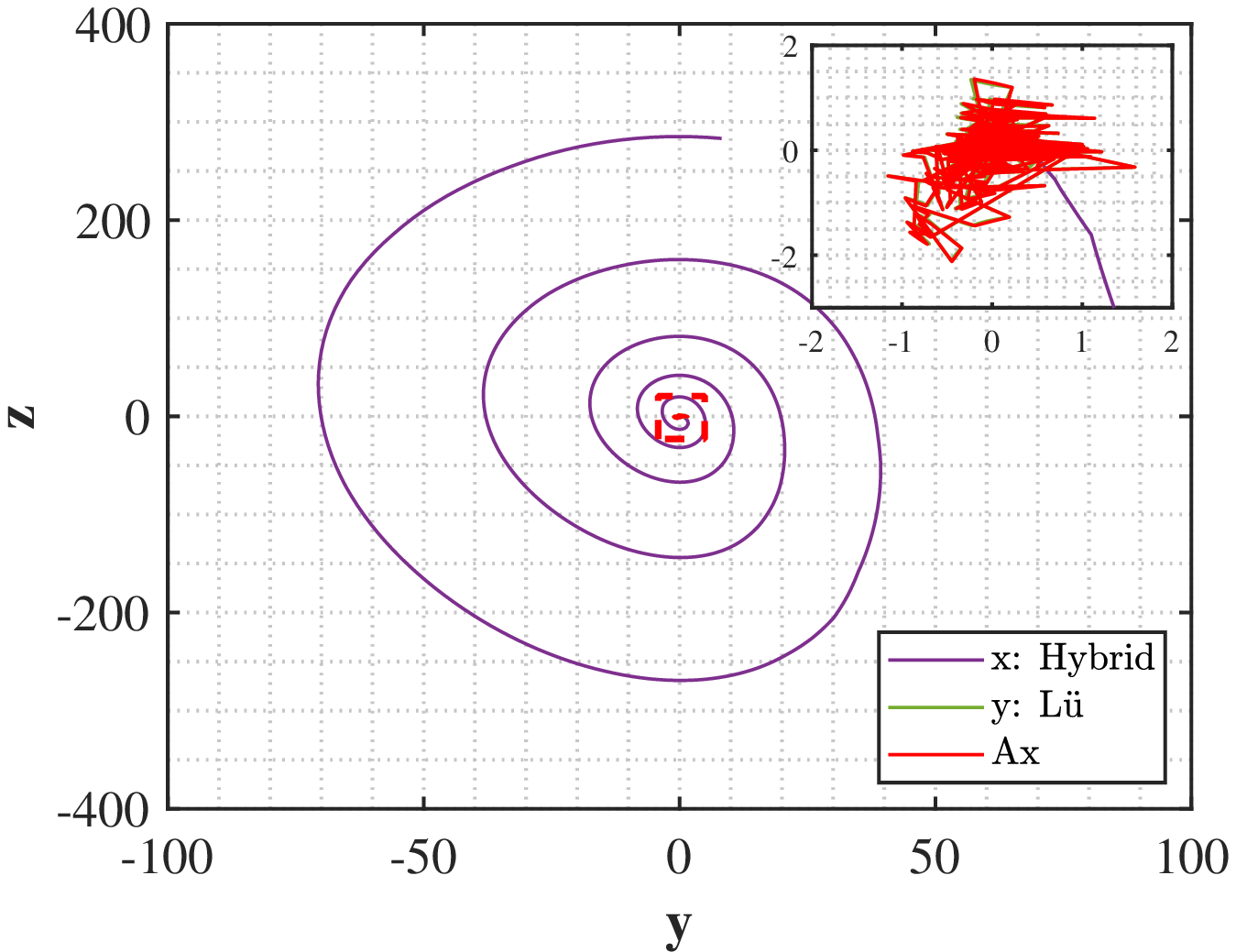}
\label{Fig.18c}}
\hfil
\subfloat[$\rm u$=-12.]{\includegraphics[scale=.28]{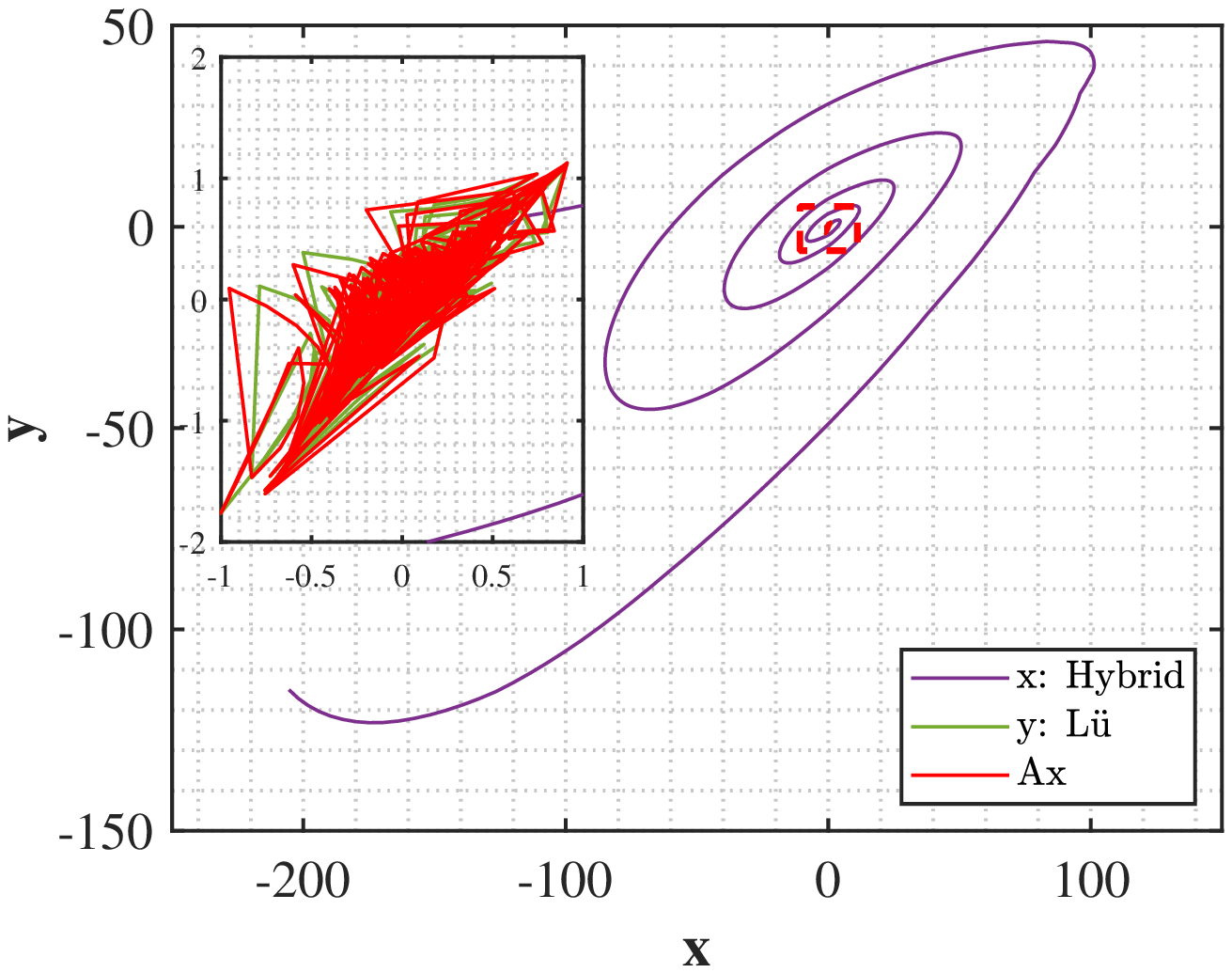}
\hfil
\includegraphics[scale=.28]{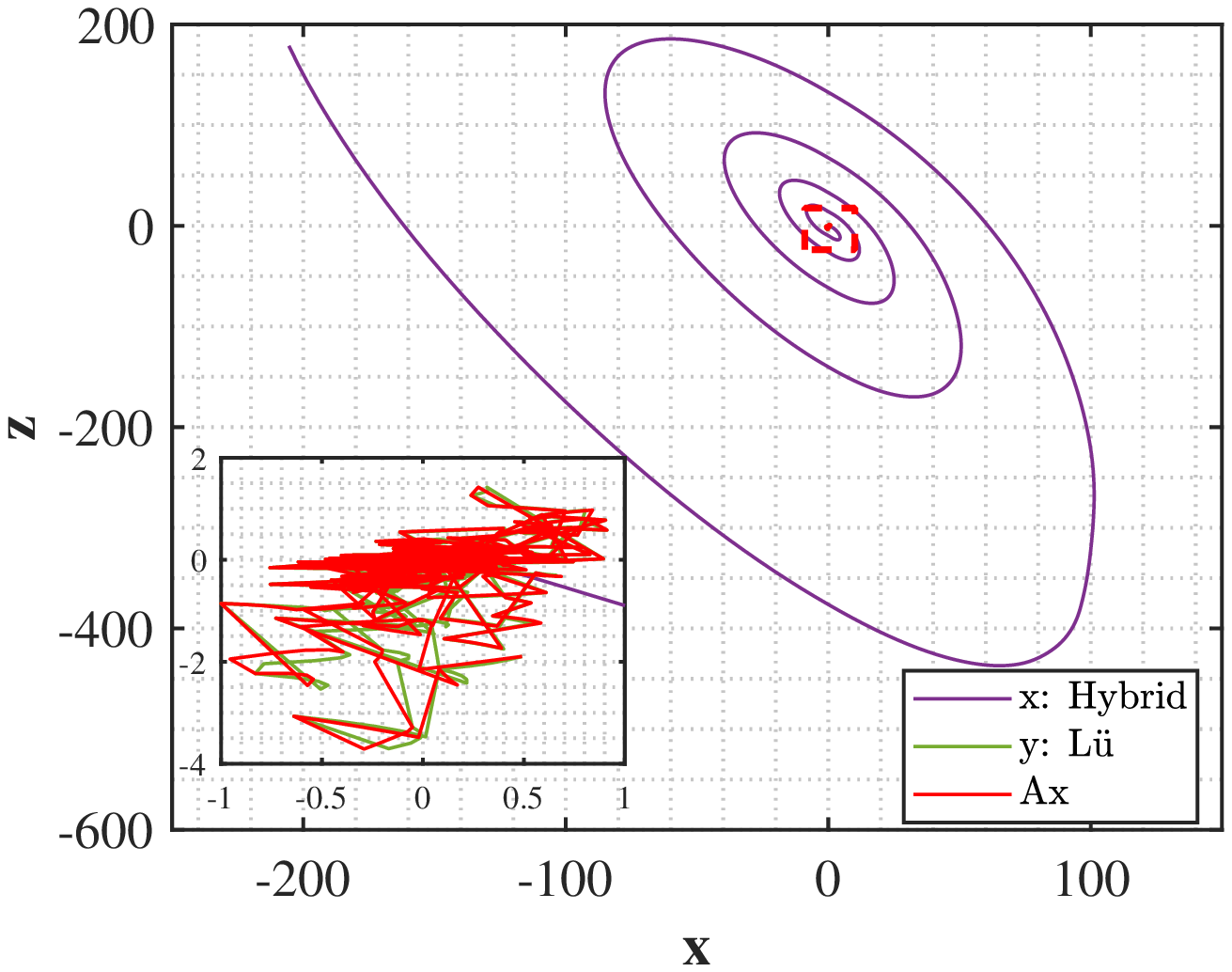}
\hfil
\includegraphics[scale=.28]{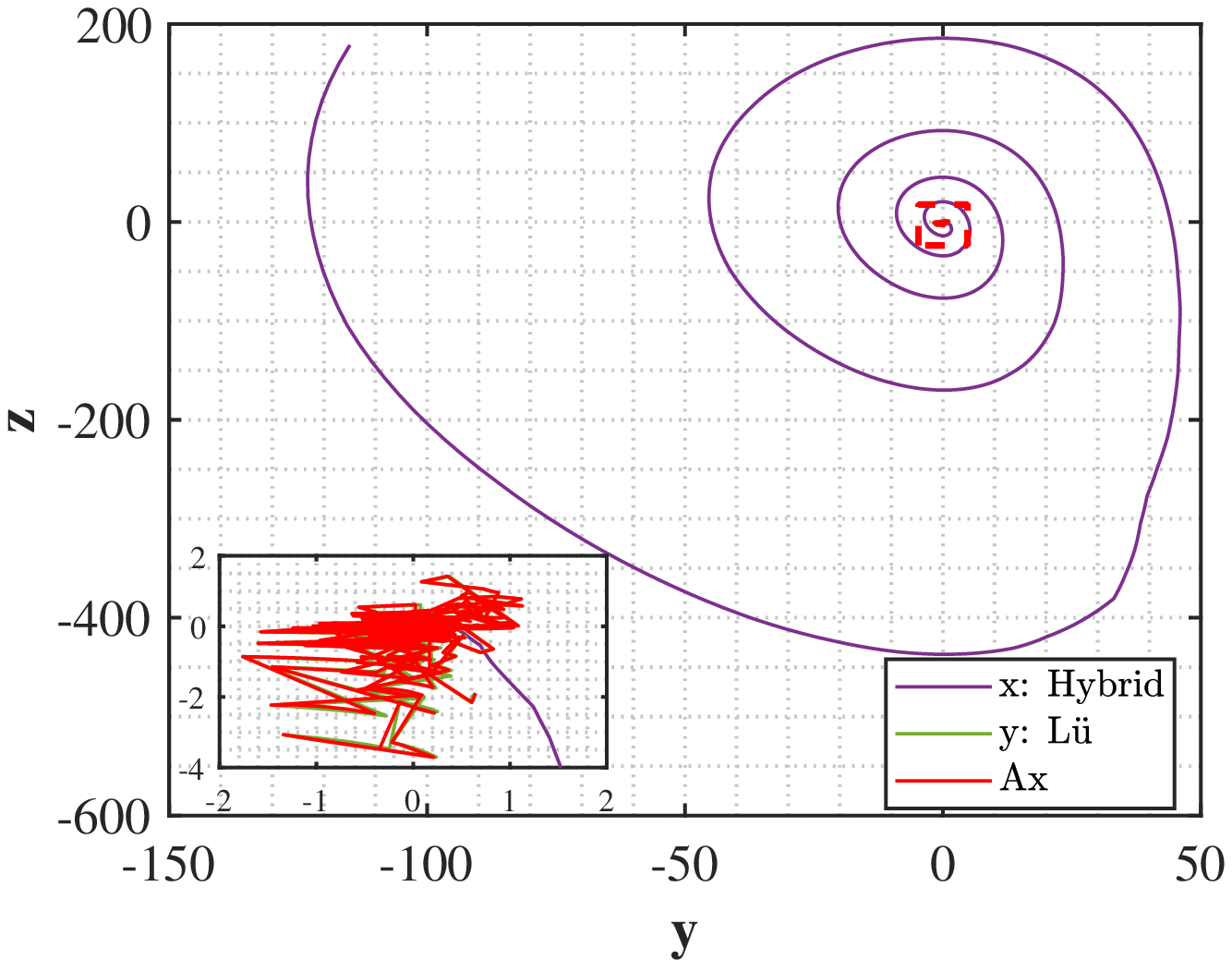}
\label{Fig.18d}}
\caption{Two dimensional plans of Example~4.5.}
\label{Fig.18}
\end{figure}

Let $\{x_k\}$ and $\{y_k\}$ be the numerical solutions obtained from hybrid system (\ref{41})  and L$\rm\ddot{u}$ attractor for 1000 time steps. Breaking the problems into 200 simple subproblems with 5 steps for each, the approximate solutions of parameter $\lambda$ and optimal similarity transformation matrix are found by (\ref{37}) and (\ref{38}).

\begin{figure}[htbp]
\centering
\subfloat[$\rm u$=-1.]{\includegraphics[scale=.47]{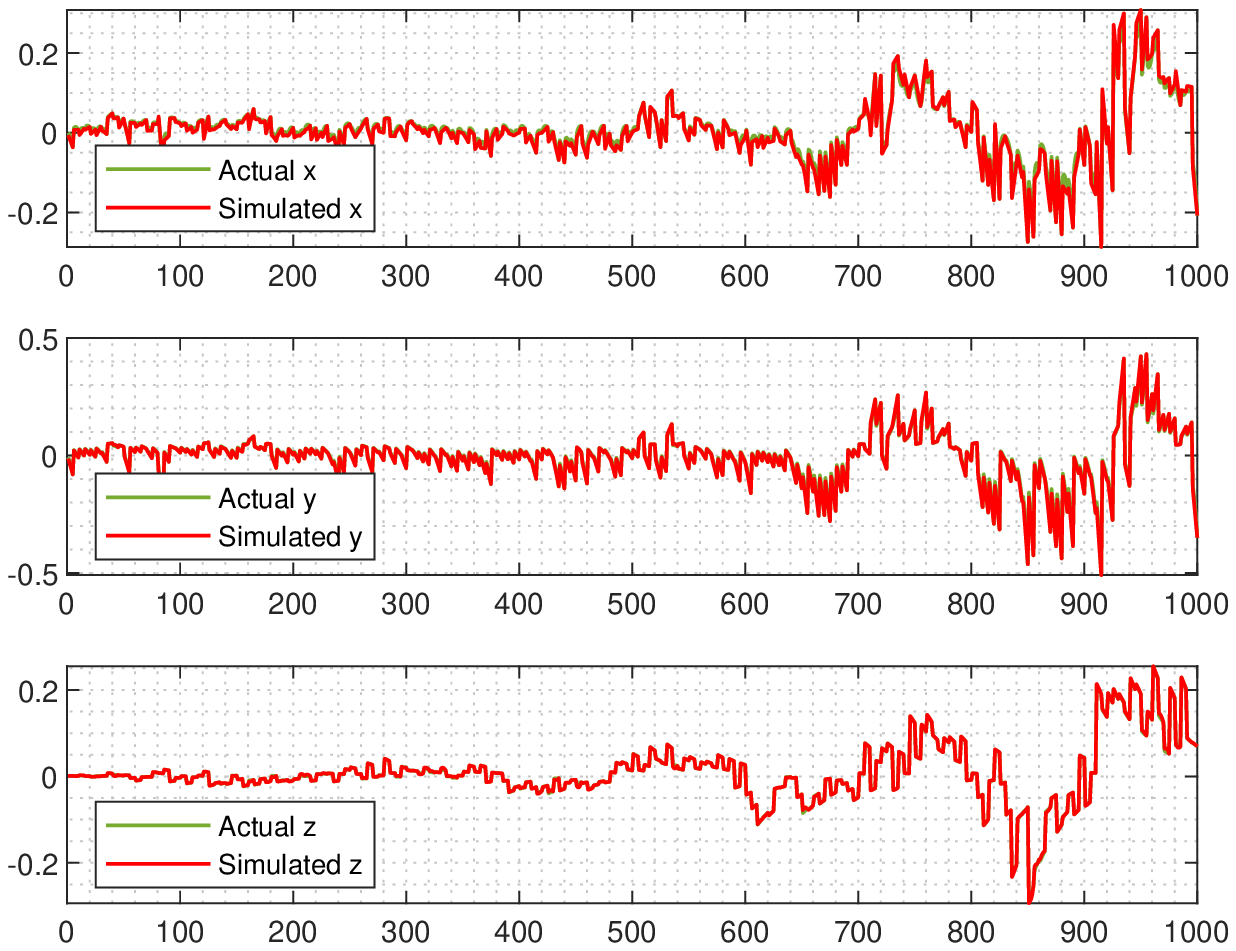}
\label{Fig.19a}}
\hfil
\subfloat[$\rm u$=8.]{\includegraphics[scale=.47]{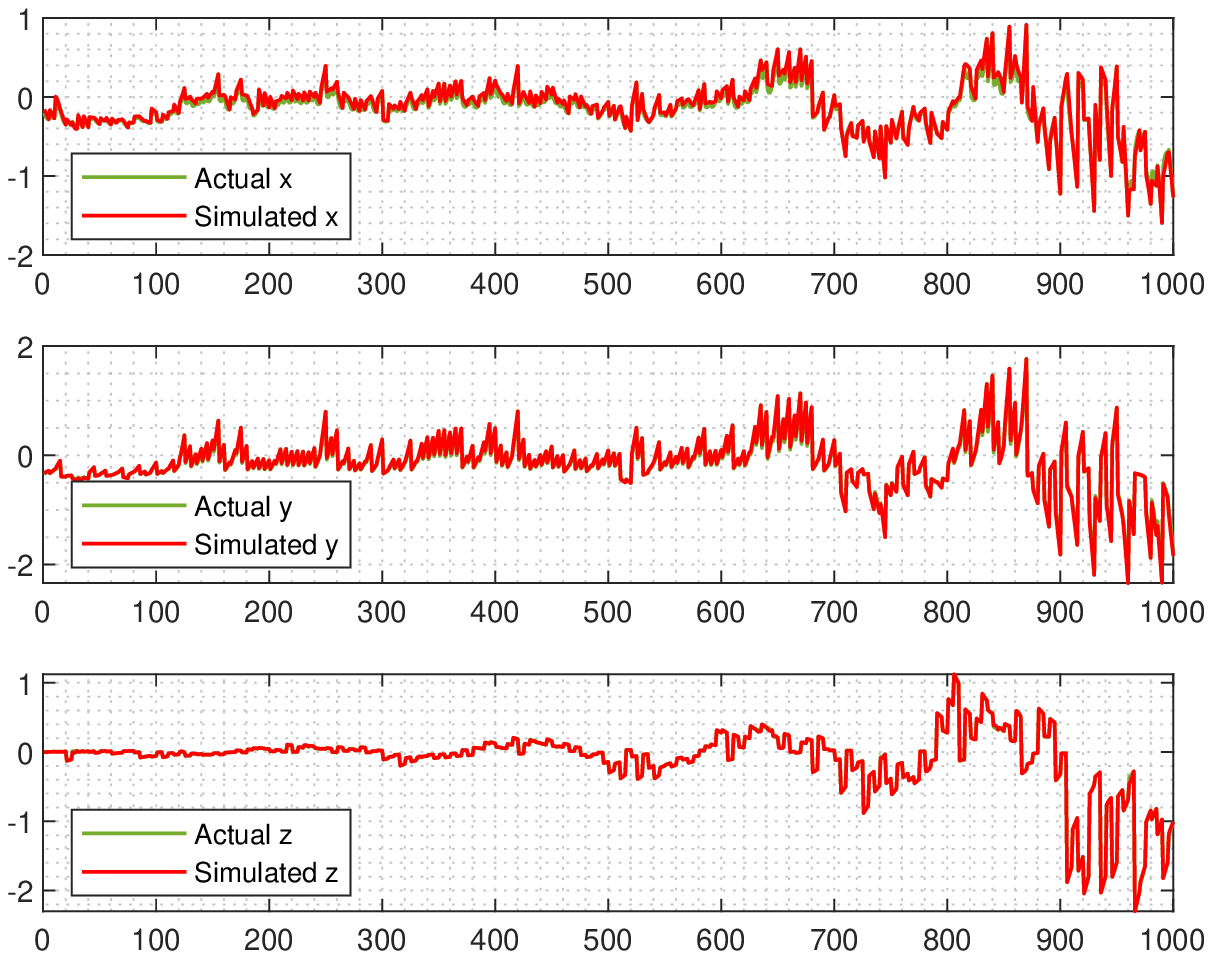}
\label{Fig.19b}}
\hfil
\subfloat[$\rm u$=12.]{\includegraphics[scale=.47]{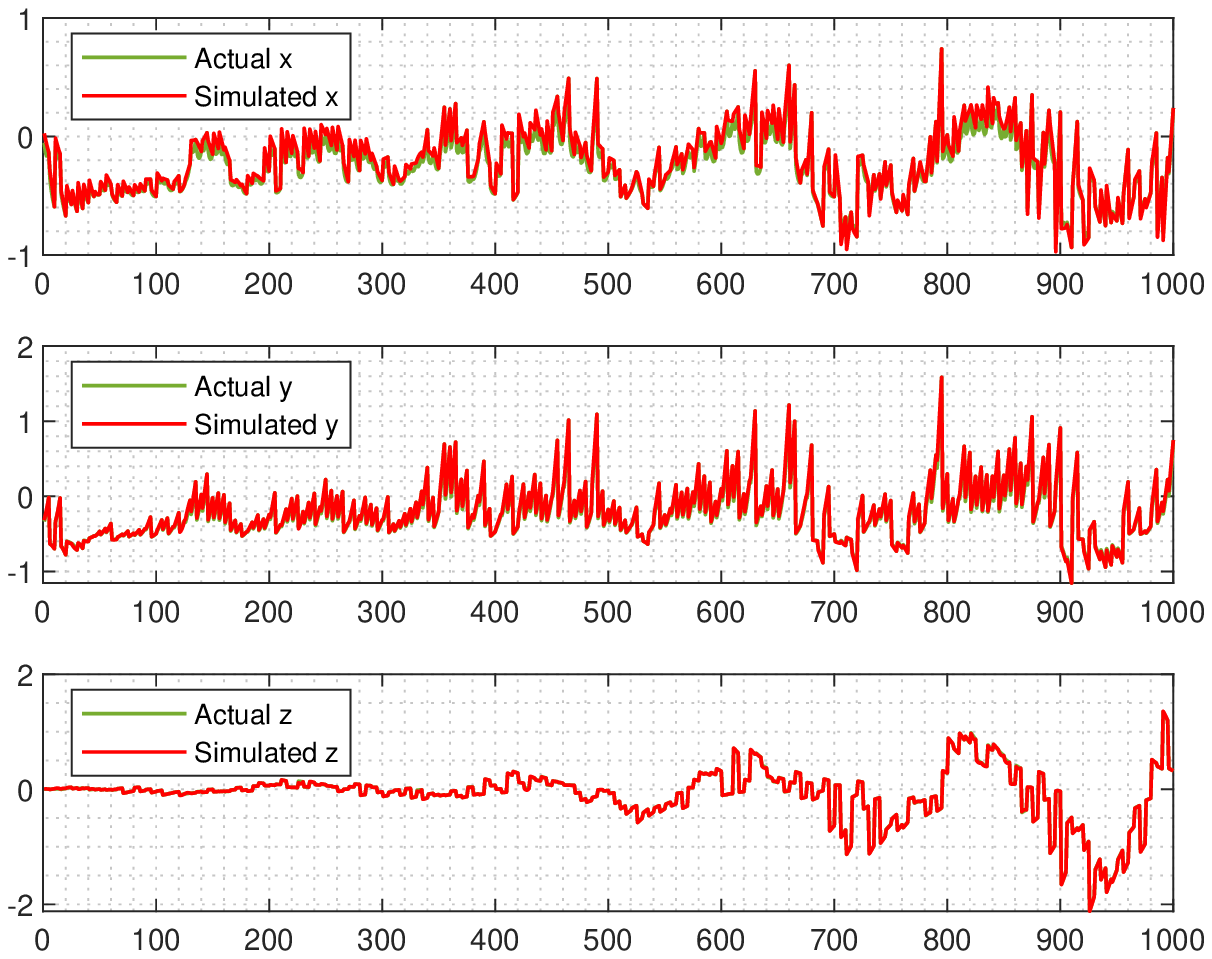}
\label{Fig.19c}}
\hfil
\subfloat[$\rm u$=-12.]{\includegraphics[scale=.47]{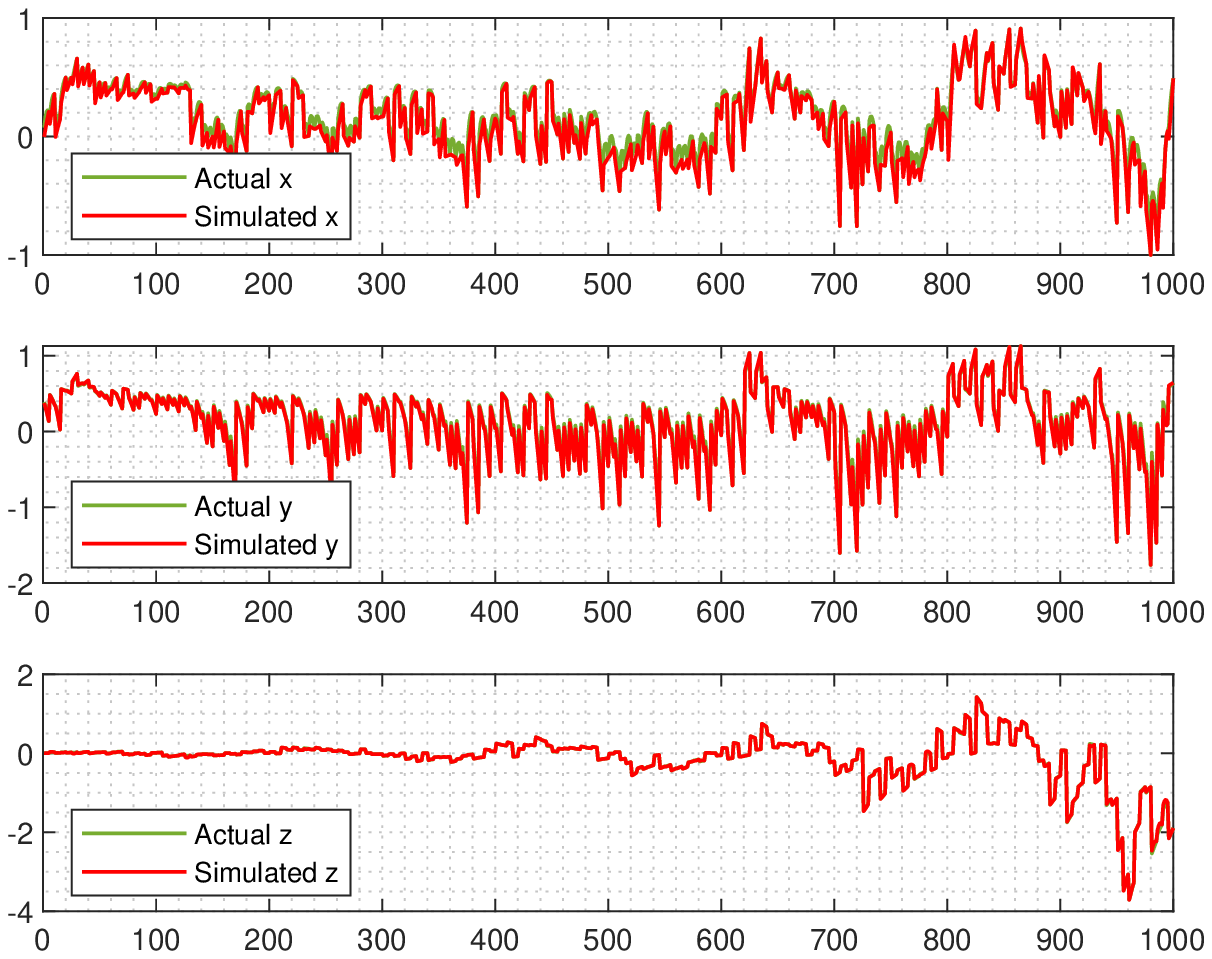}
\label{Fig.19d}}
\caption{Comparison of actual and simulated sequences of Example~4.5.}
\label{Fig.19}
\end{figure}

We show the numerical performance in Fig.\ref{Fig.18} for the four different values of $u$, respectively. For the purpose of demonstrating the optimal principle simulated effect more intuitively, Fig. \ref{Fig.19} illustrates the orbit of L$\rm\ddot{u}$ system (actual) and that of hybrid Lorenz-Chua attractor acted by the optimal similarity transformation matrix (simulated). In spite of L$\rm\ddot{u}$ system exhibits various dynamical behaviors due to varying parameter $u$, the two sequences almost completely coincide, showing the availability and universality of the proposed optimal principle.

\begin{figure}[htbp]
\centering
\includegraphics[scale=.6]{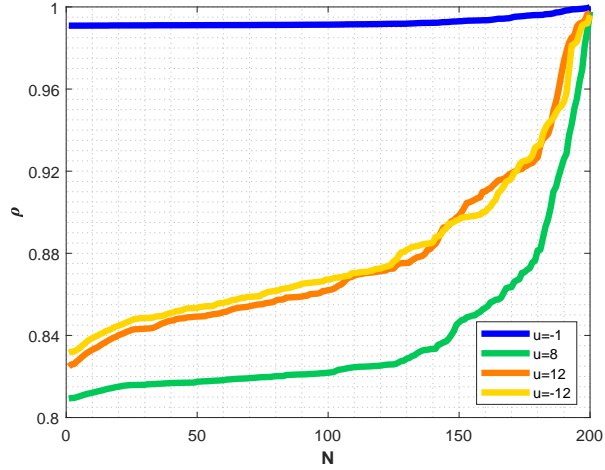}
\caption{Change in similarity degree of Example~4.5.}
\label{Fig.20}
\end{figure}

For the sake of completeness, the change in similarity degree of four cases of Example 4.5 are shown in Fig. \ref{Fig.20}, which also fulfill Bellman's principle of optimality.

Chaotic systems, are generally characterized by complex behavior and rapidly changing solutions, whose orbits become quite similar taking advantage of the general optimal principle presented in this paper.

\section{Conclusions}
\label{sec5}
In scientific research, capturing certain underlying similarity between two complex physical processes is one of the most intensively essential problems. The critical challenge for finding similarity of orbits between dynamical systems arises when facing the high sensitivity to small perturbations with respect to initial states of chaotic systems. Main contribution, in addition to proposing some concepts described what extent the orbits between two markedly different systems are similar, is the general optimal principle built up from the viewpoint of dynamical systems together with optimization theory. This optimal principle is applied to various well-known chaotic attractors and some kind of hybrid chaotic system that is constructed on the basis of homotopy idea, yielding encouraging numerical simulation results surprisingly.

Specifically, attention is paid to find similarity of orbits between dynamical systems with complex behavior, mathematically. As necessary foundations for this paper, the definitions of similarity transformation matrix and similarity degree are introduced. We present a general optimal principle based on optimality condition and variational method, finding some underlying similarity between orbits of two dynamical systems. The numerical simulations concern with both Pontryagin's maximum principle and Bellman's dynamic programming formulated in terms of the optimal principle for similarity of orbits between various chaotic systems. The orbits differed markedly become remarkably similar under action of the optimal similarity transformation matrix, and the value of similarity degree also supports this, implying significance of the optimal principle we proposed in this paper.

\section*{Acknowledgments}
This work was supported by National Basic Research Program of China (2013CB834100), National Natural Science Foundation of China (11571065, 11171132, 12071175), Project of Science and Technology Development of Jilin Province (2017C028-1, 20190201302JC), and Natural Science Foundation of Jilin Province (20200201253JC).


\end{document}